\documentclass[11pt]{amsart}
\usepackage{amssymb,amsmath,txfonts,mathrsfs,mathtools,titletoc, tikz, stmaryrd}
\usepackage[alphabetic]{amsrefs}
\usepackage{float}
\usepackage[colorlinks, linkcolor=blue,anchorcolor=Periwinkle,
    citecolor=red,urlcolor=Emerald]{hyperref}
\usepackage{xcolor}
\usepackage{enumitem}
\usepackage{geometry,array} 
\usepackage{graphicx}
\usepackage{subfigure}
\usepackage{bookmark}
\usepackage{tikz}\usetikzlibrary{matrix}
\usepackage{url}
\usepackage{braids}

\usetikzlibrary{decorations.markings}
\tikzset{->-/.style={decoration={  markings,  mark=at position #1 with
    {\arrow{>}}},postaction={decorate}}}
\tikzset{-<-/.style={decoration={  markings,  mark=at position #1 with
    {\arrow{<}}},postaction={decorate}}}
\usepackage[all]{xy}
\usetikzlibrary{arrows}
\usetikzlibrary{decorations.pathreplacing,decorations.pathmorphing,shadings,fadings,calc}
\usepackage{extarrows}
\usepackage{tikz-3dplot}
\usepackage{appendix}

\usetikzlibrary{matrix}
\usepackage[all]{xy}
\usetikzlibrary{arrows,shapes}
\usetikzlibrary{calc}
\tikzset{->-/.style={decoration={  markings,  mark=at position #1 with
    {\arrow{>}}},postaction={decorate}}}
\tikzset{-<-/.style={decoration={  markings,  mark=at position #1 with
    {\arrow{<}}},postaction={decorate}}}
\usepackage[]{hyperref}
\numberwithin{equation}{section}

\usetikzlibrary{decorations.markings}
\tikzset{->-/.style={decoration={  markings,  mark=at position #1 with
    {\arrow{>}}},postaction={decorate}}}
\tikzset{-<-/.style={decoration={  markings,  mark=at position #1 with
    {\arrow{<}}},postaction={decorate}}}

\usepackage{lscape}
\usepackage{CJK}

\newtheorem{theorem}{Theorem}[section]
\newtheorem{prop}[theorem]{Proposition}
\newtheorem{lemma}[theorem]{Lemma}
\newtheorem{remark}[theorem]{Remark}
\newtheorem{claim}[theorem]{Claim}

\newtheorem{cor}[theorem]{Corollary}

\numberwithin{equation}{section}

\def\QEDopen{{\setlength{\fboxsep}{0pt}\setlength{\fboxrule}{0.2pt}\fbox{\rule[0pt]{0pt}{1.3ex}\rule[0pt]{1.3ex}{0pt}}}}
\def\QED{\hfill\QEDopen} 
\def\proof{\noindent{\it Proof}: }

\newcommand{\R}{\mathbb{R}}

\newcommand{\Z}{\mathbb{Z}}

\newcommand{\sph}{\mathbb{S}}

\newcommand{\id}{\text{\rm id}}

\begin{document}
\def\nn{node{$\bullet$}}
\def\ww{node{$\circ$}}

\title{Blow Up of Compact Mean Curvature Flow Solutions with Bounded Mean Curvature}
\author[1]{Zichang Liu}
\address{Department of Mathematical Sciences\\Tsinghua University, Beijing\\P. R. China, 100084}
\email{liu-zc19@mails.tsinghua.edu.cn}
\date{\today}

\begin{abstract}
In 1994, Vel\'{a}zquez constructed a countable family of complete hypersurfaces flowing in $\R^{2N}$ $(N\geq 4)$ by mean curvature, each of which develops a type II singularity at the origin in finite time. Later Guo and Sesum showed that for a non-empty subset of Vel\'{a}zquez's solutions, the mean curvature blows up near the origin, at a rate smaller than that of the second fundamental form; recently Stolarski proved another subset of these solutions has bounded mean curvature up to the singular time. In this paper, we follow their arguments to construct compact mean curvature flow solutions in $\R^n$ $(n\geq 8)$ with bounded mean curvature. 
\end{abstract}
\thanks{}

\maketitle

\titlecontents{section}[0em]{}{\hspace{.5em}}{}{\titlerule*[1pc]{.}\contentspage}
\titlecontents{subsection}[1.5em]{}{\hspace{.5em}}{}{\titlerule*[1pc]{.}\contentspage}
\tableofcontents

\section{Introduction}
\label{in}
Let $M^{n-1}$ be a smooth manifold, $F:M^{n-1}\times [0,T)\rightarrow \R^n$ be a one-parameter family of immersions depending smoothly on both space and time. We say $F$ satisfies the \textbf{mean curvature flow} equation, if
\begin{equation}
    \partial_t F(x,t)^\perp=H(x,t),\ \ x\in M^{n-1},\ t\in [0,T), \label{in1}
\end{equation}
where $\partial_t=\frac{\partial}{\partial t}$, $\perp$ denotes the component of a vector perpendicular to the image of $F(\cdot,t)$, and $H(\cdot,t)$ is the mean curvature vector of the immersion $F(\cdot,t)$. Usually, we regard the image of $F(\cdot,t)$ as an immersed submanifold of $\R^n$, denoted by $\Sigma(t)$. From this viewpoint, the family of hypersurfaces $\{\Sigma(t)\}_{0\leq t<T}$ is considered to move by its mean curvature, regardless of the way it is parameterized.\\
According to \cite{hui}, when $M^{n-1}$ is compact, then for any smooth immersion $F_0=F(\cdot,0):M^{n-1}\rightarrow \R^n$ as an initial value, equation \eqref{in1} has a unique solution for a short time. A natural question is: will the solution exist permanently, or develop a singularity (blow up) in finite time? How can we characterize the first blow-up time of \eqref{in1}? Actually, when $M^{n-1}$ is compact, \eqref{in1} will always develop a finite-time singularity, and \cite{hui} tells us that the norm of the second fundamental form must become unbounded at the first singular time. In \cite{ac}, Cooper tells us that the inner product of the mean curvature vector and the second fundamental form, $H\cdot A$, must blow up as the flow becomes singular. One may naturally ask that ``does the mean curvature necessarily blow up at the first singular time", which was presented as an open problem in \cite{man}, Open Problem 2.4.10. At least for small $n$, people believe this is true, and in fact Li and Wang gave an affirmative answer in \cite{lw} for compact embedded surfaces flowing in $\R^3$. For 
$4\leq n\leq 7$, the answer remains unknown. However, in higher dimensions, things become completely different. Vel\'{a}zquez \cite{vel} constructed  a series of complete, embedded hypersurfaces flowing by mean curvature in $\R^n$ for each even number $n\geq 8$, such that they ``approach" the Simons' cone at the first singular time:
\begin{theorem} (\cite{vel})
    \label{thmin2}
    Let $N\geq 4$, $l\geq 2$ be integers. For $t_0<0$, $|t_0|\ll 1$, there exists a family of $O(N)\times O(N)$-invariant mean curvature flow solutions $\{\Sigma_{l}^{2N-1}(t)\}_{t_0\leq t<0}$ in $\R^{2N}$ such that
    \begin{enumerate}
        \item $\{\Sigma_{l}(t)\}_{t_0\leq t<0}$ develops a type II singularity at $0\in \R^{2N}$ in the sense that there exists $\sigma_l=\sigma_l(N)>0$ s.t. $A_{\Sigma_{l}(t)}$ (the second fundamental form of $\Sigma_{l}(t)$) satisfies
        $$\limsup_{t\nearrow 0} \sup_{\Sigma_{l}(t)\cap B(0,\sqrt{-t})}(-t)^{\frac{1}{2}+\sigma_l}|A_{\Sigma_{l}(t)}|>0.$$
        \item The type I rescaled hypersurfaces $(-t)^{-\frac{1}{2}}\Sigma_{l}(t)$ converges in the $C^2$ sense to Simons' cone on any compact subset of $\R^{2N}-\{0\}$ as $t\nearrow 0$.
        \item The type II rescaled hypersurfaces $(-t)^{-\frac{1}{2}-\sigma_l}\Sigma_{l}(t)$ converges uniformly to a minimal hypersurface tangent to Simons' cone at infinity, on any compact subset of $\R^{2N}$ as $t\nearrow 0$.
    \end{enumerate}
\end{theorem}
Guo and Sesum \cite{gs} performed more detailed analyses on Vel\'{a}zquez's solutions; in particular they improved the convergence results (2) (3) in Theorem \ref{thmin2} to $C^\infty$, and used those estimates to show:
\begin{theorem} (\cite{gs})
    \label{thmin3}
    For the mean curvature flows $\{\Sigma_{l}(t)\}_{t_0\leq t<0}$ in Theorem \ref{thmin2}, if $N\geq 5$ and $l=2$, then
    $$\limsup_{t\nearrow 0} \sup_{\Sigma_{2}(t)}(-t)^{\frac{1}{2}+\sigma_2}|A_{\Sigma_{2}(t)}|<+\infty.$$
    In addition, $H_{\Sigma_{2}(t)}$, the mean curvature of $\Sigma_{2}(t)$, blows up as $t\nearrow 0$ at a rate smaller than that of the second fundamental form. More precisely,
    $$\limsup_{t\nearrow 0} \sup_{\Sigma_{2}(t)\cap B(0,C(N)(-t)^{\frac{1}{2}+\sigma_2})}(-t)^{\frac{1}{2}-\sigma_2}|H_{\Sigma_{2}(t)}|>0,$$
    and there exists $0<\tilde{\sigma}<\sigma_2$ s.t.
    $$\limsup_{t\nearrow 0} \sup_{\Sigma_{2}(t)}(-t)^{\frac{1}{2}+\tilde{\sigma}}|H_{\Sigma_{2}(t)}|<+\infty.$$
\end{theorem}
Stolarski \cite{sto} showed that another subset of Vel\'{a}zquez's solutions has uniformly bounded mean curvature, using a blow-up method:
\begin{theorem} (\cite{sto})
    \label{thmin4}
    For the mean curvature flows $\{\Sigma_{l}(t)\}_{t_0\leq t<0}$ in Theorem \ref{thmin2}, if $l\geq 4$ is an even number, then
    $$\sup _{t\in [t_0,0)} \sup_{\Sigma_{l}(t)}|H_{\Sigma_{l}(t)}|<+\infty.$$
\end{theorem}

In \cite{vel}, Simons' cone and the minimal hypersurfaces tangent to it are important models for Vel\'{a}zquez's construction, but they only exist in even dimensional Euclidean spaces. To cover general cases, we use the so-called ``Lawson's cones" and related minimal hypersurfaces as models instead, which were thoroughly discussed in \cite{dav}. The author showed the ``Lawson's cones" are globally area minimizing in certain dimensions, by constructing minimal hypersurfaces tangent to them at infinity, and regarding the volume form of these hypersurfaces as a calibration form. In our paper, we obtain a compact mean curvature flow solution with an asymptotic behavior similar to \cite{vel} by writing down the initial compact hypersurface explicitly. \cite{gs} provided a simpler proof of the $C^0$ estimate, but the argument does not work when $N=4$; instead we turn to the original proof in \cite{vel}, which is inspired by \cite{hv}. Here we present a complete proof. Finally, it should be noted that the condition ``$l$ is even" in Theorem \ref{thmin4} is superfluous, if one performs a careful observation to the estimates in \cite{vel}. In brief, we improve and generalize their results as the following theorem:
\begin{theorem}
    \label{thmin1}
    For each dimension $n\geq 8$, there exists a smooth, compact, embedded mean curvature flow $\{\Sigma^{n-1}(t)\}_{t\in [t_0,0)}\subset \R^n$ such that
    $$\limsup _{t\nearrow 0}\sup_{x\in \Sigma(t)} |A_{\Sigma(t)}(x)|=+\infty,\ \ \sup _{t\in [t_0,0)}\sup_{x\in \Sigma(t)} |H_{\Sigma(t)}(x)|<+\infty.$$
\end{theorem}
\noindent{\bf Acknowledgements:} I would like to thank J.J.L. Vel\'{a}zquez, Siao-Hao Guo, Natasa Sesum, and Maxwell Stolarski for their insights into this problem. I also thank my advisor, Xiaoli Han, for helpful instructions.

\section{Minimal Hypersurfaces Tangent to Lawson's Cones at Infinity}
\label{mi}
We first take a glance at the models we shall use, i.e. Lawson`s cones and related minimal hypersurfaces. From now on, we fix $n\in \Z$, $n\geq 8$, and two integers $p,\ q\geq 2$ with $p+q=n$. If $n=8$, we require $p,\ q\geq 3$. Any point in $\R^n$ is denoted by $z=(x,y)$, $x\in \R^p$, $y\in \R^q$.\\
Let $O(p)$ be the group of isometries of $\R^p$ which fix the origin. For $(g,h)\in O(p)\times O(q)$, we define the group action of $O(p)\times O(q)$ on $\R^n$ by $(g,h)(x,y)=(g(x),h(y))$. We shall restrict our discussion to hypersurfaces in $\R^n$ which are invariant under the action of $O(p)\times O(q)$.\\
Let $\gamma$: $x^1=\xi(t)$, $y^1=\eta(t)$ be a smooth curve lying in the first quadrant of the $x^1 y^1$-plane, i.e. the points with $x^1$, $y^1\geq 0$, $x^2=...=x^p=y^2=...=y^q=0$, where $t$ belongs to some real interval, say $J$, and $\xi'^2+\eta'^2\neq 0$. Then the set $\Sigma=\{(\xi(t)\nu,\eta(t)\omega)|\ t\in J,\ \nu\in \sph^{p-1},\ \omega\in \sph^{q-1}\}$ is a smooth hypersurface in $\R^n=\R^p\times \R^q$, possibly singular at the points where $x=0$ or $y=0$. We say $\Sigma$ is the hypersurface \textbf{generated by} $\gamma$ (by the action of $O(p)\times O(q)$), and $\gamma$ is the \textbf{profile curve} of $\Sigma$.\\
Write $\mu=\sqrt{\frac{q-1}{p-1}}$, and let $C_{p,q}$ be the Lawson's cone, i.e. the hypersurface generated by the ray
\begin{equation}
    l_{p,q}=\{(x^1,y^1)|\ y^1=\mu x^1,\ x^1\geq 0\}. \label{mi25}
\end{equation}
According to \cite{dav} and \cite{vel}, there is a smooth minimal hypersurface 
$$\mathcal{M}=\{(r\nu,\hat{\psi}(r)\omega)|\ r\geq 0,\ \nu\in \sph^{p-1},\ \omega\in \sph^{q-1}\}$$
tangent to $C_{p,q}$ at infinity, with $\hat{\psi}(r)$ satisfying the equation
\begin{equation}
    \frac{\hat{\psi}''}{1+\hat{\psi}'^2}+(p-1)\frac{\hat{\psi}'}{r}-(q-1)\frac{1}{\hat{\psi}}=0, \label{mi1}
\end{equation}
and
\begin{equation}
    \begin{cases}
    \hat{\psi}''>0,\ \hat{\psi}(0)>0,\ \hat{\psi}'(0)=0,\ \hat{\psi}(r)>\mu r,\\
    \lim_{r\rightarrow +\infty} \frac{\hat{\psi}(r)-\mu r}{r^\alpha}=(1+\mu^2)^{\frac{\alpha+1}{2}},\\
    \lim_{r\rightarrow +\infty} \frac{\hat{\psi}'(r)-\mu}{r^{\alpha-1}}=\alpha(1+\mu^2)^{\frac{\alpha+1}{2}},
\end{cases}
\label{mi2}
\end{equation}
where 
$$\alpha=\frac{1}{2}(3-n+\sqrt{n^2-10 n+17})\in [-2,-1)$$
is a root of 
\begin{equation}
    x(x-1)+(n-2)(x+1)=0. \label{mi9}
\end{equation}
The profile curve
\begin{equation}
    \bar{\mathcal{M}}=\{(r,\hat{\psi}(r))|\ r\geq 0\} \label{mi3}
\end{equation}
of $\mathcal{M}$ can also be parameterized as a graph over $l_{p,q}$; more precisely,
$$\bar{\mathcal{M}}=\{\frac{r(1,\mu)}{\sqrt{1+\mu^2}}+\frac{\psi(r)(-\mu,1)}{\sqrt{1+\mu^2}}|\ r\geq \hat{\psi}(0)\frac{\mu}{\sqrt{1+\mu^2}}\}$$
\begin{equation}
    =\{(\frac{r-\mu\psi(r)}{\sqrt{1+\mu^2}},\frac{\mu r+\psi(r)}{\sqrt{1+\mu^2}})|\ r\geq \hat{\psi}(0)\frac{\mu}{\sqrt{1+\mu^2}}\}. \label{mi4}
\end{equation}
The function $\psi(r)$ satisfies 
\begin{equation}
    \frac{\psi''}{1+\psi'^2}+(p-1)\frac{\mu+\psi'}{r-\mu\psi}-(q-1)\frac{1-\mu\psi'}{\mu r+\psi}=0, \label{mi5}
\end{equation}
and
\begin{equation}
    \begin{cases}
    \psi''>0,\ \psi(\hat{\psi}(0)\frac{\mu}{\sqrt{1+\mu^2}})=\frac{\hat{\psi}(0)}{\sqrt{1+\mu^2}},\ \psi'(\hat{\psi}(0)\frac{\mu}{\sqrt{1+\mu^2}})=-\mu,\\
    \lim_{r\rightarrow +\infty} \frac{\psi(r)}{r^\alpha}=1,\\
    \lim_{r\rightarrow +\infty} \frac{\psi'(r)}{r^{\alpha-1}}=\alpha.
\end{cases}
\label{mi6}
\end{equation}
More generally, for each $k>0$, we can define
$$\mathcal{M}_k=k^{\frac{1}{1-\alpha}}\mathcal{M},$$
then $\mathcal{M}_k$ is also a smooth minimal hypersurface tangent to $C_{p,q}$ at infinity. Denote by $\bar{\mathcal{M}}_k$ the profile curve of $\mathcal{M}_k$, then $\bar{\mathcal{M}}_k=k^{\frac{1}{1-\alpha}}\bar{\mathcal{M}}$. Clearly, $\bar{\mathcal{M}}_k$ can be parameterized as \eqref{mi3} with $\hat{\psi}$ replaced by $\hat{\psi}_k$, and 
$$\hat{\psi}_k(r)=k^{\frac{1}{1-\alpha}}\hat{\psi}(k^{-\frac{1}{1-\alpha}}r).$$
$\hat{\psi}_k$ satisfies \eqref{mi1}, and \eqref{mi2} with the two limits replaced by $k(1+\mu^2)^{\frac{\alpha+1}{2}}$ and $k\alpha(1+\mu^2)^{\frac{\alpha+1}{2}}$ respectively. $\bar{\mathcal{M}}_k$ can also be parameterized as \eqref{mi4} with $\psi$ replaced by $\psi_k$, and 
$$\psi_k(r)=k^{\frac{1}{1-\alpha}}\psi(k^{-\frac{1}{1-\alpha}}r).$$
$\psi_k$ satisfies \eqref{mi5}, and \eqref{mi6} with $\hat{\psi}$ replaced by $\hat{\psi}_k$, and the two limits replaced by $k$ and $k\alpha$ respectively.\\
Moreover, as in \cite{gs}, we have
\begin{equation}
    \hat{\psi}_k(r)-r\hat{\psi}_k'(r)>0,\ \psi_k(r)-r\psi_k'(r)>0, \label{mi23}
\end{equation}
\begin{equation}
    \lim_{r\rightarrow +\infty} \frac{\hat{\psi}_k(r)-r\hat{\psi}_k'(r)}{r^\alpha}=k(1-\alpha)(1+\mu^2)^{\frac{\alpha+1}{2}},\ \lim_{r\rightarrow +\infty} \frac{\psi_k(r)-r\psi_k'(r)}{r^\alpha}=k(1-\alpha), \label{mi24}
\end{equation}
\begin{equation}
    \partial_k \hat{\psi}_k(r)=\frac{1}{k(1-\alpha)}(\hat{\psi}_k(r)-r\hat{\psi}_k'(r)) >0,\ \partial_k \psi_k(r)=\frac{1}{k(1-\alpha)}(\psi_k(r)-r\psi_k'(r))>0. \label{mi21}
\end{equation}
Also, we have the following higher order estimates:
\begin{equation}
    \lim_{r\rightarrow +\infty} \frac{\psi_k^{(m)}(r)}{r^{\alpha-m}}=k\alpha(\alpha-1)(\alpha-2)...(\alpha-m+1)\ (m\geq 2), \label{mi7}
\end{equation}
\begin{equation}
    |\psi_k^{(m)}(r)|\leq C(p,q,m)kr^{\alpha-m}\ (m\geq 0,\ r\geq \hat{\psi}_k(0)\frac{\mu}{\sqrt{1+\mu^2}}), \label{mi15}
\end{equation}
$$\lim_{r\rightarrow +\infty} \frac{\hat{\psi}_k^{(m)}(r)}{r^{\alpha-m}}=k\alpha(\alpha-1)(\alpha-2)...(\alpha-m+1)(1+\mu^2)^{\frac{\alpha+1}{2}}\ (m\geq 2),$$
\begin{equation}
    |(\hat{\psi}_k(r)-\mu r)^{(m)}|\leq C(p,q,m)k(k^{\frac{1}{1-\alpha}}+r)^{\alpha-m}\ (m\geq 0,\ r\geq 0). \label{mi22}
\end{equation}
Next, we derive an estimate of the difference between $\psi_k$ and its asymptotic function appeared in \eqref{mi7}, which was proved in \cite{gs}, Lemma 2.5 for Simons' cones. Things are slightly different for general Lawson's cones, however:
\begin{lemma}
    \label{lemmi1}
    \begin{equation}
        |\psi_k^{(m)}(r)-k\alpha(\alpha-1)(\alpha-2)...(\alpha-m+1)r^{\alpha-m}|\leq C(p,q,m)k^{\frac{1-\tilde{\alpha}}{1-\alpha}} r^{\tilde{\alpha}-m},\ (m\geq 0,\ r\geq \hat{\psi}_k(0)\frac{\mu}{\sqrt{1+\mu^2}}), \label{mi8}
    \end{equation}
    where
    \begin{equation}
        \tilde{\alpha}=\max\{2\alpha-1,\hat{\alpha}\}=\begin{cases}
        2\alpha-1\ (n\geq 9),\\
        \hat{\alpha}\ (n=8),
    \end{cases} \label{mi20}
    \end{equation}
    and $\hat{\alpha}=\frac{1}{2}(3-n-\sqrt{n^2-10 n+17})$ is another root of \eqref{mi9}.
\end{lemma}
\proof
We may assume $k=1$; the conclusion for general $k$ follows by a simple scaling argument.\\
In order to turn \eqref{mi5} into an autonomous equation, we set
\begin{equation}
    s=\ln r,\ \ W(s)=e^{-s}\psi(e^s)\ (s\geq \ln (\hat{\psi}_1(0)\frac{\mu}{\sqrt{1+\mu^2}})),\ \ Z=W',\label{mi10}
\end{equation}
and then
\begin{equation}
    \psi(r)=rW(\ln r),\ \ \psi'(r)=(W+W')(\ln r),\ \ \psi''(r)=\frac{1}{r}(W'+W'')(\ln r). \label{mi12}
\end{equation}
Now, \eqref{mi5} becomes
\begin{equation}
    \begin{cases}
    W'=Z,\\
    Z'=-Z-(n-2)(1+(W+Z)^2)(\frac{1+(\mu^{-1}-\mu)W}{(1-\mu W)(1+\mu^{-1}W)}(W+Z)+\frac{W}{(1-\mu W)(1+\mu^{-1}W)}).
    \end{cases}
    \label{mi11}
\end{equation}
Since $\psi$ and $\psi'$ are asymptotic to $r^\alpha$ and $\alpha r^{\alpha-1}$ respectively, as $r\rightarrow +\infty$, we know $W$ (resp. $Z$) is asymptotic to $W_*:=e^{(\alpha-1)s}$ (resp. $Z_*:=(\alpha-1)e^{(\alpha-1)s}$), as $s\rightarrow +\infty$, by making use of \eqref{mi10} and \eqref{mi12}. Clearly, $W_*$ and $Z_*$ satisfy the following equation
\begin{equation}
    \begin{cases}
    W_*'=Z_*,\\
    Z_*'=-Z_*-(n-2)(2W_*+Z_*).
    \end{cases}
    \label{mi13}
\end{equation}
Now subtract \eqref{mi13} from \eqref{mi11} to get
\begin{equation}
    \begin{cases}
    (W-W_*)'=Z-Z_*,\\
    (Z-Z_*)'=-(Z-Z_*)-(n-2)(2(W-W_*)+(Z-Z_*))+f(s),
    \end{cases}
    \label{mi14}
\end{equation}
where
$$f(s)=\frac{n-2}{(1-\mu W)(1+\mu^{-1}W)}((\mu^{-1}-\mu)(W^2-W(W+Z)^3)-((W+Z)^2+W^2)(2W+Z))(s),$$
with
$$f(s)=O(e^{2(\alpha-1)s})\ (s\rightarrow +\infty,\ \text{by}\ \eqref{mi15}).$$
Also, we have
$$W-W_*=o(e^{(\alpha-1)s}),\ Z-Z_*=o(e^{(\alpha-1)s})\ (s\rightarrow +\infty,\ \text{by}\ \eqref{mi6}).$$
Rewrite \eqref{mi14} further as
\begin{equation*}
    \binom{W-W_*}{Z-Z_*}'=A\binom{W-W_*}{Z-Z_*}+\binom{0}{f(s)},\ A=\left(
\begin{array}{cc}
    0 & 1 \\
    -2(n-2) & -(n-1)
\end{array}
\right).
\end{equation*}
Since $A$ has two distinct eigenvalues $\alpha-1$ and $\hat{\alpha}-1$ with $\hat{\alpha}-1<\alpha-1<0$, and
$$A=\left(
\begin{array}{cc}
    1 & 1 \\
    \alpha-1 & \hat{\alpha}-1
\end{array}
\right)
\left(
\begin{array}{cc}
    \alpha-1 & 0 \\
    0 & \hat{\alpha}-1
\end{array}
\right)
\left(
\begin{array}{cc}
    1 & 1 \\
    \alpha-1 & \hat{\alpha}-1
\end{array}
\right)^{-1},
$$
we set
$$\binom{U}{V}=\left(
\begin{array}{cc}
    1 & 1 \\
    \alpha-1 & \hat{\alpha}-1
\end{array}
\right)^{-1} \binom{W-W_*}{Z-Z_*},\ \binom{g(s)}{h(s)}=\left(
\begin{array}{cc}
    1 & 1 \\
    \alpha-1 & \hat{\alpha}-1
\end{array}
\right)^{-1} \binom{0}{f(s)},
$$
to get the equation which $U$, $V$ satisfy:
$$\begin{cases}
    U'=(\alpha-1) U+g(s),\\
    V'=(\hat{\alpha}-1) V+h(s).
\end{cases}$$
Again, we still have the estimates
$$g(s)=O(e^{2(\alpha-1)s}),\ h(s)=O(e^{2(\alpha-1)s})\ (s\rightarrow +\infty),$$
$$U=o(e^{(\alpha-1) s}),\ V=o(e^{(\alpha-1) s})\ (s\rightarrow +\infty),$$
and there exists a constant $C(p,q)>0$ s.t.
$$|g(s)|,\ |h(s)|\leq C(p,q)e^{2(\alpha-1)s}\ (s\geq \ln (\hat{\psi}_1(0)\frac{\mu}{\sqrt{1+\mu^2}})+1).$$
For $s_0>s\geq \ln (\hat{\psi}_1(0)\frac{\mu}{\sqrt{1+\mu^2}})+1$ ,
$$U(s)=e^{(\alpha-1) (s-s_0)}U(s_0)+\int_{s_0}^s e^{(\alpha-1)(s-t)}g(t)dt,$$
$$|U(s)|\leq e^{(\alpha-1) s}e^{-(\alpha-1) s_0}|U(s_0)|+C(p,q)e^{(\alpha-1) s}\int_s^{+\infty} e^{(\alpha-1) t}dt$$
\begin{equation}
    =e^{(\alpha-1) s} e^{-(\alpha-1) s_0}|U(s_0)|+C(p,q)e^{2(\alpha-1)s}. \label{mi16}
\end{equation}
Since $s_0>s$ is arbitrary, if we let $s_0\rightarrow +\infty$, the first term of the right side of \eqref{mi16} tends to zero, which implies
$$|U(s)|\leq C(p,q)e^{2(\alpha-1)s}.$$
Similarly, for $s>s_0=\ln (\hat{\psi}_1(0)\frac{\mu}{\sqrt{1+\mu^2}})+1$,
\begin{equation}
    V(s)=e^{(\hat{\alpha}-1) (s-s_0)}V(s_0)+\int_{s_0}^s e^{(\hat{\alpha}-1)(s-t)}h(t)dt. \label{mi17}
\end{equation}
Clearly, the first term of the right side of \eqref{mi17} is $O(e^{(\hat{\alpha}-1) s})$ as $s\rightarrow +\infty$, and the second term is bounded by
\begin{equation}
    C(p,q)e^{(\hat{\alpha}-1) s}\int_{s_0}^s e^{(2(\alpha-1)-(\hat{\alpha}-1))t}dt. \label{mi18}
\end{equation}
If $n=8$, then $2(\alpha-1)<\hat{\alpha}-1$, and \eqref{mi18} is bounded by 
$$C(p,q)e^{(\hat{\alpha}-1) s}\int_{s_0}^{+\infty} e^{(2(\alpha-1)-(\hat{\alpha}-1))t}dt=C(p,q)e^{(\hat{\alpha}-1) s}.$$
If $n\geq 9$, then $2(\alpha-1)>\hat{\alpha}-1$, and \eqref{mi18} is bounded by 
$$C(p,q)e^{(\hat{\alpha}-1) s}\int_{-\infty}^{s} e^{(2(\alpha-1)-(\hat{\alpha}-1))t}dt=C(p,q)e^{2(\alpha-1)s},$$
and the first term of the right side of \eqref{mi17} is also $O(e^{2(\alpha-1)s})$.
In summary, no matter which case happens, we always have
$$|U(s)|,\ |V(s)|\leq C(p,q)e^{(\tilde{\alpha}-1)s}\ (s\geq \ln (\hat{\psi}_1(0)\frac{\mu}{\sqrt{1+\mu^2}})+1),$$
namely,
$$|W-W_*|,\ |Z-Z_*|\leq C(p,q)e^{(\tilde{\alpha}-1)s}\ (s\geq \ln (\hat{\psi}_1(0)\frac{\mu}{\sqrt{1+\mu^2}})+1).$$
Applying the transform \eqref{mi10}, \eqref{mi12} again, we get
\begin{equation}
    |\psi(r)-r^{\alpha}|\leq C(p,q)r^{\tilde{\alpha}},\ |\psi'(r)-\alpha r^{\alpha-1}|\leq C(p,q)r^{\tilde{\alpha}-1} \label{mi19}
\end{equation}
for all $r\geq e\hat{\psi}_1(0)\frac{\mu}{\sqrt{1+\mu^2}}$. By continuity, \eqref{mi19} holds for all $r\geq \hat{\psi}_1(0)\frac{\mu}{\sqrt{1+\mu^2}}$.\\
The higher order estimates (i.e. \eqref{mi8} with $m\geq 2$) follows by differentiating \eqref{mi5} and induction on $m$. The start point $m=2$ is a direct consequence of \eqref{mi19}.

\QED

\section{Admissible Flow}
\label{ad}
From now on in this paper, we fix the following constants:
$$\begin{cases}
    n\in \Z,\ n\geq 8;\\
    p,\ q\in \Z,\ p+q=n,\ p,\ q\geq 2\ (\text{if}\ n=8, \text{we require further that}\ p,\ q\geq 3);\\
    l\in \Z,\ l\geq 2,
\end{cases}$$
and set the following real constants to be determined:
\begin{equation}
    \begin{cases}
    \Lambda\gg 1\ (\text{depending on}\ p,q,l);\\
    0<\rho\ll 1\ (\text{depending on}\ p,q,l,\Lambda);\\
    \beta\gg 1\ (\text{depending on}\ p,q,l,\Lambda);\\
    R\gg 1\ (\text{depending on}\ p,q,l,\Lambda);\\
    t_0<0,\ |t_0|\ll 1\ (\text{depending on}\ p,q,l,\Lambda,\rho,\beta,R).
\end{cases} \label{ad12}
\end{equation}
Assume $t_0<\mathring{t}<0$, and there is a one-parameter family of smooth hypersurfaces $\{\Sigma_t\}_{t_0\leq t\leq \mathring{t}}$ in $\R^n$. We say $\{\Sigma_t\}_{t_0\leq t\leq \mathring{t}}$ is \textbf{admissible} if
\begin{enumerate}
    \item \label{ad_1} $\{\Sigma_t\}_{t_0\leq t\leq \mathring{t}}$ depends smoothly on both space and time, moves by its mean curvature, and every time-slice $\Sigma_t$ is a compact, embedded, $O(p)\times O(q)$ invariant hypersurface.
    \item \label{ad_2} For some $\epsilon>0$ and all $t_0\leq t\leq \mathring{t}$, the profile curve of $\Sigma_t\cap B(0,2(1+\epsilon)\rho)$ can be parameterized by a single function, as
    \begin{equation}
        (x,\hat{u}(x,t)), \label{ad3}
    \end{equation}
    where $\hat{u}(x,t)$ is a positive smooth function, defined at least for $0\leq x\leq \frac{(1+\epsilon)\sqrt{-t}}{\sqrt{1+\mu^2}}$. Note that in this case (\ref{ad_1}) implies $\hat{u}(x,t)$ satisfies
    \begin{equation*}
        \partial_t \hat{u}=\frac{\hat{u}''}{1+\hat{u}'^2}+(p-1)\frac{\hat{u}'}{x}-(q-1)\frac{1}{\hat{u}}.
    \end{equation*}
    We denote by ' the derivative in space variables. Moreover, the even extension of $\hat{u}(\cdot,t)$ is smooth, and particularly $\hat{u}'(0,t)=0$.
    \item \label{ad_3} For some $\epsilon>0$ and all $t_0\leq t\leq \mathring{t}$, the profile curve of $\Sigma_t\cap (B(0,2(1+\epsilon)\rho)-\bar{B}(0,\frac{1}{2}(1-\epsilon)\beta(-t)^{\frac{1}{2}+\sigma_l}))$ (where $\bar{B}$ denotes a closed ball) can be parameterized by a single function, as
    \begin{equation}
        (\frac{x-\mu u(x,t)}{\sqrt{1+\mu^2}},\frac{\mu x+u(x,t)}{\sqrt{1+\mu^2}}), \label{ad1}
    \end{equation}
    where $u(x,t)$ is a smooth function, defined at least for $(1-\epsilon)\beta(-t)^{\frac{1}{2}+\sigma_l}\leq x\leq (1+\epsilon)\rho$, and
    $$\sigma_l=\frac{\lambda_l}{1-\alpha},\ \lambda_l=-\frac{1}{2}(1-\alpha)+l.$$
    (The meaning of $\lambda_l$ will be interpreted later, see Proposition \ref{propad1}.) In this case (\ref{ad_1}) implies $u(x,t)$ satisfies
    \begin{equation}
        \partial_t u=\frac{u''}{1+u'^2}+(p-1)\frac{\mu+u'}{x-\mu u}-(q-1)\frac{1-\mu u'}{\mu x+u}. \label{ad2}
    \end{equation}
    \item \label{ad_4} The following estimate
    \begin{equation}
        x^i|\partial_x^i u(x,t)|<\Lambda ((-t)^l x^\alpha+x^{2\lambda_l+1}) \label{ad5}
    \end{equation}
    holds for all $i=0,1,2$, $\beta(-t)^{\frac{1}{2}+\sigma_l}\leq x\leq \rho$, and $t_0\leq t\leq \mathring{t}$.
\end{enumerate}
Rescaling is a common technique for analysis of blow-up behaviors of mean curvature flows. In the later discussion, we roughly divide the space into three (time dependent) regions and dilate each region at different rates:
\begin{enumerate}
    \item The \textbf{outer region}: $\Sigma_t-B(0,R\sqrt{-t})$.
    \item The \textbf{intermediate region}: $\Sigma_t\cap (\bar{B}(0,R\sqrt{-t})-B(0,\beta(-t)^{\frac{1}{2}+\sigma_l}))$.
    \item The \textbf{tip region}: $\Sigma_t\cap \bar{B}(0,\beta(-t)^{\frac{1}{2}+\sigma_l})$.
\end{enumerate}
In the outer region, we mainly analyze the unrescaled mean curvature flow, especially, the function $u(x,t)$ defined in \eqref{ad1}, and the equation \eqref{ad2}.\\
In the intermediate region, we perform the ``type I" rescaling
\begin{equation}
    \Pi_s=(-t)^{-\frac{1}{2}}\Sigma_t|_{t=-e^{-s}}, \label{ad4}
\end{equation}
the time interval becoming $s_0\leq s\leq \mathring{s}$, where $s_0=-\ln(-t_0)$, $\mathring{s}=-\ln(-\mathring{t})$. Note that $s_0\gg 1$ iff $|t_0|\ll 1$.\\
Under rescaling \eqref{ad4}, the piece of the profile curve of $\Sigma_t$ which can be parameterized as \eqref{ad1} is now rescaled and parameterized as 
$$(\frac{y-\mu v(y,s)}{\sqrt{1+\mu^2}},\frac{\mu y+v(y,s)}{\sqrt{1+\mu^2}}),$$
with
\begin{equation}
    v(y,s)=(-t)^{-\frac{1}{2}}u(\sqrt{-t}y,t)|_{t=-e^{-s}}, \label{ad6}
\end{equation}
satisfying the equation 
\begin{equation}
        \partial_s v=\frac{v''}{1+v'^2}+(p-1)\frac{\mu+v'}{y-\mu v}-(q-1)\frac{1-\mu v'}{\mu y+v}+\frac{1}{2}(-yv'+v). \label{ad7}
\end{equation}
The condition \eqref{ad5} becomes
\begin{equation}
    y^i|\partial_y^i v(y,s)|<\Lambda e^{-\lambda_l s}(y^\alpha+y^{2\lambda_l+1}),\ i=0,1,2,\ \beta e^{-\sigma_l s}\leq y\leq \rho e^{\frac{s}{2}},\ s_0\leq s\leq \mathring{s}. \label{ad14}
\end{equation}
In this region, we mainly investigate the function $v(y,s)$ defined in \eqref{ad6}, and the equation \eqref{ad7}.\\
In the tip region, we perform the ``type II" rescaling
\begin{equation}
    \Gamma_\tau=(-t)^{-(\frac{1}{2}+\sigma_l)}\Sigma_t|_{t=-(2\sigma_l \tau)^{-\frac{1}{2\sigma_l}}}, \label{ad8}
\end{equation}
the time interval becoming $\tau_0\leq \tau\leq \mathring{\tau}$, where $\tau_0=(2\sigma_l)^{-1}(-t_0)^{-2\sigma_l}$, $\mathring{\tau}=(2\sigma_l)^{-1}(-\mathring{t})^{-2\sigma_l}$. Note that $\tau_0\gg 1$ iff $|t_0|\ll 1$.\\
Under rescaling \eqref{ad8}, the piece of the profile curve of $\Sigma_t$ which can be parameterized as \eqref{ad3} is now rescaled and parameterized as 
$$(z,\hat{w}(z,\tau)),$$
with
\begin{equation}
    \hat{w}(z,\tau)=(-t)^{-(\frac{1}{2}+\sigma_l)}\hat{u}((-t)^{\frac{1}{2}+\sigma_l}z,t)|_{t=-(2\sigma_l \tau)^{-\frac{1}{2\sigma_l}}}=e^{\sigma_l s}\hat{v}(e^{-\sigma_l s}z,s)|_{s=\frac{1}{2\sigma_l}\ln (2\sigma_l \tau)}, \label{ad9}
\end{equation}
satisfying the equation 
\begin{equation}
        \partial_\tau \hat{w}=\frac{\hat{w}''}{1+\hat{w}'^2}+(p-1)\frac{\hat{w}'}{z}-(q-1)\frac{1}{\hat{w}}+\frac{\frac{1}{2}+\sigma_l}{2\sigma_l \tau}(-z\hat{w}'+\hat{w}). \label{ad10}
\end{equation}
The piece of the profile curve of $\Sigma_t$ which can be parameterized as \eqref{ad1} is now rescaled and parameterized as 
\begin{equation}
    (\frac{z-\mu w(z,\tau)}{\sqrt{1+\mu^2}},\frac{\mu z+w(z,\tau)}{\sqrt{1+\mu^2}}), \label{ad17}
\end{equation}
with
\begin{equation}
    w(z,\tau)=(-t)^{-(\frac{1}{2}+\sigma_l)}u((-t)^{\frac{1}{2}+\sigma_l}z,t)|_{t=-(2\sigma_l \tau)^{-\frac{1}{2\sigma_l}}}=e^{\sigma_l s}v(e^{-\sigma_l s}z,s)|_{s=\frac{1}{2\sigma_l}\ln (2\sigma_l \tau)}, \label{ad16}
\end{equation}
satisfying the equation 
\begin{equation*}
        \partial_\tau w=\frac{w''}{1+w'^2}+(p-1)\frac{\mu+w'}{z-\mu w}-(q-1)\frac{1-\mu w'}{\mu z+w}+\frac{\frac{1}{2}+\sigma_l}{2\sigma_l \tau}(-zw'+w).
\end{equation*}
The condition \eqref{ad5} becomes
\begin{equation}
    z^i|\partial_z^i w(z,\tau)|<\Lambda (z^\alpha+\frac{z^{2\lambda_l+1}}{(2\sigma_l\tau)^l}),\ i=0,1,2,\ \beta\leq z\leq \rho (2\sigma_l \tau)^{\frac{1}{2}+\frac{1}{4\sigma_l}},\ \tau_0\leq \tau\leq \mathring{\tau}. \label{ad18}
\end{equation}
In this region, we mainly investigate the function $\hat{w}(z,\tau)$ defined in \eqref{ad9}, and the equation \eqref{ad10}.\\
In Vel\'{a}zquez's construction, the first step of a priori estimates is to show that the type-I rescaled flow satisfies the property that $v(y,s)$ tends exponentially to zero. For this reason, we linearize the right-hand side of \eqref{ad7} at $v=0$. More precisely, \eqref{ad7} can be rewritten as
\begin{equation}
    \partial_s v=\mathcal{L}v+\mathcal{Q}v, \label{ad13}
\end{equation}
where
$$\mathcal{L}v=v''+(\frac{n-2}{y}-\frac{y}{2})v'+(\frac{n-2}{y^2}+\frac{1}{2})v=(y^{n-2}e^{-\frac{y^2}{4}})^{-1}(y^{n-2}e^{-\frac{y^2}{4}}v')'+(\frac{n-2}{y^2}+\frac{1}{2})v$$
is the linear part of the right-hand side (note that $\mathcal{L}$ depends only on $n$), and
$$\mathcal{Q}v=-\frac{v'^2}{1+v'^2}v''+(n-2)\frac{(\frac{v}{y})^2(\frac{v}{y^2}+\frac{v'}{y})+(\mu-\mu^{-1})\frac{v}{y}\frac{v}{y^2}}{(1-\mu \frac{v}{y})(1+\mu^{-1}\frac{v}{y})}$$
is the remaining part. We will use the following properties of the linear operator $\mathcal{L}$ shown by Vel\'{a}zquez (see \cite{vel}, Proposition 2.1):
\begin{prop}
    \label{propad1}
    Let
    $$\mathbf{H}=L^2((0,+\infty);x^{n-2}e^{-\frac{x^2}{4}})\ (\text{a weighted}\ L^2\ \text{space}).$$
    For $u,v\in \mathbf{H}$, let $\langle u,v\rangle$ be the inner product in $\mathbf{H}$ and $|u|_{\mathbf{H}}$ be the norm in $\mathbf{H}$, namely,
    $$\langle u,v\rangle=\int_0^\infty u(x)v(x)x^{n-2}e^{-\frac{x^2}{4}}dx,\ |u|_{\mathbf{H}}=\sqrt{\langle u,u\rangle}.$$
    Then there exists a countable set of eigenvalues $\{\lambda_i\}_{i\geq 0}$ of $\mathcal{L}$, with corresponding normalized eigenfunctions $\{\varphi_i\}_{i\geq 0}\subset \mathbf{H}$, satisfying $\mathcal{L}\varphi_i+\lambda_i\varphi_i=0$. 
    They are given by:
    $$\lambda_i=-\frac{1}{2}(1-\alpha)+i,$$
    \begin{equation}
        \varphi_i(y)=c_i y^\alpha M(-i,\alpha+\frac{n-1}{2};\frac{y^2}{4}), \label{ad15}
    \end{equation}
    where $M(a,b;x)$ is the Kummer's function defined by
    $$M(a,b;x)=1+\sum_{j=1}^{+\infty} \frac{a(a+1)...(a+j-1)}{b(b+1)...(b+j-1)}\frac{x^j}{j!}$$
    and satisfying
    $$xM''(a,b;x)+(b-x)M'(a,b;x)-aM(a,b;x)=0,$$
    and $c_i>0$ is a normalizing constant s.t. $|\varphi_i|_{\mathbf{H}}=1$.\\
    Moreover, the set of eigenfunctions $\{\varphi_i\}_{i\geq 0}$ forms an orthonormal basis of $\mathbf{H}$.
\end{prop}
Note that
$$\lambda_0<\lambda_1<0<\lambda_2<\lambda_3<...,$$
i.e. $\lambda_2$ is the first positive eigenvalue of $\mathcal{L}$. The eigenfunction $\varphi_l$ can be written as
\begin{equation}
    \varphi_l(y)=c_ly^\alpha(1-K_{l,1}y^2+K_{l,2}y^4+...+(-1)^l K_{l,l}y^{2l}), \label{ad11}
\end{equation}
where $\{K_{l,j}\}_{1\leq j\leq l}$ are positive constants. The term of $\varphi_l$ of highest power has positive sign when $l$ is even, but has negative sign when $l$ is odd.

\section{Construction of Vel\'{a}zquez's Solutions}
The main idea of the construction can be summarized as follows. Our goal is to find a family of $O(p)\times O(q)$-invariant, smooth, compact hypersurfaces $\{\Sigma_t\}_{t_0\leq t<0}\subset \R^n$ moving by mean curvature, which blows up at $0\in \R^n$ and $t=0$ but remains smooth elsewhere; the type-I rescaled flow converges to the Lawson's cone $C_{p,q}$; the type-II rescaled flow converges to $\mathcal{M}_k$ for some $k\approx 1$, i.e. a smooth minimal hypersurface tangent to $C_{p,q}$ at infinity, introduced in Section \ref{mi}.\\
To achieve the goal, we begin with the construction of the initial hypersurface $\Sigma_{t_0}$. In fact we will choose a family of hypersurfaces $\{\Sigma_{t_0}^{\mathbf{a}}\}$, where $\mathbf{a}=(a_0,a_1,...,a_{l-1})$ is an $l$-dimensional parameter with $|\mathbf{a}|\ll 1$. The hypersurfaces $\{\Sigma_{t_0}^{\mathbf{a}}\}$ have uniformly bounded curvature outside the ball $B(0,\rho)$; they are ``close" to $C_{p,q}$ after type-I rescaling and ``close" to $\mathcal{M}_1$ after type-II rescaling. Then we can show for each $t_0<\mathring{t}<0$, there exists a parameter $\mathbf{a}$ such that the corresponding mean curvature flow $\{\Sigma_{t_0}^{\mathbf{a}}\}_{t\geq t_0}$ exists up to time $\mathring{t}$, and behaves in our prescribed way. In addition, this flow satisfies some uniform smooth estimates. Finally, by the compactness theory, we then get a solution of MCF which exists and has the desired behavior on the whole time interval $t_0\leq t<0$, and satisfies the same uniform estimates.\\
Now, let's write down $\gamma_{t_0}^\mathbf{a}$, the profile curve of $\{\Sigma_{t_0}^{\mathbf{a}}\}$ in the first quadrant, in an explicit way.\\
Let 
$$\mathbf{a}\in B(0,\beta^{\tilde{\alpha}-\alpha})\subset \R^l,$$
where $\tilde{\alpha}$ is defined in \eqref{mi20} with $\tilde{\alpha}-\alpha<0$. \\
In the region $\sqrt{(x^1)^2+(y^1)^2}>\frac{1}{2}$, $\gamma_{t_0}^\mathbf{a}$ is parameterized by a single smooth function as 
$$(x,\hat{u}(x,t_0)),\ \frac{1}{2\sqrt{1+\mu^2}}<x\leq 2,\ (\text{see \eqref{ad3}})$$
where
$$\hat{u}(x,t_0)=\begin{cases}
\begin{array}{ll}
   \frac{1}{2}x,  & \frac{1}{2\sqrt{1+\mu^2}}<x\leq \frac{2}{\sqrt{1+\mu^2}}-\delta, \\
   \text{smooth and concave}, & \frac{2}{\sqrt{1+\mu^2}}-\delta\leq x\leq \frac{2}{\sqrt{1+\mu^2}}+\delta,\\
   \sqrt{4-x^2},  & \frac{2}{\sqrt{1+\mu^2}}+\delta\leq x\leq 2,
\end{array} 
\end{cases}$$
with $0<\delta\ll 1$ (depending on $p,q$).\\
In the region $\sqrt{(x^1)^2+(y^1)^2}<1$, $\gamma_{t_0}^\mathbf{a}$ is parameterized by a single smooth function as 
$$(\frac{x-\mu u(x)}{\sqrt{1+\mu^2}},\frac{\mu x+u(x)}{\sqrt{1+\mu^2}}),\ \hat{\psi}_{1+a_0+a_1+...+a_{l-1}}(0)\frac{\mu}{\sqrt{1+\mu^2}}(-t_0)^{\frac{1}{2}+\sigma_l}\leq x<1,\ (\text{see \eqref{ad1}})$$
where 
\begin{equation}
    u(x,t_0)=u(x,t_0;\mathbf{a})=(-t_0)^{\frac{1}{2}+\sigma_l} \psi_{1+a_0+a_1+...+a_{l-1}}((-t_0)^{-(\frac{1}{2}+\sigma_l)}x)(1-\eta(\frac{(-t_0)^{-(\frac{1}{2}+\sigma_l)}x-\frac{1}{2}\beta}{\frac{1}{2}\beta})) \label{co2}
\end{equation}
$$+(-t_0)^{\frac{1}{2}+\lambda_l}(\frac{1}{c_l}\varphi_l(\frac{x}{\sqrt{-t_0}})+\sum_{j=0}^{l-1}\frac{a_j}{c_j} \varphi_j(\frac{x}{\sqrt{-t_0}}))\eta(\frac{(-t_0)^{-(\frac{1}{2}+\sigma_l)}x-\frac{1}{2}\beta}{\frac{1}{2}\beta})\eta(\frac{2\rho-x}{\rho}),$$
with $\eta:\R\rightarrow \R$ a smooth, non-decreasing function satisfying
\begin{equation}
    \eta(x)=\begin{cases}
\begin{array}{ll}
   0,  & x\leq 0, \\
   1,  & x\geq 1,
\end{array} 
\end{cases} \label{co4}
\end{equation}
and $\{c_j\}_{j\geq 0}$ the normalizing constants defined in Proposition \ref{propad1}.\\
For simplicity, we usually write $a_l=1$. According to \eqref{ad11}, 
$$(-t_0)^{\frac{1}{2}+\lambda_l}(\frac{1}{c_l}\varphi_l(\frac{x}{\sqrt{-t_0}})+\sum_{j=0}^{l-1}\frac{a_j}{c_j} \varphi_j(\frac{x}{\sqrt{-t_0}}))$$
$$=(\sum_{j=0}^l a_j)(-t_0)^l x^\alpha-(\sum_{j=1}^l K_{j,1}a_j)(-t_0)^{l-1} x^{\alpha+2}+(\sum_{j=2}^l K_{j,2}a_j)(-t_0)^{l-2} x^{\alpha+4}+...+(-1)^l K_{l,l}x^{\alpha+2l}$$
\begin{equation}
    =x^{2\lambda_l+1}\{(\sum_{j=0}^l a_j)(\frac{-t_0}{x^2})^l -(\sum_{j=1}^l K_{j,1}a_j)(\frac{-t_0}{x^2})^{l-1} +(\sum_{j=2}^l K_{j,2}a_j)(\frac{-t_0}{x^2})^{l-2} +...+(-1)^l K_{l,l}\}. \label{co16}
\end{equation}
Thus, for all $x>0$,
\begin{equation}
    x^i|\partial_x^i \{(-t_0)^{\frac{1}{2}+\lambda_l}(\frac{1}{c_l}\varphi_l(\frac{x}{\sqrt{-t_0}})+\sum_{j=0}^{l-1}\frac{a_j}{c_j} \varphi_j(\frac{x}{\sqrt{-t_0}}))\}|\leq C(n,l)((-t_0)^l x^\alpha+x^{2\lambda_l+1}),\ i=0,1,2, \label{co1}
\end{equation}
and by \eqref{mi7}, for $x\geq \hat{\psi}_{1+a_0+a_1+...+a_{l-1}}(0)\frac{\mu}{\sqrt{1+\mu^2}}(-t_0)^{\frac{1}{2}+\sigma_l}$,
\begin{equation}
    x^i|\partial_x^i (-t_0)^{\frac{1}{2}+\sigma_l} \psi_{1+a_0+a_1+...+a_{l-1}}((-t_0)^{-(\frac{1}{2}+\sigma_l)}x)|\leq C(p,q)(-t_0)^l x^\alpha,\ i=0,1,2. \label{co3}
\end{equation}
Therefore, according to \eqref{co1}, if $\Lambda\gg 1$ (depending on $n,l$), then \eqref{ad5} holds. Moreover, by \eqref{co1}, \eqref{co3}, and \eqref{co2} (the definition of $u$), for $\frac{1}{2}\beta(-t_0)^{\frac{1}{2}+\sigma_l}\leq x\leq 2\rho$,
$$|\frac{u(x,t_0)}{x}|,|u'(x,t_0)|\leq C(l,p,q)(\beta^{\alpha-1}+\rho^{2\lambda_l})\leq \frac{1}{2}\min\{\mu,\mu^{-1}\},$$
provided $\beta\gg 1$, $\rho\ll 1$ (depending on $l,p,q$). Thus the curve $\gamma_{t_0}^\mathbf{a}$ lies in the first quadrant, and the hypersurface $\Sigma_{t_0}^{\mathbf{a}}$ generated by $\gamma_{t_0}^\mathbf{a}$ is smooth, compact, embedded, satisfying the admissible condition. In addition,  $\Sigma_{t_0}^{\mathbf{a}}$ depends smoothly on $\mathbf{a}$.\\
Next, for each $\mathbf{a}\in B(0,\beta^{\tilde{\alpha}-\alpha})$, by the short-time existence (see \cite{hui}, Theorem 3.1), the mean curvature flow starting from $\Sigma_{t_0}^{\mathbf{a}}$ has a unique solution for a short time, denoted by $\{\Sigma_t^{\mathbf{a}}\}$. Define the set $\mathcal{O}\subset B(0,\beta^{\tilde{\alpha}-\alpha})\times [-t_0,0)$ as follows: $(\mathbf{a},\mathring{t})\in \mathcal{O}$ iff
\begin{enumerate}
    \item the corresponding (smooth) MCF $\{\Sigma_t^{\mathbf{a}}\}$ exists for $t_0\leq t\leq \mathring{t}$ and can be extended beyond $\mathring{t}$.
    \item $\{\Sigma_t^{\mathbf{a}}\}$ is admissible for $t_0\leq t\leq \mathring{t}$.
\end{enumerate}
From the smooth dependence of MCF on the initial data, $\mathcal{O}$ is a (relatively) open subset of $B(0,\beta^{\tilde{\alpha}-\alpha})\times [-t_0,0)$.\\
For $t_0\leq \mathring{t}<0$, let 
$$\mathcal{O}_{\mathring{t}}=\{\mathbf{a}\in B(0,\beta^{\tilde{\alpha}-\alpha})|\ (\mathbf{a},\mathring{t})\in \mathcal{O}\},$$
then $\mathcal{O}_{\mathring{t}}$ is an open subset of $B(0,\beta^{\tilde{\alpha}-\alpha})$, and is decreasing in $\mathring{t}$. Obviously, $\mathcal{O}_{t_0}=B(0,\beta^{\tilde{\alpha}-\alpha})$.\\
Recall the function $v(y,s)=v(y,s;\mathbf{a})$ \eqref{ad6} in the type-I rescaled flow $\{\Pi_s^{\mathbf{a}}\}$ \eqref{ad4}. We define a map $\Phi: \mathcal{O}\rightarrow \R^l$ by
$$\Phi(\mathbf{a},t)=e^{\lambda_l s_0}(\langle c_0  \tilde{v}(y,s;\mathbf{a}),\varphi_0(y)\rangle,\langle c_1  \tilde{v}(y,s;\mathbf{a}),\varphi_1(y)\rangle,...,\langle c_{l-1}  \tilde{v}(y,s;\mathbf{a}),\varphi_{l-1}(y)\rangle)_{s=-\ln(-t)},$$
where
\begin{equation}
    \tilde{v}(y,s;\mathbf{a})=\eta(e^{\sigma_l s}y-\beta)\eta(\rho e^{\frac{s}{2}}-y)v(y,s;\mathbf{a}), \label{co19}
\end{equation}
$\langle \cdot, \cdot \rangle$ is the inner product defined in Proposition \ref{propad1}, and $\eta$ is defined in \eqref{co4}. According to the admissible condition, $\Phi$ is well-defined, and continuous (actually smooth) on $\mathcal{O}$. For $t_0\leq t<0$, we also define
$$\Phi_t(\mathbf{a})=\Phi(\mathbf{a},t),\ \mathbf{a}\in \mathcal{O}_t.$$
When $t=t_0$, we have the following lemma (see \cite{gs}, Lemma 4.3):
\begin{lemma}
    \label{lemco1}
    If $s_0\gg 1$ (depending on $n,l,\rho,\beta$), then
    $$|\langle \eta(e^{\sigma_l s_0}y-\beta)\eta(\rho e^{\frac{s_0}{2}}-y) \varphi_i(y),\varphi_j(y) \rangle-\delta_{ij}|\leq C(n,l,\beta)e^{-(n-1+2\alpha)\sigma_l s_0},$$
    $$|(1-\eta(e^{\sigma_l s}y-\beta)\eta(\rho e^{\frac{s}{2}}-y))\varphi_i(y)|_{\mathbf{H}}\leq C(n,l,\beta)e^{-\frac{1}{2}(n-1+2\alpha)\sigma_l s_0},$$
    for all $0\leq i,j\leq l$, where $s_0=-\ln(-t_0)$.
\end{lemma}
By \eqref{co2}, 
\begin{equation}
    v(y,s_0;\mathbf{a})=e^{-\lambda_l s_0}(\frac{1}{c_l}\varphi_l(y)+\sum_{j=0}^{l-1} \frac{a_j}{c_j}\varphi_j(y)),\ \beta e^{-\sigma_l s_0}\leq y\leq \rho e^{\frac{s_0}{2}}, \label{co8}
\end{equation}
and we know from Lemma \ref{lemco1} that $\Phi_{t_0}$ converges uniformly to the identity map on $B(0,\beta^{\tilde{\alpha}-\alpha})$ as $t_0\nearrow 0$. Thus, if $|t_0|\ll 1$ (depending on $n,l,\rho,\beta$), we have 
$$\Phi_{t_0}^{-1}(0)\subset \subset B(0,\beta^{\tilde{\alpha}-\alpha}),$$
(in general, we say two sets $A\subset \subset B$ iff $A$ has compact closure, and $\bar{A}\subset B$,) and the topological degree
$$\deg(\Phi_{t_0},\mathcal{O}_{t_0},0)=\deg(\Phi_{t_0},B(0,\beta^{\tilde{\alpha}-\alpha}),0)=\deg(\id,B(0,\beta^{\tilde{\alpha}-\alpha}),0)=1.$$
When
\begin{equation}
    (\mathbf{a},\mathring{t})\in \mathcal{O}\ \text{and}\ \Phi_{\mathring{t}}(\mathbf{a})=0, \label{co5}
\end{equation}
and if \eqref{ad12} holds, we have the a priori estimates which are crucial for the extension of solution: 
\begin{prop}
    \label{propco1}
    Set
    \begin{equation}
        \varsigma=\min\{1,\frac{n-3+2\alpha}{2(1-\alpha)}\} \label{co7}
    \end{equation}
    (if $n=8$, we require further that $0<\varsigma<\frac{n-3+2\alpha}{2(1-\alpha)}$ and $\varsigma$ is fixed, e.g. $\varsigma=\frac{1}{7}$). \\
    If \eqref{co5} and \eqref{ad12} hold, then for all $t_0\leq t\leq \mathring{t}$,
    \begin{enumerate}
        \item \label{propco1_1} $|\mathbf{a}|\leq C(p,q,l,\Lambda,\beta)(-t_0)^{\varsigma \lambda_l}\leq \frac{1}{2}\beta^{\tilde{\alpha}-\alpha}$.
        \item \label{propco1_2} In $\R^n-\bar{B}(0,2\rho)$, the second fundamental form of $\Sigma_t^a$ is uniformly bounded by a constant $C(p,q,\rho)$.
        \item \label{propco1_3} The profile curve of $\Sigma_t\cap (B(0,3\rho)-\bar{B}(0,\frac{1}{3}\beta(-t)^{\frac{1}{2}+\sigma_l}))$ can be parameterized by a single function as \eqref{ad1}, with
        \begin{equation}
            x^i|\partial_x^i u(x,t)|\leq \frac{\Lambda}{2} ((-t)^l x^\alpha+x^{2\lambda_l+1}),\ i=0,1,2,\ \beta(-t)^{\frac{1}{2}+\sigma_l}\leq x\leq \rho, \label{co17}
        \end{equation}
        and for $x\geq \rho$ in this region,
        \begin{equation}
            \begin{cases}
            |\frac{u(x,t)}{x}|,|u'(x,t)|\leq \frac{1}{2}\min\{\mu,\mu^{-1}\},\\
            |u''(x,t)|\leq \frac{C(p,q)}{\rho}.
        \end{cases} \label{co11}
        \end{equation}
        \item \label{propco1_4} The profile curve of $\Sigma_t\cap B(0,3\rho)$ can be parameterized by a single function as \eqref{ad3}. Moreover, if we perform the type-II rescaling, the rescaled function $\hat{w}(z,\tau)$ defined in \eqref{ad9} satisfies
        \begin{equation}
            \begin{cases}
            \hat{\psi}_{\frac{1}{2}}(z)\leq \hat{w}(z,\tau)\leq \hat{\psi}_{2}(z),\\
            \hat{w}'(0,\tau)=0\leq \hat{w}'(z,\tau)\leq \mu+1,\\
            |\hat{w}''(z,\tau)|\leq C(p,q),
        \end{cases} \label{co18}
        \end{equation}
        for $0\leq z\leq \frac{2\beta}{\sqrt{1+\mu^2}}$.
    \end{enumerate}
\end{prop}
Furthermore, we have the following smooth estimates near the origin, which describe precisely the asymptotic blow-up behavior of the solution:
\begin{prop}
    \label{propco2}
    If \eqref{co5} and \eqref{ad12} hold, then there exists $k\in \R$ satisfying
    \begin{equation}
        |k-1|\leq C(p,q,l,\Lambda,\beta)(-t_0)^{\varsigma \lambda_l}, \label{co6}
    \end{equation}
    such that for any non-negative integers $m,r$, and all $t_0\leq t\leq \mathring{t}$, the following estimates hold:
    \begin{enumerate}
        \item In the outer region, the function $u(x,t)$ defined in \eqref{ad1} satisfies
        \begin{equation}
            x^{m+2r}|\partial_x^m \partial_t^r (u(x,t)-\frac{k}{c_l}(-t)^{\lambda_l+\frac{1}{2}}\varphi_l(\frac{x}{\sqrt{-t}}))|\leq C(p,q,l,\Lambda,m,r)(R^{-2}+\rho^{2\lambda_l})x^{2\lambda_l+1} \label{co12}
        \end{equation}
        for $R\sqrt{-t}\leq x\leq \frac{3}{4}\rho$.
        \item In the intermediate region, if we perform the type-I rescaling, then the function $v(y,s)$ defined in \eqref{ad6} satisfies
        \begin{equation}
            y^{m+2r} |\partial_y^m \partial_s^r (v(y,s)-\frac{k}{c_l}e^{-\lambda_l s}\varphi_l(y))|\leq C(p,q,l,\Lambda,\beta,R,m,r)e^{-(1+\kappa)\lambda_l s} y^\alpha \label{co13}
        \end{equation}
        for $e^{-\frac{1}{2}\sigma_l s}\leq y\leq R$, where
        \begin{equation}
            \kappa=\min\{\frac{1}{2},\frac{n-1+2\alpha}{6(1-\alpha)},\varsigma,\frac{1}{\lambda_l+1}\}, \label{co9}
        \end{equation}
        and
        \begin{equation}
            y^{m+2r} |\partial_y^m \partial_s^r (v(y,s)-e^{-\sigma_l s}\psi_k(e^{\sigma_l s}y))|\leq C(p,q,l,\Lambda,m,r)\beta^{\tilde{\alpha}-\alpha}e^{-2\varrho \sigma_l(s-s_0)}e^{-\lambda_ls} y^\alpha \label{co14}
        \end{equation}
        for $\frac{3}{2}\beta e^{-\sigma_l s}\leq y\leq e^{-\frac{1}{2}\sigma_l s}$, where
        \begin{equation}
            \varrho=\min\{\frac{\kappa}{2}(1-\alpha),\frac{1}{5}\}. \label{co10}
        \end{equation}
        \item In the tip region, if we perform the type-II rescaling, then the function $\hat{w}(z,\tau)$ defined in \eqref{ad9} satisfies
        \begin{equation}
            (1+z)^{m+2r}|\partial_z^m \partial_\tau ^r (\hat{w}(z,\tau)-\hat{\psi}_k(z))|\leq C(p,q,l,m,r)\beta^{\tilde{\alpha}-\alpha}(\frac{\tau}{\tau_0})^{-\varrho}(1+z)^\alpha \label{co15}
        \end{equation}
        for $0\leq z\leq \frac{2\beta}{\sqrt{1+\mu^2}}$.
    \end{enumerate}
\end{prop}
Note that our estimates in Proposition \ref{propco2} are slightly different from those in \cite{gs} that our ones cover the region where $t$ is close to the initial time $t_0$.

\begin{lemma}
    \label{lemco2}
    Assume \eqref{ad12} holds. If there is a sequence $\{(\mathbf{a}_j,t_j)\}_{j\geq 1}\subset \mathcal{O}$ s.t. $\Phi_{t_j}(\mathbf{a}_j)=0$ for all $j\geq 1$, and $(\mathbf{a}_j,t_j)\rightarrow (\mathbf{a}_\infty,t_\infty)$ with $t_\infty<0$, then $(\mathbf{a}_\infty,t_\infty)\in \mathcal{O}$, and $\Phi_{t_\infty}(\mathbf{a}_\infty)=0$.
\end{lemma}
\proof
First, by (\ref{propco1_1}) of Proposition \ref{propco1}, $\mathbf{a}_\infty\in B(0,\beta^{\tilde{\alpha}-\alpha})$. Then, from the smooth dependence of MCF on the initial data, the flow $\{\Sigma_t^{\mathbf{a}_\infty}\}$ exists on $t_0\leq t<t_\infty$, and satisfies all the conditions in Proposition \ref{propco1} on this time interval. But obviously those conditions imply a uniform bound on the second fundamental form of $\{\Sigma_t^{\mathbf{a}_\infty}\}_{t_0\leq t<t_\infty}$, and thus the flow can be smoothly extended past time $t_\infty$, with the help of Theorem 8.1 in \cite{hui}. The statements in Proposition \ref{propco1} still hold for $t_0\leq t\leq t_\infty$. It's straightforward to check that the conditions in Proposition \ref{propco1} imply the flow is admissible, i.e. $(\mathbf{a}_\infty,t_\infty)\in \mathcal{O}$, and $\Phi_{t_\infty}(\mathbf{a}_\infty)=0$ follows from the continuity of $\Phi$. (If $t_\infty=t_0$, things become trivial.)

\QED

A direct consequence of Lemma \ref{lemco2} is, for any $t_0\leq \mathring{t}<0$, $\Phi^{-1}(0)\cap (\R^l\times [t_0,\mathring{t}])$ is a compact subset of $\mathcal{O}$. Thus $\Phi_{\mathring{t}}^{-1}(0)\subset \subset \mathcal{O}_{\mathring{t}}$, and $\deg(\Phi_{\mathring{t}},\mathcal{O}_{\mathring{t}},0)$ is well-defined. Moreover, by the homotopy invariance of degree, 
$$\deg(\Phi_{\mathring{t}},\mathcal{O}_{\mathring{t}},0)=\deg(\Phi_{t_0},\mathcal{O}_{t_0},0)=1,\ \text{and}\ \Phi_{\mathring{t}}^{-1}(0)\neq \emptyset,$$
as long as $\mathcal{O}_{\mathring{t}}\neq \emptyset$. To show $\mathcal{O}_{\mathring{t}}$ is indeed non-empty, we define
$$t^*=\sup\{t|\ t_0\leq t<0,\ \mathcal{O}_t\neq \emptyset\}.$$
If $t^*<0$, then $\mathcal{O}_{t^*}=\emptyset$. Choose a sequence $\{(\mathbf{a}_j,t_j)\}_{j\geq 1}\subset \Phi^{-1}(0)$ such that $t_j\nearrow t^*$, and we may assume (by passing to a subsequence) $\mathbf{a}_j\rightarrow \mathbf{a}^*$. Using Lemma \ref{lemco2} again, we know $\mathbf{a}^*\in \mathcal{O}_{t^*}$, which is a contradiction, i.e. $\mathcal{O}_t\neq \emptyset$ for all $t_0\leq t<0$.\\
Now it's time to state and prove the existence result of our desired solutions.
\begin{theorem}
    \label{thmco1}
    Assuming \eqref{ad12} holds, then there exists an $\mathbf{a}\in B(0,\beta^{\tilde{\alpha}-\alpha})$ s.t. $\{\Sigma_t^{\mathbf{a}}\}$ exists for all $t_0\leq t<0$, satisfying all the conditions in Propositions \ref{propco1} and \ref{propco2} on this time interval.
\end{theorem}
\proof
Choose a sequence $\{(\mathbf{a}_j,t_j)\}_{j\geq 1}\subset \Phi^{-1}(0)$ s.t. $t_j\nearrow 0$, and we may assume $\mathbf{a}_j\rightarrow \mathbf{a}\in B(0,\beta^{\tilde{\alpha}-\alpha})$. Again, by the smooth dependence of MCF on the initial data, the flow $\{\Sigma_t^{\mathbf{a}}\}$ exists on $t_0\leq t<0$, and satisfies all the conditions in Proposition \ref{propco1} on this time interval. In addition, the estimates in Proposition \ref{propco2} hold for $\{\Sigma_t^{\mathbf{a}_j}\}_{t_0\leq t\leq t_j}$ with $k$ replaced by some $k_j$ satisfying \eqref{co6}. We may again assume $k_j\rightarrow k_\infty$ satisfying \eqref{co6}. Thus the estimates in Proposition \ref{propco2} hold for $\{\Sigma_t^{\mathbf{a}}\}_{t_0\leq t<0}$ with $k$ replaced by $k_\infty$. 

\QED

\section{$C^0$ Estimates}
In this and the next sections, we finish the proof of Propositions \ref{propco1} and \ref{propco2}. The $C^0$ estimates in these Propositions are the main topic of this section.\\
From now on in this and the next sections, we always assume \eqref{co5} holds for some $t_0<\mathring{t}<0$.\\
We start from the following estimate in the intermediate region:
\begin{prop}
    \label{propes1}
    If $0<\rho\ll 1$, $\beta\gg 1$ (depending on $p,q,l,\Lambda$), and $s_0\gg 1$ (depending on $n,l,\rho,\beta$), then (\ref{propco1_1}) of Proposition \ref{propco1} hold. Moreover, there exists $k\in \R$ satisfying \eqref{co6}, such that for any $\vartheta\in (0,1)$, any $R\geq 1$, and all $s_0\leq s\leq \mathring{s}=-\ln(-\mathring{t})$, if we require further that $s_0\gg 1$ (depending on $R,\vartheta$), then the function $v(y,s)$ defined in \eqref{ad6} satisfies
    \begin{equation}
        |v(y,s)-\frac{k}{c_l}e^{-\lambda_l s}\varphi_l(y)|\leq C(p,q,l,\Lambda,\beta,R)e^{-(1+\tilde{\kappa})\lambda_ls} y^\alpha \label{es23}
    \end{equation}
    for $\frac{1}{2} e^{-\vartheta \sigma_l s}\leq y\leq 2R$, where
    $$\tilde{\kappa}=\min\{1-\vartheta,(1-\vartheta)\frac{n-3+2\alpha}{1-\alpha},(1-\vartheta)\frac{n-1+2\alpha}{3(1-\alpha)},\varsigma,\frac{1}{\lambda_l+1}\}$$
    $$=\min\{1-\vartheta,(1-\vartheta)\frac{n-1+2\alpha}{3(1-\alpha)},\varsigma,\frac{1}{\lambda_l+1}\}.$$
\end{prop}
The idea of proving Proposition \ref{propes1} mainly comes from \cite{hv}, especially Sections 4 and 5, which were devoted to the construction of an example of rotationally symmetric type-II blow up solution of a semilinear heat equation. The computation is quite lengthy, but we present it here for the sake of completeness. \cite{gs} gave an easier proof but could not cover the case where $n=8$. We linearize the right-hand side of \eqref{ad7} at $v=0$ and do Fourier expansion (see \eqref{ad13}); the $\varphi_l$ term is regarded as the ``main frequency", and the terms of ``lower frequency", i.e. $\varphi_j$ terms with $0\leq j\leq l-1$, are controlled using \eqref{co5}.\\
Let $\mathbf{X}$ be the completion of $C_c^\infty((0,+\infty))$ w.r.t. the norm
$$|\phi|_{\mathbf{X}}=(\int_0^\infty (\phi(x)^2+\phi'(x)^2)x^{n-2}e^{-\frac{x^2}{4}}dx)^{\frac{1}{2}},\ \phi\in C_c^\infty((0,+\infty)).$$
Denote by $\mathbf{X}^*$ the dual space of $\mathbf{X}$, and define the norm on $\mathbf{X}^*$ in a standard way as
$$|u|_{\mathbf{X}^*}=\sup_{\phi\in \mathbf{X},\ |\phi|_{\mathbf{X}}=1}|u(\phi)|.$$
We need the following lemma about the eigenfunctions $\{\varphi_j\}_{j\geq 0}$ defined in Proposition \ref{propad1}:
\begin{lemma}
    \label{lemes1}
    For all $j\geq 0$, $\varphi_j\in \mathbf{X}$, and
    $$|\varphi_j|_{\mathbf{X}}^2\leq C(n)(1+|\lambda_j|).$$
    
\end{lemma}
\proof
According to \eqref{ad11}, $\varphi_j,\varphi'_j\in \mathbf{H}$. For $h\gg 1$, let
$$\varphi_{j,h}(x)=\eta(\frac{2h-x}{h})\eta(hx-1)\varphi_j(x),$$
where $\eta$ is the cut-off function defined in \eqref{co4}. Clearly, $\varphi_{j,h}\in C_c^\infty ((0,+\infty))$, and $\varphi_{j,h}\rightarrow \varphi_j$ in $\mathbf{H}$ as $h\rightarrow +\infty$. Look at its derivative:
$$\varphi_{j,h}'(x)=\eta(\frac{2h-x}{h})\eta(hx-1)\varphi_j'(x)+(h\eta'(hx-1)-\frac{1}{h}\eta'(\frac{2h-x}{h}))\varphi_j(x).$$
Again, $\eta(\frac{2h-x}{h})\eta(hx-1)\varphi_j'(x)\rightarrow \varphi_j'(x)$ in $\mathbf{H}$ as $h\rightarrow +\infty$, and
$$|(h\eta'(hx-1)-\frac{1}{h}\eta'(\frac{2h-x}{h}))\varphi_j(x)|_{\mathbf{H}}^2\leq C(h^2\int_{h^{-1}}^{2h^{-1}}x^{2\alpha} x^{n-2}e^{-\frac{x^2}{4}}dx+\frac{1}{h^2}\int_h^{2h} x^{4\lambda_l+2} x^{n-2}e^{-\frac{x^2}{4}}dx)$$
$$\leq C(h^{3-2\alpha-n}+h^{4\lambda_l+n-1}e^{-\frac{h^2}{4}})\rightarrow 0\ (h\rightarrow +\infty),$$
where $C$ is independent of $h$. Thus, $\varphi_{j,h}\rightarrow \varphi_j$ in $\mathbf{X}$ as $h\rightarrow +\infty$, i.e. $\varphi_j\in \mathbf{X}$. \\
According to Proposition 3.1 of \cite{gs}, if $u\in \mathbf{X}$, then $\frac{u(x)}{x}\in \mathbf{H}$, and we can define a bilinear form on $\mathbf{X}$ as
$$B(u,v)=\int_0^{+\infty} (u'(x)v'(x)-(\frac{n-2}{x^2}+\frac{1}{2})u(x)v(x))x^{n-2}e^{-\frac{x^2}{4}}dx, \ u,v\in \mathbf{X},$$
satisfying the estimate
\begin{equation}
    B(u,u)\geq \frac{n^2-10n+17}{(n-3)^2}|u'|_{\mathbf{H}}^2-\frac{3n-7}{2(n-3)}|u|_{\mathbf{H}}^2, \ u\in \mathbf{X}. \label{es9}
\end{equation}
Letting $u=v=\varphi_j$, and integrating $\int_0^{+\infty} \varphi_j'(x)\varphi_j'(x)x^{n-2}e^{-\frac{x^2}{4}}dx$ by parts, (noting that $\varphi_j(x)\varphi_j'(x)x^{n-2}e^{-\frac{x^2}{4}}\rightarrow 0$ as $x\searrow 0$, $x\rightarrow +\infty$) we find that
$$B(\varphi_j,\varphi_j)=-\langle \mathcal{L}\varphi_j,\varphi_j \rangle=\lambda_j|\varphi_j|_{\mathbf{H}}^2,$$
and by \eqref{es9},
$$|\varphi'_j|_{\mathbf{H}}^2\leq C(n)(1+|\lambda_j|)|\varphi_j|_{\mathbf{H}}^2=C(n)(1+|\lambda_j|).$$
Therefore,
$$|\varphi_j|_{\mathbf{X}}^2\leq C(n)(1+|\lambda_j|).$$

\QED

\noindent{\it Proof of Proposition \ref{propes1}}: 
Let
$$\tilde{v}(y,s)=\eta(e^{\sigma_l s}y-\beta)\eta(\rho e^{\frac{s}{2}}-y) v(y,s)$$
(see \eqref{co19}), then according to \eqref{ad13}, 
\begin{equation}
    (\partial_s-\mathcal{L})\tilde{v}:=f(y,s):=f_1+f_2+f_3, \label{es10}
\end{equation}
with
$$f_1(y,s)=\eta(e^{\sigma_l s}y-\beta)\eta(\rho e^{\frac{s}{2}}-y) \mathcal{Q}v(y,s),$$
$$f_2(y,s)=\eta'(e^{\sigma_l s}y-\beta)e^{\sigma_l s}(-2v'(y,s)+(-\frac{n-2}{y}+(\sigma_l+\frac{1}{2})y)v(y,s))-\eta''(e^{\sigma_l s}y-\beta)e^{2\sigma_l s}v(y,s),$$
$$f_3(y,s)=\eta'(\rho e^{\frac{s}{2}}-y)((\frac{\rho}{2}e^{\frac{s}{2}}-\frac{y}{2}+\frac{n-2}{y})v(y,s)+2v'(y,s))-\eta''(\rho e^{\frac{s}{2}}-y)v(y,s).$$
If $0<\rho\ll 1\ll \beta$ (depending on $p,q,l,\Lambda$) and $s_0\gg 1$ (depending on $n,l,\rho,\beta$), the following estimates of $f$ hold for all $s_0\leq s\leq \mathring{s}$:
\begin{align}
    |f_1(y,s)|\leq C(p,q,\Lambda) e^{-2\lambda_l s}y^{-3}(y^{2\alpha}+y^{4\lambda_l+2}) \chi_{(\beta e^{-\sigma_l s},\rho e^{\frac{s}{2}})}(y), \label{es1}\\
    |f_2(y,s)|\leq C(n,l,\Lambda,\beta) e^{-\lambda_l s}y^{\alpha-2}\chi_{(\beta e^{-\sigma_l s},(\beta+1) e^{-\sigma_l s})}(y), \label{es2}\\
    |f_3(y,s)|\leq C(n,\Lambda) e^{-\lambda_l s}y^{2\lambda_l+2}\chi_{(\rho e^{\frac{s}{2}}-1,\rho e^{\frac{s}{2}})}(y). \label{es3}
\end{align}
In fact, according to the admissible condition \eqref{ad14}, 
\begin{equation}
    |\frac{v(y,s)}{y}|,|v'(y,s)|<\Lambda(\beta^{\alpha-1}+\rho^{2\lambda_l})\leq \frac{1}{2}\min\{\mu,\mu^{-1}\},\ \beta e^{-\sigma_l s}\leq y\leq \rho e^{\frac{s}{2}},\ s_0\leq s\leq \mathring{s}, \label{es4}
\end{equation}
if $0<\rho\ll 1\ll \beta$ (depending on $p,q,l,\Lambda$). Therefore, \eqref{es1} is verified by putting \eqref{ad14} and \eqref{es4} into the definition of $\mathcal{Q}$; \eqref{es2} and \eqref{es3} are also proved directly by \eqref{ad14}. \\
We claim that
\begin{equation}
    |f(\cdot,s)|_{\mathbf{X}^*} \leq C(p,q,l,\Lambda,\beta)e^{-\lambda_l(1+\varsigma)s}, \label{es5}
\end{equation}
for all $s_0\leq s\leq \mathring{s}$, provided $s_0\gg 1$ (depending on $n,l,\rho,\beta$).\\
Take any $\phi\in C_c^\infty((0,+\infty))$. Let's estimate $f_1$ first:
$$|\int_0^{+\infty} f_1(y,s)\phi(y)y^{n-2}e^{-\frac{y^2}{4}}dy|$$
$$\leq C(p,q,\Lambda) e^{-2\lambda_l s}(\int_{\beta e^{-\sigma_l s}}^{\rho e^{\frac{s}{2}}} |\phi(y)|y^{2\alpha-3}y^{n-2}e^{-\frac{y^2}{4}}dy+\int_{\beta e^{-\sigma_l s}}^{\rho e^{\frac{s}{2}}} |\phi(y)|y^{4\lambda_l-1}y^{n-2}e^{-\frac{y^2}{4}}dy)$$
$$:=C(p,q,\Lambda) e^{-2\lambda_l s}(I_1+I_2),$$
$$I_2\leq \int_0^{+\infty} |\phi(y)|y^{4\lambda_l-1}y^{n-2}e^{-\frac{y^2}{4}}dy\leq (\int_0^{+\infty} \phi(y)^2y^{n-2}e^{-\frac{y^2}{4}}dy)^{\frac{1}{2}}(\int_0^{+\infty} y^{8\lambda_l-2}y^{n-2}e^{-\frac{y^2}{4}}dy)^{\frac{1}{2}}$$
$$=C(n,l)|\phi|_{\mathbf{H}}\leq C(n,l)|\phi|_{\mathbf{X}}.$$
To estimate $I_1$, assuming $n\geq 9$ first, we integrate by parts (by choosing the indefinite integral of $y^{2\alpha+n-5}\chi_{(\beta e^{-\sigma_l s},\rho e^{\frac{s}{2}})}(y)$ which takes 0 at 0), and use the fact that $||\phi|'(y)|=|\phi'(y)|$ for a.e. $y\in (0,+\infty)$:
$$I_1=-\frac{1}{2\alpha+n-4}\int_{\beta e^{-\sigma_l s}}^{\rho e^{\frac{s}{2}}} y^{2-n}(y^{2\alpha+n-4}-(\beta e^{-\sigma_l s})^{2\alpha+n-4}) (|\phi|'(y)-\frac{y}{2}|\phi(y)|)y^{n-2}e^{-\frac{y^2}{4}}dy$$
$$-\frac{1}{2\alpha+n-4}((\rho e^{\frac{s}{2}})^{2\alpha+n-4}-(\beta e^{-\sigma_l s})^{2\alpha+n-4})\int_{\rho e^{\frac{s}{2}}}^{+\infty} y^{2-n} (|\phi|'(y)-\frac{y}{2}|\phi(y)|)y^{n-2}e^{-\frac{y^2}{4}}dy$$
$$\leq C(n)|\phi|_{\mathbf{X}}\{(\int_{\beta e^{-\sigma_l s}}^{\rho e^{\frac{s}{2}}} y^{2-n}(y^{2\alpha+n-4}-(\beta e^{-\sigma_l s})^{2\alpha+n-4})^2(1+y^2)e^{-\frac{y^2}{4}}dy)^{\frac{1}{2}}$$
$$+(\rho e^{\frac{s}{2}})^{2\alpha+n-4}(\int_{\rho e^{\frac{s}{2}}}^{+\infty} y^{2-n} (1+y^2)e^{-\frac{y^2}{4}}dy)^{\frac{1}{2}}\}:=C(n)|\phi|_{\mathbf{X}}(I_{1,1}+I_{1,2}),$$
(if $s_0\gg 1$ (depending on $n,l,\rho,\beta$))
$$I_{1,1}\leq \sqrt{2}(\int_{\beta e^{-\sigma_l s}}^{\rho e^{\frac{s}{2}}} y^{4\alpha+n-6} (1+y^2)e^{-\frac{y^2}{4}}dy+(\beta e^{-\sigma_l s})^{4\alpha+2n-8}\int_{\beta e^{-\sigma_l s}}^{\rho e^{\frac{s}{2}}} y^{2-n} (1+y^2)e^{-\frac{y^2}{4}}dy)^{\frac{1}{2}}$$
$$\leq \begin{cases}
\begin{array}{ll}
     C(n), & 4\alpha+n-5>0\ (\text{if}\ s_0\gg 1\ (\text{depending on}\ n,l,\beta)), \\
    C(n)(\beta e^{-\sigma_l s})^{\frac{1}{2}(4\alpha+n-5)}, & 4\alpha+n-5<0.
\end{array}
\end{cases}$$
$$I_{1,2}\leq C(n)(\rho e^{\frac{s}{2}})^{2\alpha+n-4}e^{-\rho e^{\frac{s}{2}}}\leq C(n)\ (\text{if}\ s_0\gg 1\ (\text{depending on}\ n,\rho)).$$
When $n=8$, $2\alpha+n-5=-1$, and similarly we have
$$I_1=-\int_{\beta e^{-\sigma_l s}}^{\rho e^{\frac{s}{2}}} y^{2-n}(\ln y-\ln (\beta e^{-\sigma_l s})) (|\phi|'(y)-\frac{y}{2}|\phi(y)|)y^{n-2}e^{-\frac{y^2}{4}}dy$$
$$-(\ln(\rho e^{\frac{s}{2}})-\ln(\beta e^{-\sigma_l s}))\int_{\rho e^{\frac{s}{2}}}^{+\infty} y^{2-n} (|\phi|'(y)-\frac{y}{2}|\phi(y)|)y^{n-2}e^{-\frac{y^2}{4}}dy$$
$$\leq C(n)|\phi|_{\mathbf{X}}\{(\int_{\beta e^{-\sigma_l s}}^{\rho e^{\frac{s}{2}}} y^{2-n}(\ln y-\ln (\beta e^{-\sigma_l s}))^2(1+y^2)e^{-\frac{y^2}{4}}dy)^{\frac{1}{2}}$$
$$+(|\ln(\rho e^{\frac{s}{2}})|+|\ln(\beta e^{-\sigma_l s})|)(\int_{\rho e^{\frac{s}{2}}}^{+\infty} y^{2-n} (1+y^2)e^{-\frac{y^2}{4}}dy)^{\frac{1}{2}}\}:=C(n)|\phi|_{\mathbf{X}}(I'_{1,1}+I'_{1,2}).$$
For $0<\xi\ll 1$, 
$$I'_{1,1}\leq \sqrt{2}(\int_{\beta e^{-\sigma_l s}}^{\rho e^{\frac{s}{2}}} y^{2-n} (\ln y)^2(1+y^2)e^{-\frac{y^2}{4}}dy+(\ln(\beta e^{-\sigma_l s}))^2\int_{\beta e^{-\sigma_l s}}^{\rho e^{\frac{s}{2}}} y^{2-n} (1+y^2)e^{-\frac{y^2}{4}}dy)^{\frac{1}{2}}$$
$$\leq C(n,\xi)(\beta e^{-\sigma_l s})^{\xi-\frac{5}{2}}\ (\text{if}\ s_0\gg 1\ (\text{depending on}\ n,l,\beta,\xi)),$$
$$I'_{1,2}\leq C(n)(|\ln(\rho e^{\frac{s}{2}})|+|\ln(\beta e^{-\sigma_l s})|)e^{-\rho e^{\frac{s}{2}}}\leq C(n)\ (\text{if}\ s_0\gg 1\ (\text{depending on}\ n,l,\rho,\beta)).$$
No matter which case, we always have (recalling the definition of $\varsigma$, \eqref{co7})
\begin{equation}
    |f_1(\cdot,s)|_{\mathbf{X}^*} \leq C(p,q,l,\Lambda)e^{-\lambda_l(1+\varsigma)s}. \label{es6}
\end{equation}
Similarly, we can estimate $f_2$ (by choosing the indefinite integral of $y^{\alpha+n-4}\chi_{(\beta e^{-\sigma_l s},(\beta+1) e^{-\sigma_l s})}(y)$ which takes 0 at 0):
$$|\int_0^{+\infty} f_2(y,s)\phi(y)y^{n-2}e^{-\frac{y^2}{4}}dy|\leq C(n,l,\Lambda,\beta) e^{-\lambda_l s}\int_{\beta e^{-\sigma_l s}}^{(\beta+1) e^{-\sigma_l s}} |\phi(y)|y^{\alpha-2}y^{n-2}e^{-\frac{y^2}{4}}dy$$
$$=C(n,l,\Lambda,\beta) e^{-\lambda_l s}(-\frac{1}{\alpha+n-3}\int_{\beta e^{-\sigma_l s}}^{(\beta+1) e^{-\sigma_l s}} y^{2-n}(y^{\alpha+n-3}-(\beta e^{-\sigma_l s})^{\alpha+n-3}) (|\phi|'(y)-\frac{y}{2}|\phi(y)|)y^{n-2}e^{-\frac{y^2}{4}}dy$$
$$-\frac{1}{\alpha+n-3}(((\beta+1) e^{-\sigma_l s})^{\alpha+n-3}-(\beta e^{-\sigma_l s})^{\alpha+n-3})\int_{(\beta+1) e^{-\sigma_l s}}^{+\infty} y^{2-n} (|\phi|'(y)-\frac{y}{2}|\phi(y)|)y^{n-2}e^{-\frac{y^2}{4}}dy)$$
$$\leq C(n,l,\Lambda,\beta) e^{-\lambda_l s}|\phi|_{\mathbf{X}}\{(\int_{\beta e^{-\sigma_l s}}^{(\beta+1) e^{-\sigma_l s}} y^{2-n}(y^{\alpha+n-3}-(\beta e^{-\sigma_l s})^{\alpha+n-3})^2(1+y^2)e^{-\frac{y^2}{4}}dy)^{\frac{1}{2}}$$
$$+(((\beta+1) e^{-\sigma_l s})^{\alpha+n-3}-(\beta e^{-\sigma_l s})^{\alpha+n-3})(\int_{(\beta+1) e^{-\sigma_l s}}^{+\infty} y^{2-n} (1+y^2)e^{-\frac{y^2}{4}}dy)^{\frac{1}{2}}\}$$
$$:=C(n,l,\Lambda,\beta) e^{-\lambda_l s}|\phi|_{\mathbf{X}}(I_3+I_4),$$
$$I_3\leq \sqrt{2}(\int_{\beta e^{-\sigma_l s}}^{(\beta+1) e^{-\sigma_l s}} y^{2\alpha+n-4} (1+y^2)e^{-\frac{y^2}{4}}dy+(\beta e^{-\sigma_l s})^{2\alpha+2n-6}\int_{\beta e^{-\sigma_l s}}^{(\beta+1) e^{-\sigma_l s}} y^{2-n} (1+y^2)e^{-\frac{y^2}{4}}dy)^{\frac{1}{2}}$$
$$\leq C(n)\beta^{-\frac{1}{2}}(\beta e^{-\sigma_l s})^{\frac{1}{2}(2\alpha+n-3)},$$
$$I_4\leq C(n)(\beta e^{-\sigma_l s})^{\frac{1}{2}(2\alpha+n-3)},$$
and thus (recalling \eqref{co7})
\begin{equation}
    |f_2(\cdot,s)|_{\mathbf{X}^*} \leq C(n,l,\Lambda,\beta)e^{-\lambda_l(1+\varsigma)s}. \label{es7}
\end{equation}
The estimate of $f_3$ is much simpler:
$$|f_3(\cdot,s)|_{\mathbf{X}^*} \leq|f_3(\cdot,s)|_{\mathbf{H}}\leq C(n,\Lambda) e^{-\lambda_l s}(\int_{\rho e^{\frac{s}{2}}-1}^{\rho e^{\frac{s}{2}}}y^{4\lambda_l+n+2}e^{-\frac{y^2}{4}}dy)^{\frac{1}{2}}$$
\begin{equation}
    \leq C(n,\Lambda) e^{-\lambda_l s}((\rho e^{\frac{s}{2}})^{4\lambda_l+n+2}e^{-\frac{1}{4}(\rho e^{\frac{s}{2}}-1)^2})^{\frac{1}{2}}\leq C(n,\Lambda)e^{-\lambda_l(1+\varsigma)s} \label{es8}
\end{equation}
(if $s_0\gg 1$ (depending on $n,l,\rho$)). Combining \eqref{es6}, \eqref{es7}, \eqref{es8}, we get \eqref{es5}.\\
Now, we first estimate the ``lower frequency" terms in the Fourier expansion of $\tilde{v}(\cdot,s)$. For $0\leq j\leq l-1$, according to \eqref{es10} (with integration by parts) and the condition \eqref{co5}, we have
\begin{equation}
    \begin{cases}
    \partial_s \langle \tilde{v}(\cdot,s),\varphi_j \rangle +\lambda_j \langle \tilde{v}(\cdot,s),\varphi_j \rangle=\langle f(\cdot,s),\varphi_j \rangle,\\
    \langle \tilde{v}(\cdot,\mathring{s}),\varphi_j \rangle=0.
\end{cases} \label{es27}
\end{equation}
Thus, for all $s_0\leq s\leq \mathring{s}$, 
$$|\langle \tilde{v}(\cdot,s),\varphi_j \rangle|=|\int_s^{\mathring{s}} e^{\lambda_j (\xi-s)}\langle f(\cdot,\xi),\varphi_j \rangle d\xi|\leq |\varphi_j|_{\mathbf{X}}\int_s^{\mathring{s}} e^{(\lambda_l-1) (\xi-s)}|f(\cdot,\xi)|_{\mathbf{X}^*}d\xi$$
\begin{equation}
    \leq C(p,q,l,\Lambda,\beta)e^{-(\lambda_l-1)s}\int_s^{+\infty} e^{(\lambda_l-1)\xi-(1+\varsigma)\lambda_l\xi}d\xi=C(p,q,l,\Lambda,\beta)e^{-(1+\varsigma)\lambda_l s}, \label{es11}
\end{equation}
where we used \eqref{es5} and Lemma \ref{lemes1}. In addition, for $0\leq j\leq l-1$, we compute from Lemma \ref{lemco1} and \eqref{co8} that
$$|e^{\lambda_l s_0}\langle \tilde{v}(\cdot,s_0),c_j\varphi_j \rangle-a_j|=|\langle \eta(e^{\sigma_l s}y-\beta)\eta(\rho e^{\frac{s}{2}}-y)(\frac{1}{c_l}\varphi_l(y)+\sum_{m=0}^{l-1} \frac{a_m}{c_m}\varphi_m(y)),c_j\varphi_j \rangle -a_j|$$
$$\leq C(n,l,\beta)e^{-\frac{n-1+2\alpha}{1-\alpha}\lambda_l s_0}\leq C(n,l,\beta)e^{-2\varsigma\lambda_l s_0},$$
and combining with \eqref{es11} we get 
\begin{equation}
    |a_j|\leq |e^{\lambda_l s_0}\langle \tilde{v}(\cdot,s_0),c_j\varphi_j \rangle-a_j|+|e^{\lambda_l s_0}\langle \tilde{v}(\cdot,s_0),c_j\varphi_j \rangle|\leq C(p,q,l,\Lambda,\beta)e^{-\varsigma\lambda_l s_0}, \label{es21}
\end{equation}
which is exactly (\ref{propco1_1}) of Proposition \ref{propco1}.\\
Next we estimate the ``main frequency" term of $\tilde{v}(\cdot,s)$. According to \eqref{es10} (with integration by parts) and Lemma \ref{lemco1}, we have
$$\begin{cases}
    \partial_s (e^{\lambda_l s} \langle \tilde{v}(\cdot,s),\varphi_l \rangle)=e^{\lambda_l s}\langle f(\cdot,s),\varphi_l \rangle,\\
    |e^{\lambda_l s_0} \langle \tilde{v}(\cdot,s_0),c_l\varphi_l \rangle-1|\leq C(n,l,\beta)e^{-2\varsigma\lambda_l s_0}.
\end{cases}$$
Set
$$k=e^{\lambda_l \mathring{s}}\langle \tilde{v}(\cdot,\mathring{s}),c_l\varphi_l \rangle,$$
then for all $s_0\leq s\leq \mathring{s}$,
$$|e^{\lambda_l s} \langle \tilde{v}(\cdot,s),c_l\varphi_l \rangle-k|=|e^{\lambda_l s} \langle \tilde{v}(\cdot,s),c_l\varphi_l \rangle-e^{\lambda_l \mathring{s}}\langle \tilde{v}(\cdot,\mathring{s}),c_l\varphi_l \rangle|=|c_l\int_s^{\mathring{s}}e^{\lambda_l \xi}\langle f(\cdot,\xi),\varphi_l \rangle d\xi|$$
\begin{equation}
    \leq c_l |\varphi_l|_{\mathbf{X}}\int_s^{\mathring{s}} e^{\lambda_l\xi}|f(\cdot,\xi)|_{\mathbf{X}^*}d\xi\leq C(p,q,l,\Lambda,\beta)e^{-\varsigma\lambda_l s}, \label{es12}
\end{equation}
i.e.
\begin{equation}
    |\langle \tilde{v}(\cdot,s),\varphi_l \rangle-\frac{k}{c_l} e^{-\lambda_l s}|\leq C(p,q,l,\Lambda,\beta)e^{-(1+\varsigma)\lambda_l s}, \label{es35}
\end{equation}
where we used \eqref{es5} and Lemma \ref{lemes1} again. Moreover, by Lemma \ref{lemco1} and \eqref{es12},
\begin{equation}
    |k-1|\leq |e^{\lambda_l s_0} \langle \tilde{v}(\cdot,s_0),c_l\varphi_l \rangle-k|+|e^{\lambda_l s_0} \langle \tilde{v}(\cdot,s_0),c_l\varphi_l \rangle-1|\leq C(p,q,l,\Lambda,\beta)e^{-\varsigma\lambda_l s_0}, \label{es22}
\end{equation}
which is exactly \eqref{co6}.\\
The estimate of the ``higher frequency" terms of $\tilde{v}(\cdot,s)$, i.e. $\varphi_j$ terms with $j\geq l+1$, is much more complicated. It is divided into two parts, a ``short time" estimate and a ``long time" estimate; the short time part is achieved by writing down an integral representation of $\tilde{v}$, involving $f$, the initial value $\tilde{v}(\cdot,s_0)$, and the ``heat kernel" of the operator $\mathcal{L}$.
\begin{lemma}
    \label{lemes2}
    For all $s_0<s\leq \mathring{s}$ and $y>0$, we have the integral representation formula:
    $$\tilde{v}(y,s)=\int_0^{+\infty} K(y,z,s-s_0)\tilde{v}(z,s_0)dz+\int_{s_0}^s \int_0^{+\infty} K(y,z,s-\tau)f(z,\tau)dzd\tau,$$
    where
    \begin{equation}
        K(y,z,s)=(\frac{z}{y})^{\frac{n}{2}-1}\sqrt{yz}\frac{e^{\frac{n-1}{4}s}}{2(1-e^{-s})}I_{\frac{n-3}{2}+\alpha}(\frac{e^{-\frac{s}{2}}yz}{2(1-e^{-s})})\exp(-\frac{e^{-s}y^2+z^2}{4(1-e^{-s})}) \label{es20}
    \end{equation}
    is the ``heat kernel" of $\mathcal{L}$, and $I_\nu$ ($\nu\geq 0$) is the modified Bessel function of the first kind, satisfying
    \begin{equation}
        x^2 I_\nu''(x)+x I_\nu'(x)-(x^2+\nu^2)I_\nu(x)=0 \label{es13}
    \end{equation}
    and 
    \begin{equation}
        I_\nu(x)\sim \frac{1}{\Gamma(\nu+1)}(\frac{x}{2})^\nu,\ x\searrow 0. \label{es15}
    \end{equation}
\end{lemma}
\proof
Let
\begin{equation}
    V(x,t)=x^{\frac{n}{2}-1}\sqrt{-t}\tilde{v}(\frac{x}{\sqrt{-t}},-\ln(-t)), \label{es19}
\end{equation}
then $V$ satisfies
$$\partial_t V=\partial_x^2 V+\frac{\frac{1}{4}-\gamma^2}{x^2}V+F(x,t):=\Delta_\gamma V+F(x,t),$$
where
$$\gamma=\frac{n-3}{2}+\alpha=\frac{1}{2}\sqrt{n^2-10n+17},\ F(x,t)=\frac{x^{\frac{n}{2}-1}}{\sqrt{-t}}f(\frac{x}{\sqrt{-t}},-\ln(-t)).$$
One may first let
$$\tilde{u}(x,t)=\sqrt{-t}\tilde{v}(\frac{x}{\sqrt{-t}},-\ln(-t))$$
(see \eqref{ad6}), to get 
$$\partial_t \tilde{u}=\partial_x^2 \tilde{u}+\frac{n-2}{x}\partial_x \tilde{u}+\frac{n-2}{x^2}\tilde{u}+\frac{1}{\sqrt{-t}}f(\frac{x}{\sqrt{-t}},-\ln(-t)),$$
and then let $V=x^{\frac{n}{2}-1}\tilde{u}$ to eliminate the $\partial_x \tilde{u}$ term.\\
The heat kernel of the operator $\Delta_\gamma$ is
$$B(x,w,t)=\frac{\sqrt{xw}}{2t}e^{-\frac{x^2+w^2}{4t}}I_\gamma(\frac{xw}{2t}),\ x,w,t>0,$$
which is known as the ``Bessel heat kernel". Actually, this kernel can be derived with the help of Hankel transform, see \cite{bd}, Section 1. (To compute the inverse transform, one may use Weber's formula, see \cite{wat}, 13.31, Formula (1).)\\
Now, $V:(0,+\infty)\times [t_0,\mathring{t}]\rightarrow \R$ is smooth and compactly supported. We claim that for all $t_0<t\leq \mathring{t}$ and $x>0$,
\begin{equation}
    V(x,t)=\int_0^{+\infty} B(x,w,t-t_0)V(w,t_0)dw+\int_{t_0}^t \int_0^{+\infty} B(x,w,t-\xi)F(w,\xi)dwd\xi. \label{es14}
\end{equation}
By direct computation using the equation \eqref{es13}, for $x,w,t>0$,
$$\partial_t B(x,w,t)=(\partial_w^2+\frac{\frac{1}{4}-\gamma^2}{w^2})B(x,w,t)=\Delta_{\gamma,w}B(x,w,t).$$
Differentiating the right-hand side of \eqref{es14} (with $t_1\geq t_0$ in the place of $t_0$) w.r.t. $t_1$ yields
$$-\int_0^{+\infty} \partial_t B(x,w,t-t_1)V(w,t_1)dw+\int_0^{+\infty}  B(x,w,t-t_1)\partial_t V(w,t_1)dw-\int_0^{+\infty} B(x,w,t-t_1)F(w,t_1)dw$$
$$=-\int_0^{+\infty} \Delta_{\gamma,w} B(x,w,t-t_1)V(w,t_1)dw+\int_0^{+\infty}  B(x,w,t-t_1)\partial_t V(w,t_1)dw-\int_0^{+\infty} B(x,w,t-t_1)F(w,t_1)dw$$
$$=-\int_0^{+\infty} B(x,w,t-t_1)\Delta_{\gamma} V(w,t_1)dw+\int_0^{+\infty}  B(x,w,t-t_1)\partial_t V(w,t_1)dw-\int_0^{+\infty} B(x,w,t-t_1)F(w,t_1)dw=0,$$
i.e. the right-hand side of \eqref{es14} is independent of $t_1$, when $t$ is fixed. In order to compute the limit when $t_1\nearrow t$, we use variable changing:
$$r=\frac{w-x}{2\sqrt{t-t_1}},$$
$$\int_0^{+\infty} B(x,w,t-t_1)V(w,t_1)dw=\int_{-\frac{x}{2\sqrt{t-t_1}}}^{+\infty} B(x,2\sqrt{t-t_1}r+x,t-t_1)V(2\sqrt{t-t_1}r+x,t_1)2\sqrt{t-t_1}dr$$
\begin{equation}
    =\int_{-\frac{x}{2\sqrt{t-t_1}}}^{+\infty} \frac{\sqrt{x(2\sqrt{t-t_1}r+x)}}{2(t-t_1)}e^{-\frac{x^2+(2\sqrt{t-t_1}r+x)^2}{4(t-t_1)}}I_\gamma(\frac{x(2\sqrt{t-t_1}r+x)}{2(t-t_1)})V(2\sqrt{t-t_1}r+x,t_1)2\sqrt{t-t_1}dr. \label{es18}
\end{equation}
According to Formula (2.13) of \cite{bd},
\begin{equation}
    \lim_{x\rightarrow +\infty} I_\gamma(x) (\frac{e^x}{\sqrt{2\pi x}})^{-1}=1, \label{es16}
\end{equation}
thus by \eqref{es15}, there exists $C(\gamma)>0$ s.t.
\begin{equation}
    |I_\gamma(x)|\leq C(\gamma)\frac{x^\gamma e^x}{(x+1)^{\gamma+\frac{1}{2}}},\ x>0. \label{es17}
\end{equation}
By \eqref{es16}, as $t_1\nearrow t$, the integrand of \eqref{es18} tends to $\frac{1}{\sqrt{\pi}}e^{-r^2}V(x,t)$, and is dominated by
$$C(\gamma) e^{-r^2}\sup_{w>0,\ t_0\leq t_1\leq t}|V(w,t_1)|,$$
and the integral domain tends to $(-\infty,+\infty)$. Therefore, by Dominated Convergence Theorem, 
$$\lim_{t_1\nearrow t}\int_0^{+\infty} B(x,w,t-t_1)V(w,t_1)dw=V(x,t).$$
Similarly, the integral
$$\int_0^{+\infty} B(x,w,t-\xi)F(w,\xi)dw$$
is bounded for $t_0\leq \xi<t$, and 
$$\lim_{t_1\nearrow t}\int_{t_1}^t \int_0^{+\infty} B(x,w,t-\xi)F(w,\xi)dwd\xi=0.$$
Thus, \eqref{es14} is proved.\\
Applying the inverse transform of \eqref{es19}, we have
$$\tilde{v}(y,s)=y^{1-\frac{n}{2}}e^{\frac{n}{4}s}V(e^{-\frac{s}{2}}y,-e^{-s})$$
$$=y^{1-\frac{n}{2}}e^{\frac{n}{4}s}(\int_0^{+\infty} B(e^{-\frac{s}{2}}y,w,-e^{-s}+e^{-s_0})V(w,-e^{-s_0})dw+\int_{-e^{-s_0}}^{-e^{-s}} \int_0^{+\infty} B(e^{-\frac{s}{2}}y,w,-e^{-s}-\xi)F(w,\xi)dwd\xi)$$
$$=y^{1-\frac{n}{2}}e^{\frac{n}{4}s}(\int_0^{+\infty} B(e^{-\frac{s}{2}}y,e^{-\frac{s_0}{2}}z,-e^{-s}+e^{-s_0})V(e^{-\frac{s_0}{2}}z,-e^{-s_0})e^{-\frac{s_0}{2}}dz$$
$$+\int_{s_0}^s \int_0^{+\infty} B(e^{-\frac{s}{2}}y,e^{-\frac{\tau}{2}}z,-e^{-s}+e^{-\tau})F(e^{-\frac{\tau}{2}}z,-e^{-\tau})e^{-\frac{3}{2}\tau}dzd\tau)$$
$$=\int_0^{+\infty} (\frac{z}{y})^{\frac{n}{2}-1}e^{\frac{n}{4}(s-s_0)}e^{-\frac{s_0}{2}} B(e^{-\frac{s}{2}}y,e^{-\frac{s_0}{2}}z,-e^{-s}+e^{-s_0})\tilde{v}(z,s_0)dz$$
$$+\int_{s_0}^s \int_0^{+\infty} (\frac{z}{y})^{\frac{n}{2}-1}e^{\frac{n}{4}(s-\tau)}e^{-\frac{\tau}{2}} B(e^{-\frac{s}{2}}y,e^{-\frac{\tau}{2}}z,-e^{-s}+e^{-\tau})f(z,\tau)dzd\tau$$
$$=\int_0^{+\infty} K(y,z,s-s_0)\tilde{v}(z,s_0)dz+\int_{s_0}^s \int_0^{+\infty} K(y,z,s-\tau)f(z,\tau)dzd\tau.$$

\QED

Parallel to Lemma \ref{lemes2}, we have the integral formula of eigenfunctions $\varphi_j$:
\begin{lemma}
    \label{lemes3}
    For $j\geq 0$, $y>0$, $s>s_0$,
    $$\varphi_j(y)=e^{\lambda_j(s-s_0)}\int_0^{+\infty} K(y,z,s-s_0)\varphi_j(z)dz,$$
    where $K(y,z,s)$ is defined in \eqref{es20}.
\end{lemma}
\proof
The proof is similar to that of Lemma \ref{lemes2}, with $\tilde{v}(y,s)$ replaced by $e^{-\lambda_j s}\varphi_j(y)$ and $V(x,t)$ replaced by $V_j(x,t):=x^{\frac{n}{2}-1}(-t)^{\lambda_j+\frac{1}{2}}\varphi_j(\frac{x}{\sqrt{-t}})$, except that $V_j$ here is not compactly supported. The transform \eqref{es19} is still used here. Note that $(\partial_s-\mathcal{L})(e^{-\lambda_j s}\varphi_j(y))=0$. It suffices to check if the computation above is still valid here.\\
According to \eqref{ad11}, for $i=0,1,2$,
$$x^i\partial_x^i V_j(x,t)=\begin{cases}
    \begin{array}{ll}
      O(x^{\frac{n}{2}-1+\alpha}),   & x\searrow 0, \\
      O(x^{\frac{n}{2}+2\lambda_j}),   & x\rightarrow +\infty,
    \end{array}
\end{cases}$$
$$\partial_t V_j(x,t)=\begin{cases}
    \begin{array}{ll}
      O(x^{\frac{n}{2}-1+\alpha}),   & x\searrow 0, \\
      O(x^{\frac{n}{2}+2\lambda_j-2}),   & x\rightarrow +\infty,
    \end{array}
\end{cases}$$
uniformly when $t$ is bounded. Also, by the derivative formula of $I_\gamma$ (\cite{bd}, Formula (2.14)), for $i=0,1,2$,
$$x^i\partial_x^i I_\gamma(x)=O(x^\gamma),\ x\searrow 0,$$
$$\partial_x^i I_\gamma(x)=O(\frac{e^x}{\sqrt{x}}),\ x\rightarrow +\infty,$$
and thus when $x>0$ is fixed,
$$w^i\partial_w^i B(x,w,t)=O(w^{\frac{n}{2}-1+\alpha}),\ w\searrow 0,$$
$$\partial_w^i B(x,w,t)=O(e^{-\frac{(w-x)^2}{4t}} w^i),\ w\rightarrow +\infty,$$
$$\partial_t B(x,w,t)=\begin{cases}
    \begin{array}{ll}
      O(w^{\frac{n}{2}-1+\alpha}),   & w\searrow 0, \\
     O(e^{-\frac{(w-x)^2}{4t}} w^2),   & w\rightarrow +\infty,
    \end{array}
\end{cases}$$
uniformly when $t$ is bounded and away from 0. Therefore the integral 
$$\int_0^{+\infty} B(x,w,t-t_1)V_j(w,t_1)dw$$
converges absolutely ($t_0\leq t_1<t$), and all the steps below are legal, including differentiating under the integral sign and integrating by parts:
$$\partial_{t_1} \int_0^{+\infty} B(x,w,t-t_1)V_j(w,t_1)dw$$
$$=-\int_0^{+\infty} \partial_t B(x,w,t-t_1)V_j(w,t_1)dw+\int_0^{+\infty}  B(x,w,t-t_1)\partial_t V_j(w,t_1)dw$$
$$=-\int_0^{+\infty} \Delta_{\gamma,w} B(x,w,t-t_1)V_j(w,t_1)dw+\int_0^{+\infty}  B(x,w,t-t_1)\partial_t V_j(w,t_1)dw$$
$$=-\int_0^{+\infty} B(x,w,t-t_1)\Delta_{\gamma} V_j(w,t_1)dw+\int_0^{+\infty}  B(x,w,t-t_1)\partial_t V_j(w,t_1)dw=0.$$
In the dominated convergence argument, the integrand
$$|B(x,2\sqrt{t-t_1}r+x,t-t_1)V_j(2\sqrt{t-t_1}r+x,t_1)2\sqrt{t-t_1}|\leq Ce^{-r^2}(1+|r|^{\frac{n}{2}+2\lambda_j}),$$
where $C$ is independent of $t_1$ if $t_0\leq t_1<t$.

\QED

Next, we estimate the evolution of the ``non-homogeneous" term $f$:
\begin{lemma}
    \label{lemes4}
    Let
    $$S(y,s)=\int_{s_0}^s \int_0^{+\infty} K(y,z,s-\tau)f(z,\tau)dzd\tau.$$
    If \eqref{es1}, \eqref{es2}, \eqref{es3} hold, and $s_0\gg 1$ (depending on $n,l,\rho,\beta,R,\vartheta$), then for any $s_0<s\leq \min\{\mathring{s},s_0+1\}$, $\frac{1}{2} e^{-\vartheta \sigma_l s}\leq y\leq 2R$,
    $$|S(y,s)|\leq C(p,q,l,\Lambda,\beta,R)e^{-(1+\tilde{\kappa})\lambda_l s_0}y^\alpha.$$
\end{lemma}
\proof
According to \eqref{es17}, for $y,z>0$, $0<s\leq 1$, there exists a constant $C(n)$ s.t.
\begin{equation}
    |K(y,z,s)|\leq C(n)y^\alpha s^{-(\frac{n-1}{2}+\alpha)} z^{n-2+\alpha}(1+c\frac{yz}{s})^{-(\frac{n}{2}-1+\alpha)}\exp(-\frac{(e^{-\frac{s}{2}}y-z)^2}{4s}), \label{es24}
\end{equation}
where $c=\frac{1}{2}e^{-\frac{1}{2}}$. Write
$$H(y,z,s):=(1+c\frac{yz}{s})^{-(\frac{n}{2}-1+\alpha)}\exp(-\frac{(e^{-\frac{s}{2}}y-z)^2}{4s}).$$
By \eqref{es10}, 
$$S(y,s)=\sum_{i=1}^3 S_i(y,s),\ S_i(y,s)=\int_{s_0}^s \int_0^{+\infty} K(y,z,s-\tau)f_i(z,\tau)dzd\tau.$$
Let's start from the estimate of $S_1$. By \eqref{es1} and \eqref{es24}, 
$$|S_1(y,s)|\leq C(p,q,\Lambda) e^{-2\lambda_l s_0}y^\alpha\int_{s_0}^s (s-\tau)^{-(\frac{n-1}{2}+\alpha)}\int_{\beta e^{-\sigma_l \tau}}^{\rho e^{\frac{\tau}{2}}}z^{n-2+\alpha}(z^{2\alpha-3}+z^{4\lambda_l-1})H(y,z,s-\tau) dzd\tau$$
$$\leq C(p,q,\Lambda) e^{-2\lambda_l s_0}y^\alpha\int_{s_0}^s (s-\tau)^{-(\frac{n-1}{2}+\alpha)}\int_{\check{c}\beta e^{-\sigma_l s_0}}^{+\infty}z^{n-2+\alpha}(z^{2\alpha-3}+z^{4\lambda_l-1})H(y,z,s-\tau) dzd\tau$$
($\check{c}=e^{-\sigma_l}$)
$$:=C(p,q,\Lambda) e^{-2\lambda_l s_0}y^\alpha(S_{1,1}+S_{1,2}),$$
$$S_{1,j}(y,s)=\int_{s_0}^s (s-\tau)^{-(\frac{n-1}{2}+\alpha)}\int_{\check{c}\beta e^{-\sigma_l s_0}}^{+\infty}z^{n-2+\alpha+b_j}(1+c\frac{yz}{s-\tau})^{-(\frac{n}{2}-1+\alpha)}\exp(-\frac{(e^{-\frac{s-\tau}{2}}y-z)^2}{4(s-\tau)}) dzd\tau,\ j=1,2,$$
$$b_1=2\alpha-3,\ b_2=4\lambda_l-1.$$
We introduce new variables
\begin{equation}
    w=\frac{z}{\sqrt{s-\tau}},\ \xi=\frac{y^2}{s-\tau}, \label{es25}
\end{equation}
to get
$$S_{1,j}(y,s)$$
$$=y^{2+b_j-\alpha}\int_{\frac{y^2}{s-s_0}}^{+\infty}\xi^{\frac{1}{2}\alpha-\frac{1}{2}b_j-2} \int_{\check{c}\beta e^{-\sigma_l s_0}\frac{\sqrt{\xi}}{y}}^{+\infty} w^{n-2+\alpha+b_j}(1+cw\sqrt{\xi})^{-(\frac{n}{2}-1+\alpha)}\exp(-\frac{1}{4}(e^{-\frac{y^2}{2\xi}}\sqrt{\xi}-w)^2)dwd\xi$$
$$\leq y^{2+b_j-\alpha}\int_{y^2}^{+\infty}\xi^{\frac{1}{2}\alpha-\frac{1}{2}b_j-2} \int_{\check{c}\beta e^{-\sigma_l s_0}\frac{\sqrt{\xi}}{y}}^{+\infty} w^{n-2+\alpha+b_j}(1+cw\sqrt{\xi})^{-(\frac{n}{2}-1+\alpha)}\exp(-\frac{1}{4}(e^{-\frac{y^2}{2\xi}}\sqrt{\xi}-w)^2)dwd\xi$$
$$=y^{2+b_j-\alpha}\int_{y^2}^{+\infty}\xi^{\frac{1}{2}\alpha-\frac{1}{2}b_j-2} (\int_{\check{c}\beta e^{-\sigma_l s_0}\frac{\sqrt{\xi}}{y}}^{\frac{1}{2}\sqrt{\xi}}...dw+\int_{\frac{1}{2}\sqrt{\xi}}^{2\sqrt{\xi}}...dw+\int_{2\sqrt{\xi}}^{+\infty}...dw)d\xi$$
$$:=S_{1,j,1}+S_{1,j,2}+S_{1,j,3}.$$
In the rest argument, we denote by $c$ a positive universal constant, which may change from line to line.
$$S_{1,j,1}\leq y^{2+b_j-\alpha}\int_0^{+\infty}\xi^{\frac{1}{2}\alpha-\frac{1}{2}b_j-2} e^{-c\xi}\int_{\check{c}\beta e^{-\sigma_l s_0}\frac{\sqrt{\xi}}{y}}^{\frac{1}{2}\sqrt{\xi}}w^{n-2+\alpha+b_j}(1+cw^2)^{-(\frac{n}{2}-1+\alpha)}dwd\xi$$
(noting $\frac{y^2}{\xi}=s-\tau\leq 1$).
$$S_{1,1,1}\leq y^{\alpha-1}\int_0^{+\infty}\xi^{-\frac{1}{2}\alpha-\frac{1}{2}} e^{-c\xi}\int_{\check{c}\beta e^{-\sigma_l s_0}\frac{\sqrt{\xi}}{y}}^{+\infty}w^{n-5+3\alpha}(1+cw^2)^{-(\frac{n}{2}-1+\alpha)}dwd\xi.$$
If $n-5+3\alpha>-1$, then
$$S_{1,1,1}\leq y^{\alpha-1}\int_0^{+\infty}\xi^{-\frac{1}{2}\alpha-\frac{1}{2}} e^{-c\xi}\int_0^{+\infty}w^{n-5+3\alpha}(1+cw^2)^{-(\frac{n}{2}-1+\alpha)}dwd\xi$$
$$=C(n)y^{\alpha-1}\leq  C(n,l)e^{\vartheta \lambda_l s_0}\leq C(n,l)e^{(1-\tilde{\kappa})\lambda_l s_0}\ (\text{since}\ y\geq \frac{1}{2} e^{-\vartheta \sigma_l s}\geq \frac{1}{2}e^{-\sigma_l}e^{-\vartheta \sigma_l s_0}).$$
If $n-5+3\alpha<-1$, then
$$S_{1,1,1}\leq y^{\alpha-1}\int_0^{+\infty}\xi^{-\frac{1}{2}\alpha-\frac{1}{2}} e^{-c\xi}\int_{\check{c}\beta e^{-\sigma_l s_0}\frac{\sqrt{\xi}}{y}}^{+\infty}w^{n-5+3\alpha}dwd\xi$$
$$=C(n,\beta)y^{\alpha-1}(\frac{y}{e^{-\sigma_l s_0}})^{-(n-4+3\alpha)}\int_0^{+\infty} \xi^{\frac{n-5}{2}+\alpha}e^{-c\xi}d\xi=C(n,\beta)e^{\lambda_l s_0}(\frac{y}{e^{-\sigma_l s_0}})^{-(n-3+2\alpha)}$$
$$\leq C(n,l,\beta)e^{\lambda_l s_0}e^{-(1-\vartheta)\frac{n-3+2\alpha}{1-\alpha}\lambda_l s_0}\leq C(n,l,\beta)e^{(1-\tilde{\kappa})\lambda_l s_0}.$$
$$S_{1,2,1}\leq y^{4\lambda_l+1-\alpha}\int_0^{+\infty}\xi^{\frac{1}{2}\alpha-2\lambda_l-\frac{3}{2}} e^{-c\xi}\int_0^{\frac{1}{2}\sqrt{\xi}}w^{n-3+\alpha+4\lambda_l}dwd\xi$$
$$=C(n,l)y^{4\lambda_l+1-\alpha}\int_0^{+\infty}\xi^{\frac{n-5}{2}+\alpha} e^{-c\xi}d\xi=C(n,l)y^{4\lambda_l+1-\alpha}\leq C(n,l,R).$$
$$S_{1,j,2}\leq C(n,b_j) y^{2+b_j-\alpha}\int_0^{+\infty}\xi^{\frac{1}{2}\alpha-\frac{1}{2}b_j-2}\xi^{\frac{1}{2}(n-2+\alpha+b_j)}(1+c\xi)^{-(\frac{n}{2}-1+\alpha)} \sqrt{\xi}d\xi$$
$$=C(n,b_j) y^{2+b_j-\alpha}\int_0^{+\infty}\xi^{\frac{n-5}{2}+\alpha}(1+c\xi)^{-(\frac{n}{2}-1+\alpha)}d\xi=C(n,b_j) y^{2+b_j-\alpha},$$
$$S_{1,1,2}\leq C(n) y^{\alpha-1}\leq  C(n,l)e^{\vartheta \lambda_l s_0}\leq C(n,l)e^{(1-\tilde{\kappa})\lambda_l s_0},$$
$$S_{1,2,2}\leq C(n,l) y^{4\lambda_l+1-\alpha}\leq  C(n,l,R).$$
$$S_{1,j,3}\leq y^{2+b_j-\alpha}\int_{y^2}^{+\infty}\xi^{\frac{1}{2}\alpha-\frac{1}{2}b_j-2} e^{-c\xi} \int_{2\sqrt{\xi}}^{+\infty} w^{n-2+\alpha+b_j}e^{-cw^2}dwd\xi$$
(since $w\geq 2\sqrt{\xi}$, we have $\exp(-\frac{1}{4}(e^{-\frac{y^2}{2\xi}}\sqrt{\xi}-w)^2)\leq e^{-\frac{1}{16}w^2}\leq e^{-\frac{1}{32}w^2}e^{-\frac{1}{8}\xi}$).
$$S_{1,1,3}\leq y^{\alpha-1}\int_{y^2}^{+\infty}\xi^{-\frac{1}{2}\alpha-\frac{1}{2}} e^{-c\xi} \int_{2\sqrt{\xi}}^{+\infty} w^{n-5+3\alpha}e^{-cw^2}dwd\xi$$
If $n-5+3\alpha>-1$, then
$$S_{1,1,3}\leq y^{\alpha-1}\int_0^{+\infty}\xi^{-\frac{1}{2}\alpha-\frac{1}{2}} e^{-c\xi} \int_0^{+\infty} w^{n-5+3\alpha}e^{-cw^2}dwd\xi$$
$$=C(n)y^{\alpha-1}\leq  C(n,l)e^{\vartheta \lambda_l s_0}\leq C(n,l)e^{(1-\tilde{\kappa})\lambda_l s_0}.$$
If $n-5+3\alpha<-1$, then
$$S_{1,1,3}\leq y^{\alpha-1}\int_0^{+\infty}\xi^{-\frac{1}{2}\alpha-\frac{1}{2}} e^{-c\xi} \int_{2\sqrt{\xi}}^{+\infty} w^{n-5+3\alpha}dwd\xi=C(n)y^{\alpha-1}\int_0^{+\infty}\xi^{\frac{n-5}{2}+\alpha} e^{-c\xi}d\xi$$
$$=C(n)y^{\alpha-1}\leq C(n,l)e^{(1-\tilde{\kappa})\lambda_l s_0}.$$
$$S_{1,2,3}\leq y^{4\lambda_l+1-\alpha}\int_{y^2}^{+\infty}\xi^{\frac{1}{2}\alpha-2\lambda_l-\frac{3}{2}} \int_0^{+\infty} w^{n-3+\alpha+4\lambda_l}e^{-cw^2}dwd\xi$$
$$=C(n,l)y^{4\lambda_l+1-\alpha}\int_{y^2}^{+\infty}\xi^{\frac{1}{2}\alpha-2\lambda_l-\frac{3}{2}}d\xi=C(n,l).$$
Now we have obtained
$$|S_1(y,s)|\leq C(p,q,l,\Lambda,\beta,R)e^{-(1+\tilde{\kappa})\lambda_l s_0}y^\alpha.$$
Next we estimate $S_2$. By \eqref{es2} and \eqref{es24}, 
$$|S_2(y,s)|\leq C(n,l,\Lambda,\beta) e^{-\lambda_l s_0}y^\alpha$$
$$\times\int_{s_0}^s (s-\tau)^{-(\frac{n-1}{2}+\alpha)}\int_0^{(\beta+1) e^{-\sigma_l s_0}}z^{n-4+2\alpha}(1+c\frac{yz}{s-\tau})^{-(\frac{n}{2}-1+\alpha)}\exp(-\frac{(e^{-\frac{s-\tau}{2}}y-z)^2}{4(s-\tau)}) dzd\tau$$
$$=C(n,l,\Lambda,\beta) e^{-\lambda_l s_0}y^\alpha$$
$$\times \int_{\frac{y^2}{s-s_0}}^{+\infty}\xi^{-1} \int_0^{(\beta+1) e^{-\sigma_l s_0}\frac{\sqrt{\xi}}{y}} w^{n-4+2\alpha}(1+cw\sqrt{\xi})^{-(\frac{n}{2}-1+\alpha)}\exp(-\frac{1}{4}(e^{-\frac{y^2}{2\xi}}\sqrt{\xi}-w)^2)dwd\xi$$
(using \eqref{es25} again)
$$\leq C(n,l,\Lambda,\beta) e^{-\lambda_l s_0}y^\alpha\int_0^{+\infty}\xi^{-1}e^{-c\xi} \int_0^{(\beta+1) e^{-\sigma_l s_0}\frac{\sqrt{\xi}}{y}} w^{n-4+2\alpha}dwd\xi$$
($\frac{y^2}{\xi}=s-\tau\leq 1$, and $(\beta+1) e^{-\sigma_l s_0}\frac{\sqrt{\xi}}{y}\leq \frac{1}{2}\sqrt{\xi}$ provided $s_0\gg 1$ (depending on $n,l,\beta,\vartheta$))
$$=C(n,l,\Lambda,\beta) e^{-\lambda_l s_0}y^\alpha(\frac{y}{e^{-\sigma_l s_0}})^{-(n-3+2\alpha)}\int_0^{+\infty}\xi^{\frac{n-5}{2}+\alpha}e^{-c\xi}d\xi$$
$$\leq C(n,l,\Lambda,\beta)e^{-\lambda_l s_0}y^\alpha e^{-(1-\vartheta)\frac{n-3+2\alpha}{1-\alpha}\lambda_l s_0}\leq C(n,l,\Lambda,\beta)e^{-(1+\tilde{\kappa})\lambda_l s_0}y^\alpha.$$
Finally we estimate $S_3$. By \eqref{es3} and \eqref{es24}, 
$$|S_3(y,s)|\leq C(n,\Lambda) e^{-\lambda_l s_0}y^\alpha$$
$$\times\int_{s_0}^s (s-\tau)^{-(\frac{n-1}{2}+\alpha)}\int_{\rho e^{\frac{s_0}{2}}-1}^{+\infty}z^{n+2\alpha+2\lambda_l}(1+c\frac{yz}{s-\tau})^{-(\frac{n}{2}-1+\alpha)}\exp(-\frac{(e^{-\frac{s-\tau}{2}}y-z)^2}{4(s-\tau)}) dzd\tau\ (\hat{c}=e^{\frac{1}{2}})$$
$$=C(n,\Lambda) e^{-\lambda_l s_0}y^{2\lambda_l+4}$$
$$\times \int_{\frac{y^2}{s-s_0}}^{+\infty}\xi^{\frac{1}{2}\alpha-\lambda_l-3} \int_{(\rho e^{\frac{s_0}{2}}-1)\frac{\sqrt{\xi}}{y}}^{+\infty} w^{n+2\alpha+2\lambda_l}(1+cw\sqrt{\xi})^{-(\frac{n}{2}-1+\alpha)}\exp(-\frac{1}{4}(e^{-\frac{y^2}{2\xi}}\sqrt{\xi}-w)^2)dwd\xi$$
(using \eqref{es25} again)
$$\leq C(n,\Lambda) e^{-\lambda_l s_0}y^{2\lambda_l+4}\int_{y^2}^{+\infty}\xi^{\frac{1}{2}\alpha-\lambda_l-3}\int_{\rho e^{\frac{s_0}{2}}-1}^{+\infty} w^{n+2\alpha+2\lambda_l} e^{-cw^2}dwd\xi$$
($(\rho e^{\frac{s_0}{2}}-1)\frac{\sqrt{\xi}}{y}\geq 2\sqrt{\xi}$ provided $s_0\gg 1$ (depending on $\rho,R$))
$$\leq C(n,l,\Lambda) e^{-\lambda_l s_0}y^{2\lambda_l+4}e^{-\frac{1}{2}c(\rho e^{\frac{s_0}{2}}-1)^2}\int_{y^2}^{+\infty}\xi^{\frac{1}{2}\alpha-\lambda_l-3}d\xi$$
$$=C(n,l,\Lambda) e^{-\lambda_l s_0}y^\alpha e^{-\frac{1}{2} c(\rho e^{\frac{s_0}{2}}-1)^2}\leq C(n,l,\Lambda) e^{-(1+\tilde{\kappa})\lambda_l s_0}y^\alpha,$$
if $s_0\gg 1$ (depending on $n,l,\rho,\vartheta$).\\
The proof of this lemma is finished by combining the estimates of $S_1,S_2,S_3$.

\QED

The following lemma deals with the evolution of the initial value:
\begin{lemma}
    \label{lemes5}
    Let
    $$g(y)=\tilde{v}(y,s_0)-\langle \tilde{v}(\cdot,s_0),\varphi_l \rangle \varphi_l(y).$$
    $$T(y,s)=\int_0^{+\infty} K(y,z,s-s_0)g(z)dz.$$
    If \eqref{es21} holds, and $s_0\gg 1$ (depending on $n,l,\rho,\beta,R,\vartheta$), then for any $s_0<s\leq \min\{\mathring{s},s_0+1\}$, $\frac{1}{2}e^{-\vartheta \sigma_l s}\leq y\leq 2R$,
    $$|T(y,s)|\leq C(p,q,l,\Lambda,\beta,R)e^{-(1+\tilde{\kappa})\lambda_ls_0} y^\alpha.$$
\end{lemma}
\proof
For the function $g(y)$, there holds:
\begin{equation}
    |g(y)|\leq \begin{cases}
    \begin{array}{ll}
      C(p,q,l,\Lambda,\beta) e^{-(1+\varsigma)\lambda_l s_0}(y^\alpha+y^{2\lambda_l+1}),   & (\beta+1)e^{-\sigma_l s_0}\leq y \leq \rho e^{\frac{s_0}{2}}-1, \\
       C(n,l)e^{-\lambda_l s_0}(y^\alpha+y^{2\lambda_l+1}),  & \text{otherwise}.
    \end{array}
\end{cases} \label{es26}
\end{equation}
Actually, by Lemma \ref{lemco1}, \eqref{co8}, \eqref{es21}, and \eqref{ad11}, for $(\beta+1)e^{-\sigma_l s_0}\leq y \leq \rho e^{\frac{s_0}{2}}-1$,
$$|g(y)|\leq e^{-\lambda_l s_0}\sum_{j=0}^{l-1} \frac{|a_j|}{c_j}|\varphi_j(y)|+|\frac{1}{c_l}e^{-\lambda_l s_0}-\langle \tilde{v}(\cdot,s_0),\varphi_l \rangle|\ |\varphi_l(y)|$$
$$\leq C(p,q,l,\Lambda,\beta) e^{-(1+\varsigma)\lambda_l s_0}\sum_{j=0}^{l-1}|\varphi_j(y)|+C(n,l,\beta) e^{-(1+2\varsigma)\lambda_l s_0}|\varphi_l(y)|$$
$$\leq  C(p,q,l,\Lambda,\beta) e^{-(1+\varsigma)\lambda_l s_0}(y^\alpha+y^{2\lambda_l+1}),$$
and similarly, for other $y$,
$$|g(y)|\leq e^{-\lambda_l s_0}\sum_{j=0}^{l-1} \frac{|a_j|}{c_j}|\varphi_j(y)|+(\frac{1}{c_l}e^{-\lambda_l s_0}+|\langle \tilde{v}(\cdot,s_0),\varphi_l \rangle|)|\varphi_l(y)|\leq C(n,l)e^{-\lambda_l s_0}(y^\alpha+y^{2\lambda_l+1}).$$
Thus, we can write
$$T(y,s)=(\int_0^{(\beta+1)e^{-\sigma_l s_0}}+\int_{(\beta+1)e^{-\sigma_l s_0}}^{\rho e^{\frac{s_0}{2}}-1}+\int_{\rho e^{\frac{s_0}{2}}-1}^{+\infty}) K(y,z,s-s_0)g(z)dz:=T_1(y,s)+T_2(y,s)+T_3(y,s).$$
Then we estimate $T_j$ ($1\leq j\leq 3$) one by one. If $s_0\gg 1$ (depending on $n,l,\rho,\beta$), then $(\beta+1)e^{-\sigma_l s_0}\leq 1 \leq \rho e^{\frac{s_0}{2}}-1$. By \eqref{es26} and \eqref{es24},
$$|T_1(y,s)|\leq C(n,l)e^{-\lambda_l s_0}y^\alpha (s-s_0)^{-(\frac{n-1}{2}+\alpha)}$$
$$\times \int_0^{(\beta+1)e^{-\sigma_l s_0}}  z^{n-2+2\alpha}(1+c\frac{yz}{s-s_0})^{-(\frac{n}{2}-1+\alpha)}\exp(-\frac{(e^{-\frac{s-s_0}{2}}y-z)^2}{4(s-s_0)})dz.$$
If $s_0\gg 1$ (depending on $n,l,\beta,\vartheta$), then
$$e^{-\frac{s-s_0}{2}}y\geq \frac{1}{2}e^{-(\frac{1}{2}+\sigma_l)}e^{-\vartheta \sigma_l s_0}\geq 2(\beta+1)e^{-\sigma_l s_0}e^{\frac{2}{3}(1-\vartheta)\sigma_l s_0}\geq 2e^{\frac{2}{3}(1-\vartheta)\sigma_l s_0}z,$$
and therefore
$$|T_1(y,s)|\leq C(n,l)e^{-\lambda_l s_0}y^\alpha (s-s_0)^{-(\frac{n-1}{2}+\alpha)}$$
$$\times \int_0^{(\beta+1)e^{-\sigma_l s_0}}  z^{n-2+2\alpha}(1+c\frac{e^{\frac{2}{3}(1-\vartheta)\sigma_l s_0}z^2}{s-s_0})^{-(\frac{n}{2}-1+\alpha)}\exp(-\frac{e^{\frac{2}{3}(1-\vartheta)\sigma_l s_0}z^2}{4(s-s_0)})dz$$
By introducing a new variable
$$w=e^{\frac{1}{3}(1-\vartheta)\sigma_l s_0}\frac{z}{\sqrt{s-s_0}},$$
we get
$$|T_1(y,s)|\leq C(n,l)e^{-\lambda_l s_0}e^{-(1-\vartheta)\frac{n-1+2\alpha}{3(1-\alpha)}\lambda_l s_0}y^\alpha \int_0^{+\infty}  w^{n-2+2\alpha}(1+cw^2)^{-(\frac{n}{2}-1+\alpha)}e^{-\frac{w^2}{4}}dw$$
$$\leq C(n,l)e^{-(1+\tilde{\kappa})\lambda_l s_0}y^\alpha.$$
Similarly,
$$|T_2(y,s)|\leq C(p,q,l,\Lambda,\beta) e^{-(1+\varsigma)\lambda_l s_0}y^\alpha (s-s_0)^{-(\frac{n-1}{2}+\alpha)}$$
$$\times \int_0^{+\infty}  z^{n-2+2\alpha}(1+z^{2\lambda_l+1-\alpha})(1+c\frac{yz}{s-s_0})^{-(\frac{n}{2}-1+\alpha)}\exp(-\frac{(e^{-\frac{s-s_0}{2}}y-z)^2}{4(s-s_0)})dz.$$
$$=C(p,q,l,\Lambda,\beta) e^{-(1+\varsigma)\lambda_l s_0}y^\alpha (s-s_0)^{-(\frac{n-1}{2}+\alpha)}(\int_0^{2y}...dz+\int_{2y}^{+\infty} ...dz)$$
$$:=C(p,q,l,\Lambda,\beta) e^{-(1+\varsigma)\lambda_l s_0}y^\alpha(T_{2,1}+T_{2,2}),$$
$$T_{2,1}\leq C(n)(s-s_0)^{-(\frac{n-1}{2}+\alpha)}\int_0^{2y}z^{n-2+2\alpha}(1+y^{2\lambda_l+1-\alpha})(\frac{z^2}{s-s_0})^{-(\frac{n}{2}-1+\alpha)}\exp(-\frac{(e^{-\frac{s-s_0}{2}}y-z)^2}{4(s-s_0)})dz$$
$$=C(n)(1+y^{2\lambda_l+1-\alpha})(s-s_0)^{-\frac{1}{2}}\int_0^{2y}\exp(-\frac{(e^{-\frac{s-s_0}{2}}y-z)^2}{4(s-s_0)})dz$$
$$\leq C(n)(1+y^{2\lambda_l+1-\alpha})\int_{-\infty}^{+\infty} e^{-\frac{w^2}{4}}dw\ (w=\frac{z-e^{-\frac{s-s_0}{2}}y}{\sqrt{s-s_0}})$$
$$=C(n)(1+y^{2\lambda_l+1-\alpha}) \leq C(n,l,R).$$
$$T_{2,2}\leq (s-s_0)^{-(\frac{n-1}{2}+\alpha)}\int_0^{+\infty}z^{n-2+2\alpha}(1+z^{2\lambda_l+1-\alpha})e^{-\frac{z^2}{16(s-s_0)}}dz$$
$$=\int_0^{+\infty}w^{n-2+2\alpha}(1+(w\sqrt{s-s_0})^{2\lambda_l+1-\alpha})e^{-\frac{w^2}{16}}dw\ (w=\frac{z}{\sqrt{s-s_0}})$$
$$\leq \int_0^{+\infty}w^{n-2+2\alpha}(1+w^{2\lambda_l+1-\alpha})e^{-\frac{w^2}{16}}dw=C(n,l).$$
Therefore,
$$|T_2(y,s)|\leq C(p,q,l,\Lambda,\beta,R) e^{-(1+\varsigma)\lambda_l s_0}y^\alpha\leq C(p,q,l,\Lambda,\beta,R) e^{-(1+\tilde{\kappa})\lambda_l s_0}y^\alpha.$$
Also,
$$|T_3(y,s)|\leq C(n,l)e^{-\lambda_l s_0}y^\alpha (s-s_0)^{-(\frac{n-1}{2}+\alpha)}$$
$$\times \int_{\rho e^{\frac{s_0}{2}}-1}^{+\infty}  z^{n-1+\alpha+2\lambda_l}(1+c\frac{yz}{s-s_0})^{-(\frac{n}{2}-1+\alpha)}\exp(-\frac{(e^{-\frac{s-s_0}{2}}y-z)^2}{4(s-s_0)})dz$$
$$\leq C(n,l)e^{-\lambda_l s_0}y^\alpha (s-s_0)^{-(\frac{n-1}{2}+\alpha)}\int_{\rho e^{\frac{s_0}{2}}-1}^{+\infty}  z^{n-1+\alpha+2\lambda_l}e^{-\frac{z^2}{16(s-s_0)}}dz$$
($z\geq \rho e^{\frac{s_0}{2}}-1\geq 4R\geq 2y$ provided $s_0\gg 1$ (depending on $\rho,R$))
$$=C(n,l)e^{-\lambda_l s_0}y^\alpha (s-s_0)^{\lambda_l+\frac{1-\alpha}{2}}\int_{(\rho e^{\frac{s_0}{2}}-1)(s-s_0)^{-\frac{1}{2}}}^{+\infty}  w^{n-1+\alpha+2\lambda_l}e^{-\frac{w^2}{16}}dw\ (w=\frac{z}{\sqrt{s-s_0}})$$
$$\leq C(n,l)e^{-\lambda_l s_0}y^\alpha \int_{\rho e^{\frac{s_0}{2}}-1}^{+\infty}  w^{n-1+\alpha+2\lambda_l}e^{-\frac{w^2}{16}}dw\leq C(n,l)e^{-\lambda_l s_0}y^\alpha e^{-\frac{1}{32}(\rho e^{\frac{s_0}{2}}-1)^2}\leq C(n,l)e^{-(1+\tilde{\kappa})\lambda_l s_0}y^\alpha,$$
if $s_0\gg 1$ (depending on $n,l,\rho,\vartheta$).\\
The proof of this lemma is finished by combining the estimates of $T_1,T_2,T_3$.

\QED

As a consequence, we are able to prove the short-time estimate of Proposition \ref{propes1}:
\begin{cor}
    \label{cores1}
    Under the assumptions of Proposition \ref{propes1}, the estimate \eqref{es23} holds for all $s_0\leq s\leq \min\{\mathring{s},s_0+1\}$, $\frac{1}{2} e^{-\vartheta \sigma_l s}\leq y\leq 2R$.
\end{cor}
\proof
If $s_0\gg 1$ (depending on $n,l,\rho,\beta,R,\vartheta$), then for $\frac{1}{2} e^{-\vartheta \sigma_l s}\leq y\leq 2R$, $v(y,s)=\tilde{v}(y,s)$, and according to Lemmas \ref{lemes2} and \ref{lemes3}, 
$$v(y,s)-\frac{k}{c_l}e^{-\lambda_l s}\varphi_l(y)=\tilde{v}(y,s)-\frac{k}{c_l}e^{-\lambda_l s}\varphi_l(y)$$
$$=\int_0^{+\infty} K(y,z,s-s_0)\tilde{v}(z,s_0)dz+\int_{s_0}^s \int_0^{+\infty} K(y,z,s-\tau)f(z,\tau)dzd\tau-\frac{k}{c_l}e^{-\lambda_l s}\varphi_l(y)$$
$$=S(y,s)+T(y,s)+\langle \tilde{v}(\cdot,s_0),\varphi_l \rangle\int_0^{+\infty} K(y,z,s-s_0)\varphi_l(z)dz-\frac{k}{c_l}e^{-\lambda_l s}\varphi_l(y)$$
$$=S(y,s)+T(y,s)+(\langle \tilde{v}(\cdot,s_0),\varphi_l \rangle e^{\lambda_l s_0}-\frac{1}{c_l})e^{-\lambda_l s}\varphi_l(y)-\frac{k-1}{c_l}e^{-\lambda_l s}\varphi_l(y).$$
By Lemmas \ref{lemes4} and \ref{lemes5}, if $s_0\leq s\leq \min\{\mathring{s},s_0+1\}$,
$$|S(y,s)|+|T(y,s)|\leq C(p,q,l,\Lambda,\beta,R)e^{-(1+\tilde{\kappa})\lambda_ls_0} y^\alpha\leq C(p,q,l,\Lambda,\beta,R)e^{-(1+\tilde{\kappa})\lambda_ls} y^\alpha,$$
(if $s=s_0$, then the estimate holds obviously.) By Lemma \ref{lemco1} and \eqref{co8},
$$|\langle \tilde{v}(\cdot,s_0),\varphi_l \rangle e^{\lambda_l s_0}-\frac{1}{c_l}|\leq C(n,l,\beta)e^{-2\varsigma\lambda_l s_0}\leq C(n,l,\beta)e^{-2\varsigma\lambda_l s}.$$
By \eqref{es22}, 
$$|k-1|\leq C(p,q,l,\Lambda,\beta)e^{-\varsigma\lambda_l s_0}\leq C(p,q,l,\Lambda,\beta)e^{-\varsigma\lambda_l s}.$$
And we know from \eqref{ad11} that
$$|\varphi_l(y)|\leq C(n,l,R)y^\alpha,\ 0<y\leq 2R.$$
Combining all the estimates above, the Corollary is proved.

\QED

The long time estimate ($s>s_0+1$) of Proposition \ref{propes1} is obtained by writing down the Fourier expansion of $\tilde{v}$ under the basis $\{\varphi_j\}_{j\geq 0}$, and directly estimating the infinite series.
\begin{lemma}
    \label{lemes6}
    For all $s_0\leq s\leq \mathring{s}$,
    \begin{equation}
        \tilde{v}(\cdot,s)=\sum_{j=0}^{+\infty} e^{-\lambda_j(s-s_0)}\langle \tilde{v}(\cdot,s_0),\varphi_j\rangle \cdot \varphi_j+\sum_{j=0}^{+\infty} \int_{s_0}^s e^{-\lambda_j(s-\tau)}\langle f(\cdot,\tau),\varphi_j\rangle d\tau \cdot \varphi_j, \label{es28}
    \end{equation}
    where both of the two infinite sums converge in $\mathbf{H}$. Moreover, for $s_0<s\leq \mathring{s}$,
    \begin{equation}
        \int_0^{+\infty} K(y,z,s-s_0)\tilde{v}(z,s_0)dz=\sum_{j=0}^{+\infty} e^{-\lambda_j(s-s_0)}\langle \tilde{v}(\cdot,s_0),\varphi_j\rangle \cdot \varphi_j(y), \label{es31}
    \end{equation}
    \begin{equation}
        \int_{s_0}^s \int_0^{+\infty} K(y,z,s-\tau)f(z,\tau)dzd\tau=\sum_{j=0}^{+\infty} \int_{s_0}^s e^{-\lambda_j(s-\tau)}\langle f(\cdot,\tau),\varphi_j\rangle d\tau \cdot \varphi_j(y), \label{es32}
    \end{equation}
    as elements in $\mathbf{H}$.
\end{lemma}
\proof
$\{\langle \tilde{v}(\cdot,s_0),\varphi_j\rangle\}_{j\geq 0}$ is the Fourier coefficient of the function $\tilde{v}(\cdot,s_0)$, so
$$\sum_{j=2}^{+\infty} (e^{-\lambda_j(s-s_0)}\langle \tilde{v}(\cdot,s_0),\varphi_j\rangle)^2\leq \sum_{j=0}^{+\infty} \langle \tilde{v}(\cdot,s_0),\varphi_j\rangle^2<+\infty.$$
Similarly, $\{\langle f(\cdot,\tau),\varphi_j\rangle\}_{j\geq 0}$ is the Fourier coefficient of the function $f(\cdot,\tau)$, and
$$\sum_{j=2}^{+\infty} (\int_{s_0}^s e^{-\lambda_j(s-\tau)}\langle f(\cdot,\tau),\varphi_j\rangle d\tau)^2\leq \sum_{j=0}^{+\infty} (\int_{s_0}^s |\langle f(\cdot,\tau),\varphi_j\rangle| d\tau)^2\leq \{\int_{s_0}^s (\sum_{j=0}^{+\infty}\langle f(\cdot,\tau),\varphi_j\rangle^2)^{\frac{1}{2}}d\tau\}^2$$
\begin{equation}
    \leq ((s-s_0)\sup_{s_0\leq \tau\leq s}|f(\cdot,\tau)|_{\mathbf{H}})^2<+\infty, \label{es34}
\end{equation}
thus both series converge in $\mathbf{H}$.\\
For all $j\geq 0$, the function $\langle \tilde{v}(\cdot,s),\varphi_j \rangle$ satisfies the following ODE (see \eqref{es27}):
$$\partial_s \langle \tilde{v}(\cdot,s),\varphi_j \rangle +\lambda_j \langle \tilde{v}(\cdot,s),\varphi_j \rangle=\langle f(\cdot,s),\varphi_j \rangle,$$
thus we have
$$\langle\tilde{v}(\cdot,s),\varphi_j \rangle=e^{-\lambda_j(s-s_0)}\langle \tilde{v}(\cdot,s_0),\varphi_j\rangle+\int_{s_0}^s e^{-\lambda_j(s-\tau)}\langle f(\cdot,\tau),\varphi_j\rangle d\tau.$$
So \eqref{es28} is indeed the Fourier expansion of $\tilde{v}$ under the basis $\{\varphi_j\}_{j\geq 0}$.\\
To prove \eqref{es31}, we first observe from \eqref{es20} and \eqref{es17} that for $s>0$, $y>0$, $K(y,z,s)(z^{n-2}e^{-\frac{z^2}{4}})^{-1}\in \mathbf{H}$, where the integral variable is $z$. Furthermore, by Lemma \ref{lemes3}, for any $M\geq 0$,
\begin{equation}
    \int_0^{+\infty} K(y,z,s-s_0)(\sum_{j=0}^{M} \langle \tilde{v}(\cdot,s_0),\varphi_j\rangle \cdot \varphi_j(z))dz=\sum_{j=0}^{M} e^{-\lambda_j(s-s_0)}\langle \tilde{v}(\cdot,s_0),\varphi_j\rangle \cdot \varphi_j(y). \label{es33}
\end{equation}
As $M\rightarrow +\infty$,
$$\sum_{j=0}^{M} \langle \tilde{v}(\cdot,s_0),\varphi_j\rangle \cdot \varphi_j\rightarrow \tilde{v}(\cdot,s_0)$$
in $\mathbf{H}$, so the left hand side of \eqref{es33} converges to 
$$\int_0^{+\infty} K(y,z,s-s_0)\tilde{v}(z,s_0)dz$$
in a pointwise manner, for $y>0$. On the other hand, the right hand side of \eqref{es33} converges to 
$$\sum_{j=0}^{+\infty} e^{-\lambda_j(s-s_0)}\langle \tilde{v}(\cdot,s_0),\varphi_j\rangle \cdot \varphi_j$$
in $\mathbf{H}$, thus the two sides of \eqref{es31} must be the same element in $\mathbf{H}$.\\
\eqref{es32} is a direct consequence of \eqref{es28}, \eqref{es31}, and Lemma \ref{lemes2}.

\QED

To estimate the infinite sum directly, we first derive some uniform bounds on the eigenfunctions $\{\varphi_j\}_{j\geq 0}$, using the properties of Kummer's functions:
\begin{lemma}
    \label{lemes7}
    For any $R\geq 1$, $j\geq 1$,
    $$|\varphi_j(y)|\leq C(n,R)j^{\frac{n-3+2\alpha}{4}} y^\alpha,\ 0<y\leq 2R.$$
\end{lemma}
\proof
Recall \eqref{ad15}, the expression of $\varphi_j$. Set $b=\alpha+\frac{n-1}{2}$. According to \cite{tem}, (10.3.58), for $j\geq 1$ and $0<y\leq 2R$,
\begin{equation}
    |M(-j,b,\frac{y^2}{4})|\leq C(n,R). \label{es29}
\end{equation}
On the other hand, the normalizing constant $c_j$ satisfies
$$c_j^{-2}=\int_0^{+\infty}y^{n-2+2\alpha}e^{-\frac{y^2}{4}} M(-j,b,\frac{y^2}{4})^2dy=2^{n-2+2\alpha}\int_0^{+\infty}x^{b-1}e^{-x} M(-j,b,x)^2dx.$$
By \cite{gp}, (A.150),
$$\int_0^{+\infty}x^{b-1}e^{-x} M(-j,b,x)^2dx=\frac{\Gamma(b)j!}{b(b+1)...(b+j-1)}=\frac{\Gamma(b)^2 \Gamma(j+1)}{\Gamma(b+j)},$$
i.e.
$$c_j=\frac{2^{-\frac{1}{2}(n-2+2\alpha)}}{\Gamma(b)}\sqrt{\frac{\Gamma(b+j)}{\Gamma(j+1)}}\sim \frac{2^{-\frac{1}{2}(n-2+2\alpha)}}{\Gamma(b)} j^{\frac{b-1}{2}}\ (j\rightarrow +\infty),$$
using Stirling's formula $\Gamma(x+1)\sim \sqrt{2\pi x}(\frac{x}{e})^x$ ($x\rightarrow +\infty$). Therefore,
\begin{equation}
    c_j\leq C(n)j^{\frac{b-1}{2}}, \label{es30}
\end{equation}
and the lemma follows from \eqref{ad15}, \eqref{es29}, and \eqref{es30}.

\QED

The estimate of ``higher frequency" terms in a long time period is also divided into two parts, the part caused by evolution of the initial value $\tilde{v}(\cdot,s_0)$ and the ``non-homogeneous" term $f$, respectively. We proceed the initial value part first:
\begin{lemma}
    \label{lemes8}
    Let 
    $$P(y,s)=\sum_{j=l+1}^{+\infty} e^{-\lambda_j(s-s_0)}\langle \tilde{v}(\cdot,s_0),\varphi_j\rangle \cdot \varphi_j(y).$$
    If $\mathring{s}>s_0+1$, then the series converges for all $y>0$, $s_0+1\leq s\leq \mathring{s}$. Moreover, if $s_0\gg 1$ (depending on $n,l,\rho,\beta$), then for any $R\geq 1$, there holds
    $$|P(y,s)|\leq C(n,l,\beta,R)e^{-(1+\tilde{\kappa})\lambda_l s}y^\alpha,\ 0<y\leq 2R,\ s_0+1\leq s\leq \mathring{s}.$$
\end{lemma}
\proof
Since $s\geq s_0+1$ and $\lambda_j=\lambda_l+j-l$ for all $j$, 
$$|P(y,s)|\leq \sum_{j=l+1}^{+\infty} e^{-(\lambda_l+1)(s-s_0)}e^{-(j-l-1)(s-s_0)}|\tilde{v}(\cdot,s_0)-e^{-\lambda_l s_0}(\frac{1}{c_l}\varphi_l+\sum_{m=0}^{l-1} \frac{a_m}{c_m}\varphi_m)|_{\mathbf{H}}\ |\varphi_j|_{\mathbf{H}}\ |\varphi_j(y)|$$
$$\leq |\tilde{v}(\cdot,s_0)-e^{-\lambda_l s_0}(\frac{1}{c_l}\varphi_l+\sum_{m=0}^{l-1} \frac{a_m}{c_m}\varphi_m)|_{\mathbf{H}}\sum_{j=l+1}^{+\infty} e^{-(1+\tilde{\kappa})\lambda_l(s-s_0)}e^{-(j-l-1)}|\varphi_j(y)|$$
$$\leq C(n,l,\beta,R)e^{-(1+\tilde{\kappa})\lambda_l(s-s_0)}e^{-(1+\tilde{\kappa})\lambda_l s_0}y^\alpha\sum_{j=l+1}^{+\infty}e^{-(j-l-1)}j^{\frac{n-3+2\alpha}{4}}\ (\text{by Lemmas \ref{lemco1} and \ref{lemes7}})$$
$$=C(n,l,\beta,R)e^{-(1+\tilde{\kappa})\lambda_l s}y^\alpha.$$
Since $R\geq 1$ is arbitrary, the series of course converges for all $y>0$.

\QED

The following lemma deals with the ``non-homogeneous" part:
\begin{lemma}
    \label{lemes9}
    Let 
    $$Q(y,s)=\sum_{j=l+1}^{+\infty} \int_{s_0}^s e^{-\lambda_j(s-\tau)}\langle f(\cdot,\tau),\varphi_j\rangle d\tau \cdot \varphi_j(y).$$
    If $\mathring{s}>s_0+1$, \eqref{es1}, \eqref{es2}, \eqref{es3} hold, and $s_0\gg 1$ (depending on $n,l,\rho,\beta,R,\vartheta$), then for any $s_0+1\leq s\leq \mathring{s}$, any $R\geq 1$, and almost all $\frac{1}{2}e^{-\vartheta \sigma_l s}\leq y\leq 2R$, there holds
    $$|Q(y,s)|\leq C(p,q,l,\Lambda,\beta,R)e^{-(1+\tilde{\kappa})\lambda_l s}y^\alpha.$$
\end{lemma}
\proof
Formally we can write (noting $s\geq s_0+1$)
$$Q(y,s)=\sum_{j=l+1}^{+\infty} \int_{s_0}^{s-1} e^{-\lambda_j(s-\tau)}\langle f(\cdot,\tau),\varphi_j\rangle d\tau \cdot \varphi_j(y)+\sum_{j=l+1}^{+\infty} \int_{s-1}^s e^{-\lambda_j(s-\tau)}\langle f(\cdot,\tau),\varphi_j\rangle d\tau \cdot \varphi_j(y)$$
$$:=Q_1(y,s)+Q_2(y,s).$$
From \eqref{es34} we know both two series above converge in $\mathbf{H}$. We estimate $Q_1$ using Cauchy-Schwartz inequality:
$$|Q_1(y,s)|\leq \int_{s_0}^{s-1} (\sum_{j=l+1}^{+\infty} \lambda_j^3 e^{-2\lambda_j(s-\tau)} \varphi_j(y)^2)^{\frac{1}{2}}(\sum_{j=l+1}^{+\infty} \frac{\langle f(\cdot,\tau),\varphi_j\rangle^2}{\lambda_j^3})^{\frac{1}{2}}d\tau$$
$$=\int_{s_0}^{s-1} e^{-\lambda_{l+1}(s-\tau)}(\sum_{j=l+1}^{+\infty} \lambda_j^3 e^{-2(\lambda_j-\lambda_{l+1})(s-\tau)} \varphi_j(y)^2)^{\frac{1}{2}}(\sum_{j=l+1}^{+\infty} \frac{\langle f(\cdot,\tau),\varphi_j\rangle^2}{\lambda_j^3})^{\frac{1}{2}}d\tau$$
$$\leq(\sum_{j=l+1}^{+\infty} \lambda_j^3 e^{-2(j-l-1)} \varphi_j(y)^2)^{\frac{1}{2}} \int_{s_0}^{s-1} e^{-(\lambda_l+1)(s-\tau)} |f(\cdot,\tau)|_{\mathbf{X}^*}(\sum_{j=l+1}^{+\infty} \frac{|\varphi_j|_{\mathbf{X}}^2}{\lambda_j^3})^{\frac{1}{2}}d\tau.$$
Using Lemma \ref{lemes7}, the uniform estimate of $\{\varphi_j\}_{j\geq 0}$, and noting $\{\lambda_j\}_{j\geq 0}$ is of linear growth, we get (for $0<y\leq 2R$)
$$(\sum_{j=l+1}^{+\infty} \lambda_j^3 e^{-2(j-l-1)} \varphi_j(y)^2)^{\frac{1}{2}}\leq C(n,l,R)y^\alpha<+\infty,$$
and by Lemma \ref{lemes1},
$$(\sum_{j=l+1}^{+\infty} \frac{|\varphi_j|_{\mathbf{X}}^2}{\lambda_j^3})^{\frac{1}{2}}\leq C(n) (\sum_{j=0}^{+\infty} \frac{1+|\lambda_j|}{|\lambda_j|^3})^{\frac{1}{2}}=C(n)<+\infty$$
(again, $\{\lambda_j\}_{j\geq 0}$ is of linear growth). Thus
$$|Q_1(y,s)|\leq C(n,l,R)y^\alpha \int_{s_0}^{s-1} e^{-(\lambda_l+1)(s-\tau)} |f(\cdot,\tau)|_{\mathbf{X}^*}d\tau$$
$$\leq C(p,q,l,\Lambda,\beta,R)y^\alpha e^{-(\lambda_l+1)s}\int_{s_0}^{s-1}e^{(1-\tilde{\kappa}\lambda_l)\tau}d\tau\ (\text{by \eqref{es5}})$$
$$\leq C(p,q,l,\Lambda,\beta,R)y^\alpha e^{-(1+\tilde{\kappa})\lambda_l s}$$
for $0<y\leq 2R$, since $1-\tilde{\kappa}\lambda_l\geq \frac{1}{\lambda_l+1}>0$.\\
As for the estimate of $Q_2$, we write
$$Q_2(y,s)=\sum_{j=0}^{+\infty} \int_{s-1}^s e^{-\lambda_j(s-\tau)}\langle f(\cdot,\tau),\varphi_j\rangle d\tau \cdot \varphi_j(y)-\sum_{j=0}^l \int_{s-1}^s e^{-\lambda_j(s-\tau)}\langle f(\cdot,\tau),\varphi_j\rangle d\tau \cdot \varphi_j(y)$$
$$=\int_{s-1}^s \int_0^{+\infty} K(y,z,s-\tau)f(z,\tau)dzd\tau-\sum_{j=0}^l \int_{s-1}^s e^{-\lambda_j(s-\tau)}\langle f(\cdot,\tau),\varphi_j\rangle d\tau \cdot \varphi_j(y)\ (\text{by \eqref{es32}})$$
$$:=Q_{2,1}+Q_{2,2}.$$
According to Lemma \ref{lemes4} (with $s-1$ in the place of $s_0$), 
$$|Q_{2,1}|\leq C(p,q,l,\Lambda,\beta,R)e^{-(1+\tilde{\kappa})\lambda_l (s-1)}y^\alpha\leq  C(p,q,l,\Lambda,\beta,R)e^{-(1+\tilde{\kappa})\lambda_l s}y^\alpha,\ \frac{1}{2}e^{-\vartheta \sigma_l s}\leq y\leq 2R,$$
and using \eqref{es5} again,
$$|Q_{2,2}|\leq C(n) \sup_{s-1\leq \tau\leq s}|f(\cdot,\tau)|_{\mathbf{X}^*} \sum_{j=0}^l |\varphi_j|_{\mathbf{X}}|\varphi_j(y)|\leq  C(p,q,l,\Lambda,\beta,R)e^{-(1+\tilde{\kappa})\lambda_l s}y^\alpha,\ 0<y\leq 2R.$$
The lemma is proved by adding up the estimates of $Q_1,Q_{2,1},Q_{2,2}$.

\QED

Now it's time to prove the long-time estimate of Proposition \ref{propes1}:
\begin{cor}
    \label{cores2}
    Under the assumptions of Proposition \ref{propes1}, if $\mathring{s}>s_0+1$, then the estimate \eqref{es23} holds for all $s_0+1\leq s\leq \mathring{s}$, $\frac{1}{2} e^{-\vartheta \sigma_l s}\leq y\leq 2R$.
\end{cor}
\proof
If $s_0\gg 1$ (depending on $n,l,\rho,\beta,R,\vartheta$), then for $\frac{1}{2} e^{-\vartheta \sigma_l s}\leq y\leq 2R$, $v(y,s)=\tilde{v}(y,s)$. According to Lemma \ref{lemes6},
$$v(y,s)-\frac{k}{c_l}e^{-\lambda_l s}\varphi_l(y)=\tilde{v}(y,s)-\frac{k}{c_l}e^{-\lambda_l s}\varphi_l(y)$$
$$=\sum_{j=0}^{l-1} \langle \tilde{v}(\cdot,s),\varphi_j\rangle \cdot \varphi_j(y)+(\langle \tilde{v}(\cdot,s),\varphi_l\rangle-\frac{k}{c_l}e^{-\lambda_l s})\varphi_l(y)$$
$$+\sum_{j=l+1}^{+\infty} e^{-\lambda_j(s-s_0)}\langle \tilde{v}(\cdot,s_0),\varphi_j\rangle \cdot \varphi_j(y)+\sum_{j=l+1}^{+\infty} \int_{s_0}^s e^{-\lambda_j(s-\tau)}\langle f(\cdot,\tau),\varphi_j\rangle d\tau \cdot \varphi_j(y)$$
$$=\sum_{j=0}^{l-1} \langle \tilde{v}(\cdot,s),\varphi_j\rangle \cdot \varphi_j(y)+(\langle \tilde{v}(\cdot,s),\varphi_l\rangle-\frac{k}{c_l}e^{-\lambda_l s})\varphi_l(y)+P(y,s)+Q(y,s).$$
All the equalities above, except the first one, mean the two sides represent the same element in $\mathbf{H}$. By Lemmas \ref{lemes8} and \ref{lemes9}, if $s_0+1\leq s\leq \mathring{s}$, then for almost all $\frac{1}{2} e^{-\vartheta \sigma_l s}\leq y\leq 2R$,
$$|P(y,s)|+|Q(y,s)|\leq C(p,q,l,\Lambda,\beta,R)e^{-(1+\tilde{\kappa})\lambda_l s}y^\alpha.$$
According to \eqref{es11} and \eqref{es35}, for all $0<y\leq 2R$,
$$\sum_{j=0}^{l-1} |\langle \tilde{v}(\cdot,s),\varphi_j\rangle| \cdot |\varphi_j(y)|+|\langle \tilde{v}(\cdot,s),\varphi_l\rangle-\frac{k}{c_l}e^{-\lambda_l s}|\cdot |\varphi_l(y)|\leq C(p,q,l,\Lambda,\beta,R)e^{-(1+\tilde{\kappa})\lambda_l s}y^\alpha.$$
Thus \eqref{es23} holds for all $\frac{1}{2} e^{-\vartheta \sigma_l s}\leq y\leq 2R$, and actually all these $y$, because $v(\cdot,s)$ is continuous.

\QED

The proof of Proposition \ref{propes1} is now finished, by putting Corollaries \ref{cores1} and \ref{cores2} together.

\QED

For brevity, in the rest part of this and the next sections, we always fix $\vartheta=\frac{1}{2}$, and therefore $\tilde{\kappa}=\kappa$ (see \eqref{co9}). It's noteworthy that the choice of $\vartheta\in (\frac{-1-\alpha}{1-\alpha},1)$ makes no essential difference in the following discussion. The lower bound $\vartheta>\frac{-1-\alpha}{1-\alpha}$ is necessary in the proof of Proposition \ref{propes3}; see \cite{gs}, Proposition 6.6.\\
Next, we provide an estimate in the outer region, for the function $u(x,t)$ with $R\sqrt{-t}\lesssim x\lesssim \rho$. The original proof appears in \cite{vel}, Lemma 4.3. By a simple observation, it's not hard to see the argument applies to both even and odd $l$.
\begin{prop}
    \label{propes2}
    If $0<\rho\ll 1$ (depending on $p,q,l,\Lambda$), $R\gg 1$ (depending on $n,l$), and $|t_0|\ll 1$ (depending on $p,q,l,\Lambda,\rho,\beta,R$), then
    $$|u(x,t)-\frac{k}{c_l}(-t)^{\lambda_l+\frac{1}{2}}\varphi_l(\frac{x}{\sqrt{-t}})|\leq C(n,l)R^{-2}x^{2\lambda_l+1}$$
        for $2R\sqrt{-t}\leq x\leq \rho$, $t_0\leq t\leq \mathring{t}$.
\end{prop}
\proof
We prove it by constructing sub- and supersolutions.\\
To investigate the equation \eqref{ad2} further, we rewrite it as
\begin{equation}
    \partial_t u=\tilde{\mathcal{L}}u+\mathcal{Q}u, \label{es36}
\end{equation}
where (see \eqref{ad13})
$$\tilde{\mathcal{L}}u=u''+\frac{n-2}{x}u'+\frac{n-2}{x^2}u,$$
$$\mathcal{Q}u=-\frac{u'^2}{1+u'^2}u''+(n-2)\frac{(\frac{u}{x})^2(\frac{u}{x^2}+\frac{u'}{x})+(\mu-\mu^{-1})\frac{u}{x}\frac{u}{x^2}}{(1-\mu \frac{u}{x})(1+\mu^{-1}\frac{u}{x})}.$$
Assume $l$ is an even number first. Define
$$u^\pm(x,t)=C_0^\pm(x^{2\lambda_l+1}-C^\pm(-t)x^{2\lambda_l-1}),$$
where $C_0^\pm>0$, $C^\pm\geq 0$ are constants to be determined later, with $C_0^\pm$ bounded and away from 0 (depending on $n,l$). Direct computation yields
$$(\partial_t-\tilde{\mathcal{L}})u^\pm=C_0^\pm x^{2\lambda_l-1}(C^\pm-M_1+C^\pm M_2\frac{-t}{x^2}),$$
$$M_1=(2\lambda_l+1)(2\lambda_l)+(n-2)(2\lambda_l+2)>0,$$
$$M_2=(2\lambda_l-1)(2\lambda_l-2)+(n-2)(2\lambda_l)>0.$$
Set $C^+=2M_1$, $C^-=0$, then
$$(\partial_t-\tilde{\mathcal{L}})u^+\geq C_0^\pm M_1 x^{2\lambda_l-1}>0,\ (\partial_t-\tilde{\mathcal{L}})u^-=-C_0^\pm M_1 x^{2\lambda_l-1}<0.$$
Since $C_0^\pm$ is bounded and away from 0 (depending on $n,l$), if $t<0$, $\sqrt{-t}\leq x\leq \rho$ and $\rho\ll 1$ (depending on $p,q,l$), then
$$|\mathcal{Q}u^\pm|\leq C(p,q,l)x^{4\lambda_l-1}\leq C_0^\pm M_1 x^{2\lambda_l-1},$$
i.e. $u^+$ (resp. $u^-$) is a supersolution (resp. subsolution) of \eqref{ad2}.\\
Next we verify the initial and boundary values of $u$. Under the transform \eqref{ad6}, Proposition \ref{propes1} implies
\begin{equation}
    |u(x,t)-\frac{k}{c_l}(-t)^{\lambda_l+\frac{1}{2}}\varphi_l(\frac{x}{\sqrt{-t}})|\leq C(p,q,l,\Lambda,\beta,R)(-t)^{\kappa \lambda_l}x^{2\lambda_l+1},\ x=2R\sqrt{-t},\ t_0\leq t\leq \mathring{t}. \label{es38}
\end{equation}
The formula of $v(\cdot,s_0)$, \eqref{co8}, together with the estimates \eqref{es21}, \eqref{es22}, indicates
$$|v(y,s_0)-\frac{k}{c_l}e^{-\lambda_l s_0}\varphi_l(y)|\leq C(p,q,l,\Lambda,\beta)e^{-(1+\kappa)\lambda_l s_0}y^{2\lambda_l+1},\ 2R\leq y\leq \rho e^{\frac{s_0}{2}},$$
i.e.
\begin{equation}
    |u(x,t_0)-\frac{k}{c_l}(-t_0)^{\lambda_l+\frac{1}{2}}\varphi_l(\frac{x}{\sqrt{-t_0}})|\leq C(p,q,l,\Lambda,\beta)(-t_0)^{\kappa \lambda_l}x^{2\lambda_l+1},\ 2R\sqrt{-t_0}\leq x\leq \rho. \label{es37}
\end{equation}
By \eqref{es36} and the admissible condition, \eqref{ad5}, if $t_0\leq t\leq \mathring{t}$, $\sqrt{-t}\leq x\leq \rho$ and $\rho\ll 1$ (depending on $p,q,l,\Lambda$), then
$$|\partial_t u|\leq |u''|+(n-2)|\frac{u'}{x}|+(n-2)|\frac{u}{x^2}|+|\mathcal{Q}u|\leq C(n,\Lambda)x^{2\lambda_l-1}.$$
Moreover, by \eqref{ad11},
$$\partial_t (\frac{k}{c_l}(-t)^{\lambda_l+\frac{1}{2}}\varphi_l(\frac{x}{\sqrt{-t}}))=k\partial_t ((-t)^l x^\alpha-K_{l,1}(-t)^{l-1}x^{\alpha+2}+K_{l,2}(-t)^{l-2}x^{\alpha+4}+...+(-1)^l K_{l,l}x^{\alpha+2l})$$
$$=-k(l(-t)^{l-1} x^\alpha-(l-1)K_{l,1}(-t)^{l-2}x^{\alpha+2}+(l-2)K_{l,2}(-t)^{l-3}x^{\alpha+4}+...+(-1)^{l-1} K_{l,l-1}x^{\alpha+2(l-1)}).$$
If $t<0$ and $x\geq \sqrt{-t}$, then
$$|\partial_t (\frac{k}{c_l}(-t)^{\lambda_l+\frac{1}{2}}\varphi_l(\frac{x}{\sqrt{-t}}))|\leq C(n,l)x^{2\lambda_l-1},$$
and in particular, for $t_0\leq t\leq \mathring{t}$,
$$|\partial_t (u(\rho,t)-\frac{k}{c_l}(-t)^{\lambda_l+\frac{1}{2}}\varphi_l(\frac{\rho}{\sqrt{-t}}))|\leq C(n,l,\Lambda)x^{2\lambda_l-1}=C(n,l,\Lambda)\rho^{-2}x^{2\lambda_l+1}.$$
Taking \eqref{es37} into account, we get
\begin{equation}
    |u(\rho,t)-\frac{k}{c_l}(-t)^{\lambda_l+\frac{1}{2}}\varphi_l(\frac{\rho}{\sqrt{-t}})|\leq C(p,q,l,\Lambda,\rho,\beta)(-t_0)^{\kappa \lambda_l}x^{2\lambda_l+1},\ t_0\leq t\leq \mathring{t}. \label{es39}
\end{equation}
Set
$$\Omega=\{(x,t)\subset \R\times \R|\ t_0<t<\mathring{t},\ 2R\sqrt{-t}<x<\rho\},$$
then \eqref{es38}, \eqref{es37}, \eqref{es39} together imply that on $\mathscr{P}\Omega$, the parabolic boundary of $\Omega$ (see \ref{pesub1}),
$$|u(x,t)-\frac{k}{c_l}(-t)^{\lambda_l+\frac{1}{2}}\varphi_l(\frac{x}{\sqrt{-t}})|\leq C(p,q,l,\Lambda,\rho,\beta,R)(-t_0)^{\kappa \lambda_l}x^{2\lambda_l+1}.$$
On the other hand, using \eqref{ad11} again, 
\begin{equation}
    |\frac{k}{c_l}(-t)^{\lambda_l+\frac{1}{2}}\varphi_l(\frac{x}{\sqrt{-t}})-kK_{l,l}x^{2\lambda_l+1}|\leq C(n,l)R^{-2}x^{2\lambda_l+1},\ t<0,\ x\geq 2R\sqrt{-t}, \label{es41}
\end{equation}
namely, on $\mathscr{P}\Omega$,
$$|u(x,t)-kK_{l,l}x^{2\lambda_l+1}|\leq (C(p,q,l,\Lambda,\rho,\beta,R)(-t_0)^{\kappa \lambda_l}+C(n,l)R^{-2})x^{2\lambda_l+1}$$
\begin{equation}
   \leq C(n,l)R^{-2}x^{2\lambda_l+1}\ \text{if}\ |t_0|\ll 1\ (\text{depending on}\ p,q,l,\Lambda,\rho,\beta,R). \label{es40}
\end{equation}
Now we can determine the constants $C_0^\pm$ in $u^\pm$. Let
$$C_0^+=\{kK_{l,l}+C_{\eqref{es40}}(n,l)R^{-2}\}(1-\frac{M_1}{2R^2})^{-1},$$
$$C_0^-=kK_{l,l}-C_{\eqref{es40}}(n,l)R^{-2}.$$
Here the subscript indicates we choose exactly that constant appeared in \eqref{es40}. Obviously, if $R\gg 1$ (depending on $n,l$), and $|t_0|\ll 1$ (depending on $p,q,l,\Lambda,\beta$), then (recall \eqref{es22})
$$\frac{1}{2}K_{l,l}\leq C_0^-<C_0^+\leq 2K_{l,l},$$
i.e. they are indeed bounded and away from 0 (depending on $n,l$). And on $\mathscr{P}\Omega$,
$$u^+=C_0^+ x^{2\lambda_l+1}(1-2M_1\frac{-t}{x^2})\geq C_0^+ x^{2\lambda_l+1}(1-\frac{M_1}{2R^2})$$
$$=\{kK_{l,l}+C_{\eqref{es40}}(n,l)R^{-2}\}x^{2\lambda_l+1}\geq u,$$
$$u^-=\{kK_{l,l}-C_{\eqref{es40}}(n,l)R^{-2}\}x^{2\lambda_l+1}\leq u.$$
Making use of Theorem \ref{thmpe2} (the comparison principle), we deduce for all $(x,t)\in \bar{\Omega}$,
$$C_0^- x^{2\lambda_l+1}=u^-\leq u\leq u^+\leq C_0^+ x^{2\lambda_l+1},$$
in other words, \eqref{es40} holds on the whole $\bar{\Omega}$ with a probably larger constant, by noting an elementary fact that $(1-\frac{M_1}{2R^2})^{-1}\leq 1+\frac{M_1}{R^2}$ when $R\gg 1$ (depending on $n,l$). Combining this with \eqref{es41}, we finally get
$$|u(x,t)-\frac{k}{c_l}(-t)^{\lambda_l+\frac{1}{2}}\varphi_l(\frac{x}{\sqrt{-t}})|\leq C(n,l)R^{-2}x^{2\lambda_l+1},\ (x,t)\in \bar{\Omega}.$$
The situation is very similar if $l$ is odd, except that the subsolution and the supersolution are reversed. More precisely, we shall let $C_0^\pm<0$ now, still bounded and away from 0 (depending on $n,l$). If we set $C^+=0$, $C^-=2M_1$ in this case, then we can get a subsolution $u^+$ and a supersolution $u^-$ in the same way. \eqref{es41} and \eqref{es40} remain valid, if we replace $K_{l,l}$ with $-K_{l,l}$. The final choice of $C_0^\pm$ is
$$C_0^+=-kK_{l,l}+C_{\eqref{es40}}(n,l)R^{-2},$$
$$C_0^-=\{-kK_{l,l}-C_{\eqref{es40}}(n,l)R^{-2}\}(1-\frac{M_1}{2R^2})^{-1}.$$

\QED

Below is an estimate in the tip region, for the function $\hat{w}(z,\tau)$ with $0\leq z\lesssim (2\sigma_l \tau)^{\frac{1}{4}}$. Again, this is proved by the method of sub- and supersolutions, which is exactly the same as \cite{gs}, Proposition 6.6, and thus we omit the proof here.
\begin{prop}
    \label{propes3}
    If $\beta\gg 1$ (depending on $p,q,l$), and $\tau_0\gg 1$ (depending on $p,q,l,\Lambda,\rho,\beta,R$), then
    $$|\hat{w}(z,\tau)-\hat{\psi}_k(z)|\leq C(p,q)\beta^{\tilde{\alpha}-\alpha}(\frac{\tau}{\tau_0})^{-\varrho}(1+z)^\alpha,\ 0\leq z\leq \frac{2(2\sigma_l \tau)^{\frac{1}{4}}}{\sqrt{1+\mu^2}},\ \tau_0\leq \tau\leq \mathring{\tau},$$
    where $\tilde{\alpha}$ is defined in \eqref{mi20}, and $\varrho$ is defined in \eqref{co10}.
\end{prop}

Finally, we provide an estimate of the remaining part, by virtue of a well-known interior estimate for hypersurfaces moving by mean curvature in $\R^n$, due to Ecker and Huisken (\cite{eh}, Theorems 2.1 and 3.1). The original estimate is ``interior" in both space and time. We first state a slightly variant version of the theorem, which can extend the estimate to the initial time, and simplify subsequent arguments.
\begin{theorem}
    \label{thmes1}
    Let $\{M^{n-1}_t\}_{0\leq t\leq T}\subset \R^n$  be a smooth family of embedded hypersurfaces moving by mean curvature. Assume $E=\{(\mathbf{x},t)|\ \mathbf{x}\in M_t,0\leq t\leq T,r(\mathbf{x},t)\leq L^2\}$ is compact, and on $E$, $\langle \mathbf{v},\omega \rangle>0$, where $L>0$ is a constant, $r(\mathbf{x},t)=|\mathbf{x}-\mathbf{x}_0|^2+2(n-1)t$, $\mathbf{x}_0\in \R^n$ is a fixed point, $\omega\in \R^n$ is a fixed vector, and $\mathbf{v}=\mathbf{v}(\mathbf{x},t)$ is a unit normal vector of $M_t$ at $\mathbf{x}$. (In other words, $E$ can be regarded as graphs over the hyperplane perpendicular to $\omega$.) Then for any $t_0\in [0,T]$ and $\theta\in (0,1)$,
    \begin{equation}
        \sup_{\mathbf{x}\in M_{t_0},r(\mathbf{x},t_0)\leq \theta L^2} v^2\leq (1-\theta)^{-2}\sup_{\mathbf{x}\in M_0,r(\mathbf{x},0)\leq L^2} v^2, \label{es52}
    \end{equation}
    \begin{equation}
        \sup_{\mathbf{x}\in M_{t_0},r(\mathbf{x},t_0)\leq \theta L^2} |A|^2\leq 8(1-\theta)^{-2}\sup_{\mathbf{x}\in M_0,r(\mathbf{x},0)\leq L^2} |A|^2 v^2+C(n)L^{-2}(1-\theta)^{-6}\sup_{\mathbf{x}\in M_0,r(\mathbf{x},0)\leq L^2} v^4, \label{es53}
    \end{equation}
    where $v=\langle \mathbf{v},\omega \rangle^{-1}$ is the ``gradient function", and $A$ is the second fundamental form of $M_t$.
\end{theorem}

Actually, Theorem \ref{thmes1} remains true if the mean curvature flow is merely immersed.

\proof
\eqref{es52} follows immediately from Theorem 2.1 of \cite{eh}.\\
The proof of \eqref{es53} is essentially the same as that of Theorem 3.1 of \cite{eh}, except that we do not cut off along the $t$ direction. More precisely, in the original proof, the author derived the following inequality:
$$(\frac{d}{dt}-\Delta)(g\eta)\leq -2kg^2\eta -2(\varphi v^{-3} \nabla v+\eta^{-1}\nabla \eta)\cdot \nabla (g\eta) +C(n)((1+\frac{1}{kv^2})r+L^2)g,$$
where $\Delta$ and $\nabla$ are the Laplacian and gradient on the hypersurface $M_t$, $\varphi=\varphi(v^2)=\frac{v^2}{1-kv^2}$, $g=|A|^2 \varphi(v^2)$, $\eta=(L^2-r)^2$, $\frac{d}{dt}$ is the derivative along a curve $\gamma(t)\in M_t$ whose velocity equals the mean curvature of $M_t$, and 
$$k=\frac{1}{2}\inf_{\begin{subarray}{c}
    \mathbf{x}\in M_t,r(\mathbf{x},t)\leq L^2,\\
    0\leq t\leq t_0
\end{subarray}} v^{-2}>0.$$
At a point $(\mathbf{x}_1,t_1)$ where $\sup_{\begin{subarray}{c}
    \mathbf{x}\in M_t,r(\mathbf{x},t)\leq L^2,\\
    0\leq t\leq t_0
\end{subarray}} g\eta$ is attained, if $t_1=0$, then
$$\sup_{\mathbf{x}\in M_{t_0},r(\mathbf{x},t_0)\leq L^2} g\eta\leq \sup_{\mathbf{x}\in M_0,r(\mathbf{x},0)\leq L^2} g\eta.$$
If $t_1>0$, then at $(\mathbf{x}_1,t_1)$ we have ($\frac{d}{dt}(g\eta)\geq 0$, $\Delta(g\eta)\leq 0$, $\nabla (g\eta)=0$)
$$2kg^2\eta\leq C(n)((1+\frac{1}{kv^2})r+L^2)g,$$
i.e. (note $r\leq L^2$, $kv^2\leq \frac{1}{2}$, $v\geq 1$)
$$g\eta\leq \frac{C(n)}{k}(1+\frac{1}{kv^2})L^2\leq \frac{C(n)}{k^2 v^2}L^2\leq C(n)L^2\sup_{\begin{subarray}{c}
    \mathbf{x}\in M_t,r(\mathbf{x},t)\leq L^2,\\
    0\leq t\leq t_0
\end{subarray}} v^4.$$
Thus,
$$\sup_{\mathbf{x}\in M_{t_0},r(\mathbf{x},t_0)\leq L^2} g\eta\leq \sup_{\mathbf{x}\in M_0,r(\mathbf{x},0)\leq L^2} g\eta+C(n)L^2\sup_{\begin{subarray}{c}
    \mathbf{x}\in M_t,r(\mathbf{x},t)\leq L^2,\\
    0\leq t\leq t_0
\end{subarray}} v^4.$$
Since
$$g\eta=|A|^2 \frac{v^2}{1-kv^2}(L^2-r)^2\leq 2|A|^2 v^2 L^4,$$
$$g\eta\geq |A|^2(1-\theta)^2 L^4\ (r\leq \theta L^2),$$
we have
\begin{equation}
    \sup_{\mathbf{x}\in M_{t_0},r(\mathbf{x},t_0)\leq \theta L^2} |A|^2\leq 2(1-\theta)^{-2}\sup_{\mathbf{x}\in M_0,r(\mathbf{x},0)\leq L^2} |A|^2 v^2+C(n)L^{-2}(1-\theta)^{-2}\sup_{\begin{subarray}{c}
    \mathbf{x}\in M_t,r(\mathbf{x},t)\leq L^2,\\
    0\leq t\leq t_0
\end{subarray}} v^4. \label{es54}
\end{equation}
Replace $L^2$ in \eqref{es54} by $\frac{1+\theta}{2}L^2$ and $\theta$ by $\frac{2\theta}{1+\theta}$:
\begin{equation}
    \sup_{\mathbf{x}\in M_{t_0},r(\mathbf{x},t_0)\leq \theta L^2} |A|^2\leq 8(1-\theta)^{-2}\sup_{\mathbf{x}\in M_0,r(\mathbf{x},0)\leq L^2} |A|^2 v^2+C(n)L^{-2}(1-\theta)^{-2}\sup_{\begin{subarray}{c}
    \mathbf{x}\in M_t,r(\mathbf{x},t)\leq \frac{1+\theta}{2}L^2,\\
    0\leq t\leq t_0
\end{subarray}} v^4. \label{es55}
\end{equation}
Putting \eqref{es52} (with $\frac{1+\theta}{2}$ in the place of $\theta$) into the last term of \eqref{es55}, we arrive at \eqref{es53}.

\QED

The following statement justifies a part of the estimates in Proposition \ref{propco1}:
\begin{prop}
    \label{propes4}
    If $\rho\ll 1$ (depending on $p,q,l$), and $|t_0|\ll 1$ (depending on $p,q,l,\rho,\beta$), then for all $t_0\leq t\leq \mathring{t}$,
    \begin{enumerate}
        \item \label{propes4_1} The statement (\ref{propco1_2}) of Proposition \ref{propco1} holds.
        \item \label{propes4_2} The profile curve of $\Sigma_t\cap (B(0,3\rho)-\bar{B}(0,\frac{1}{3}\rho))$ can be parameterized by a single function as \eqref{ad1}, and \eqref{co11} holds for all $x\geq \frac{1}{2}\rho$ in this region.
    \end{enumerate}
\end{prop}
\proof
Let's first parameterize (a part of) the initial hypersurface $\Sigma_{t_0}^\mathbf{a}=\Sigma_{t_0}$ as \eqref{ad1}, i.e. 
\begin{equation}
    \mathbf{x}(u(x,t_0),\omega,\phi)=(\frac{x-\mu u(x,t_0)}{\sqrt{1+\mu^2}}\omega,\frac{\mu x+u(x,t_0)}{\sqrt{1+\mu^2}}\phi),\ \omega\in \sph^{p-1},\ \phi\in \sph^{q-1}, \label{es56}
\end{equation}
where $u(x,t_0)$ is constructed in \eqref{co2}, with $\frac{1}{5}\rho\leq x\leq 1$. The (outward) unit normal vector of $\Sigma_{t_0}$ at $\mathbf{x}=\mathbf{x}(u(x,t_0),\omega,\phi)$ is
$$N_{t_0}(\mathbf{x})=\frac{((-\mu-u'(x,t_0))\omega,(1-\mu u'(x,t_0))\phi)}{\sqrt{1+\mu^2}\sqrt{1+u'(x,t_0)^2}}.$$
According to \eqref{co1}, if $\rho\ll 1$ (depending on $p,q,l$) and $|t_0|\ll 1$ (depending on $n,l,\rho,\beta$), then for all $\frac{1}{6}\rho\leq x\leq 1$,
\begin{equation}
    |\frac{u(x,t_0)}{x}|,|u'(x,t_0)|,|xu''(x,t_0)|\leq \frac{1}{4}\min\{\mu,\mu^{-1}\}, \label{es61}
\end{equation}
and every point in $\Sigma_{t_0}\cap (\bar{B}(0,1)-B(0,\frac{1}{5}\rho))$ can be parameterized as \eqref{es56} with $\frac{1}{6}\rho\leq x\leq 1$. Now, for $\frac{1}{6}\rho\leq x_j\leq 1$, $\omega_j\in \sph^{p-1}$, $\phi_j\in \sph^{q-1}$, $\mathbf{x}_j=\mathbf{x}(u(x_j,t_0),\omega_j,\phi_j)$, $j=1,2$, 
$$|\mathbf{x}_1-\mathbf{x}_2|^2$$
$$=\frac{1}{1+\mu^2}(|(x_1-\mu u(x_1,t_0))\omega_1-(x_2-\mu u(x_2,t_0))\omega_2|^2+|(\mu x_1+u(x_1,t_0))\phi_1-(\mu x_2+u(x_2,t_0))\phi_2|^2)$$
$$=\frac{1}{1+\mu^2}\{((x_1-\mu u(x_1,t_0))-(x_2-\mu u(x_2,t_0)))^2+2(x_1-\mu u(x_1,t_0))(x_2-\mu u(x_2,t_0))(1-\langle \omega_1,\omega_2 \rangle)$$
$$+((\mu x_1+u(x_1,t_0))-(\mu x_2+u(x_2,t_0)))^2+2(\mu x_1+u(x_1,t_0))(\mu x_2+u(x_2,t_0))(1-\langle \phi_1,\phi_2 \rangle)\}$$
$$\geq \frac{2}{1+\mu^2}\{(x_1-\mu u(x_1,t_0))(x_2-\mu u(x_2,t_0))(1-\langle \omega_1,\omega_2 \rangle)+(\mu x_1+u(x_1,t_0))(\mu x_2+u(x_2,t_0))(1-\langle \phi_1,\phi_2 \rangle)\}$$
\begin{equation}
    \geq \frac{1}{32}\rho^2(1-\max\{\langle \omega_1,\omega_2 \rangle,\langle \phi_1,\phi_2 \rangle\}). \label{es57}
\end{equation}
For the fixed unit vector $\mathbf{e}=\frac{(-\mu\omega_1,\phi_1)}{\sqrt{1+\mu^2}}$, the ``gradient function" of $\Sigma_{t_0}$ (see Theorem \ref{thmes1}) at $\mathbf{x}_2$ is
\begin{equation}
    \langle N_{t_0}(\mathbf{x}_2),\mathbf{e} \rangle^{-1}=\frac{(1+\mu^2)\sqrt{1+u'(x_2,t_0)^2}}{\mu(\mu+u'(x_2,t_0))\langle \omega_1,\omega_2 \rangle+(1-\mu u'(x_2,t_0))\langle \phi_1,\phi_2 \rangle}. \label{es58}
\end{equation}
By \eqref{es57} and \eqref{es58}, there exists $0<\delta=\delta(p,q)<\frac{1}{20}$ s.t. if $|\mathbf{x}_1-\mathbf{x}_2|\leq \delta \rho,$ then
\begin{equation}
    \langle N_{t_0}(\mathbf{x}_2),\mathbf{e} \rangle^{-1}\leq \sqrt{1+(\frac{1}{3}\min\{\mu,\mu^{-1}\})^2}, \label{es60}
\end{equation}
namely, for any $\mathbf{x}_1=\mathbf{x}(u(x_1,t_0),\omega_1,\phi_1)\in \Sigma_{t_0}\cap (\bar{B}(0,\frac{3}{4})-B(0,\frac{1}{4}\rho))$ and any $\mathbf{x}_2\in \Sigma_{t_0}\cap \bar{B}(\mathbf{x}_1,\delta \rho)$, \eqref{es60} holds.
Then \eqref{es52} implies, there exists $0<\bar{\delta}(p,q)<\delta(p,q)$ s.t. for all $\mathbf{x}'\in \Sigma_t\cap \bar{B}(\mathbf{x}_1,\bar{\delta} \rho)$ and all $t_0\leq t\leq t_0+ \frac{(\bar{\delta}\rho)^2}{2(n-1)}$, the (outward) unit normal vector of $\Sigma_t$ at $\mathbf{x}'$ satisfies
\begin{equation}
    |\langle N_t(\mathbf{x}'),\mathbf{e} \rangle|^{-1}\leq \sqrt{1+(\frac{1}{2}\min\{\mu,\mu^{-1}\})^2}. \label{es59}
\end{equation}
If $|t_0|\leq \frac{(\bar{\delta}\rho)^2}{2(n-1)}$, then \eqref{es59} holds for all $t_0\leq t\leq \mathring{t}$.\\
By direct computation, the norm of the second fundamental form of $\Sigma_{t_0}$ at $\mathbf{x}=\mathbf{x}(u(x,t_0),\omega_1,\phi_1)$ is
\begin{equation}
    |A_{t_0}(\mathbf{x})|^2=\frac{1}{1+u'^2}((\frac{u''}{1+u'^2})^2+(p-1)(\frac{\mu+u'}{x-\mu u})^2+(q-1)(\frac{1-\mu u'}{\mu x+u})^2),\ u=u(x,t_0). \label{es62}
\end{equation}
According to \eqref{es61}, for $\mathbf{x}_1$ as above, if $\mathbf{x}_2\in \Sigma_{t_0}\cap \bar{B}(\mathbf{x}_1,\delta \rho)$,
$$|A_{t_0}(\mathbf{x}_2)|^2\leq \frac{C(p,q)}{\rho^2}.$$
Applying \eqref{es53} and \eqref{es60}, we get, for all $\mathbf{x}'\in \Sigma_t\cap \bar{B}(\mathbf{x}_1,\bar{\delta} \rho)$ and all $t_0\leq t\leq \mathring{t}$, 
\begin{equation}
    |A_t(\mathbf{x}')|^2\leq \frac{C(p,q)}{\rho^2}. \label{es63}
\end{equation}
On the other hand, $\Sigma_{t_0}-B(0,\frac{1}{2})$ is a compact hypersurface (with boundary) depending only on $p,q$,  so there exists $0<\varepsilon=\varepsilon(p,q)<\frac{1}{4}$ s.t. for every $\mathbf{x}_1\in \Sigma_{t_0}-B(0,\frac{3}{4})$, and $\mathbf{x}_2\in \Sigma_{t_0}\cap \bar{B}(\mathbf{x}_1,\varepsilon)$, then 
$$|A_{t_0}(\mathbf{x}_2)|^2\leq C(p,q),\ \langle N_{t_0}(\mathbf{x}_2),N_{t_0}(\mathbf{x}_1)\rangle^{-1}\leq C(p,q).$$
Again, \eqref{es52} implies, there exists $0<\bar{\varepsilon}(p,q)<\varepsilon(p,q)$ s.t. for all $\mathbf{x}'\in \Sigma_t\cap \bar{B}(\mathbf{x}_1,\bar{\varepsilon})$ and all $t_0\leq t\leq t_0+ \bar{\varepsilon}^2$, 
$$|A_t(\mathbf{x}')|^2\leq C(p,q).$$
Assuming $\bar{\delta}\rho\leq \bar{\varepsilon}$, this estimate holds for all $t_0\leq t\leq \mathring{t}$.\\
It's known that $\Sigma_t$ lies in the (closed) $\sqrt{2(n-1)(t-t_0)}$ neighborhood of $\Sigma_{t_0}$ from \cite{wan}, Corollary 2.1, and hence, for all $t_0\leq t\leq t_0+\frac{(\bar{\delta}\rho)^2}{2(n-1)}$, $\Sigma_t-B(0,\frac{1}{3}\rho)$ lies in the (closed) $\bar{\delta}\rho$ neighborhood of $\Sigma_{t_0}-B(0,(\frac{1}{3}-\bar{\delta})\rho)$. Since $|t_0|\leq \frac{(\bar{\delta}\rho)^2}{2(n-1)}$ and $\bar{\delta}<\frac{1}{12}$, for all $t_0\leq t\leq \mathring{t}$, $\Sigma_t-B(0,\frac{1}{3}\rho)$ lies in the (closed) $\bar{\delta}\rho$ neighborhood of $\Sigma_{t_0}-B(0,\frac{1}{4}\rho)$. Therefore, for any $t_0\leq t\leq \mathring{t}$, and $\mathbf{x}\in \Sigma_t-B(0,\frac{1}{3}\rho)$, there is $$|A_t(\mathbf{x})|\leq \frac{C(p,q)}{\rho},$$
which verifies (\ref{propes4_1}) of this Proposition.\\
To prove (\ref{propes4_2}), we first notice that for $x_j,y_j\geq 0$, $\omega_j\in \sph^{p-1}$, $\phi_j\in \sph^{q-1}$, $j=1,2$,
$$|(x_1\omega_1+y_1\phi_1)-(x_2\omega_2+y_2\phi_2)|^2=|x_1\omega_1-x_2\omega_2|^2+|y_1\phi_1-y_2\phi_2|^2$$
$$=(x_1-x_2)^2+2x_1x_2(1-\langle \omega_1,\omega_2 \rangle)+(y_1-y_2)^2+2y_1y_2(1-\langle \phi_1,\phi_2 \rangle)$$
$$\geq (x_1-x_2)^2+(y_1-y_2)^2=|(x_1,y_1)-(x_2,y_2)|^2.$$
From the discussion above, for all $t_0\leq t\leq \mathring{t}$, $\Sigma_t\cap (\bar{B}(0,\frac{1}{2})-B(0,\frac{1}{3}\rho))$ lies in the (closed) $\bar{\delta}\rho$ neighborhood of $\Sigma_{t_0}\cap (\bar{B}(0,\frac{3}{4})-B(0,\frac{1}{4}\rho))$. Taking the profile curve $\gamma_t^\mathbf{a}=\gamma_t$ of $\Sigma_t$ into account, we deduce that for all $t_0\leq t\leq \mathring{t}$, $\gamma_t\cap (\bar{B}^2(0,\frac{1}{2})-B^2(0,\frac{1}{3}\rho))\subset U$, where $U$ is the (closed) $\bar{\delta}\rho$ neighborhood of $\gamma_{t_0}\cap (\bar{B}^2(0,\frac{3}{4})-B^2(0,\frac{1}{4}\rho))$, and $B^2$ and $\bar{B}^2$ are the open and closed balls on the 2-dimensional plane respectively. If $\gamma_t$ is (locally) parameterized as \eqref{ad1}, with
$$(\frac{x-\mu u(x,t)}{\sqrt{1+\mu^2}},\frac{\mu x+u(x,t)}{\sqrt{1+\mu^2}})\in U,$$
then we may assume
$$|(\frac{x-\mu u(x,t)}{\sqrt{1+\mu^2}},\frac{\mu x+u(x,t)}{\sqrt{1+\mu^2}})-(\frac{x_0-\mu u(x_0,t_0)}{\sqrt{1+\mu^2}},\frac{\mu x_0+u(x_0,t_0)}{\sqrt{1+\mu^2}})|\leq \bar{\delta}\rho,\ \frac{1}{6}\rho\leq x_0\leq 1,$$
i.e.
$$|(x,u(x,t))-(x_0,u(x_0,t_0))|\leq \bar{\delta}\rho.$$
By \eqref{es61},
$$|\frac{u(x,t)}{x}|\leq \frac{|u(x_0,t_0)|+\bar{\delta}\rho}{x_0-\bar{\delta}\rho}\leq \frac{\frac{1}{4}\min\{\mu,\mu^{-1}\}x_0+\bar{\delta}\rho}{x_0-\bar{\delta}\rho}\leq \frac{\frac{1}{4}\min\{\mu,\mu^{-1}\}\frac{1}{6}\rho+\bar{\delta}\rho}{\frac{1}{6}\rho-\bar{\delta}\rho}.$$
If we choose $\bar{\delta}\ll 1$ (depending on $p,q$), then 
$$|\frac{u(x,t)}{x}|\leq \frac{1}{2}\min\{\mu,\mu^{-1}\}.$$
Moreover, by \eqref{es58} and \eqref{es59}, (letting $\omega_1=\omega_2$, $\phi_1=\phi_2$ in \eqref{es58},) 
$$|u'(x,t)|\leq \frac{1}{2}\min\{\mu,\mu^{-1}\},$$
and from \eqref{es62}, \eqref{es63} we know
$$|u''(x,t)|\leq |A_t(\mathbf{x})|(1+u'(x,t)^2)^{\frac{3}{2}}\leq \frac{C(p,q)}{\rho},\ \mathbf{x}=\mathbf{x}(u(x,t),\omega,\phi).$$
Finally, during the flow, the curve $\gamma_t\cap (\bar{B}^2(0,\frac{1}{2})-B^2(0,\frac{1}{3}\rho))$ never leaves $U\cap (\bar{B}^2(0,\frac{1}{2})-B^2(0,\frac{1}{3}\rho))$, and has a uniformly bounded gradient as a graph over the ray $l_{p,q}$ (see \eqref{mi25}). At the initial time $t_0$, $\gamma_{t_0}\cap (\bar{B}^2(0,\frac{1}{2})-B^2(0,\frac{1}{3}\rho))$ is written globally as a (single-valued) graph over $l_{p,q}$, so $\gamma_t\cap (\bar{B}^2(0,\frac{1}{2})-B^2(0,\frac{1}{3}\rho))$ is always a graph over $l_{p,q}$, for $t_0\leq t\leq \mathring{t}$. The proof of (\ref{propes4_2}) is now complete (provided $\rho\leq \frac{1}{6}$).

\QED

\section{Smooth Estimates and Determination of the Constant $\Lambda$}
In this section, we will first derive the higher order estimates in Propositions \ref{propco1} and \ref{propco2}, and then describe how to determine the constant $\Lambda$ defined in the admissible condition. The main tools come from the ``standard" theory of parabolic equations, including Schauder's estimates and H\"{o}lder continuity estimates for linear equations; see \cite{lie} for a reference.
\begin{prop}
    \label{propsm1}
    If \eqref{ad12} holds, then we have the estimates \eqref{co12}, \eqref{co13}, \eqref{co14} and \eqref{co15}.
\end{prop}
The key steps of proof are (taking \eqref{co13} for example): estimate the H\"{o}lder continuity of $v$ and its spatial derivative, use Schauder's theory to obtain a smooth estimate of $v$, and use Schauder's theory again together with $C^0$ estimates obtained in the last section to get the desired result. See Section 7 of \cite{gs} for details. One may use the estimates near the bottom in addition to the interior estimates, in order to let our results ``global" in time. Before deriving \eqref{co15}, a gradient estimate is required, which can be derived using maximum principle. The equation \eqref{ad10} is singular at $z=0$, so to obtain \eqref{co15} one should regard $\hat{w}$ as a radially symmetric function of $p$ variables and apply the corresponding estimates to it.\\
Until now, the proof of Proposition \ref{propco2} is complete. In order to justify \eqref{co17}, we need one more estimate for the function $u(x,t)$, $x\approx \rho$. The argument differs from that in \cite{gs}.
\begin{lemma}
    \label{lemsm2}
    If $|t_0|\ll 1$ (depending on $n,l,\rho,\beta$), then for $\frac{3}{4}\rho\leq x\leq \rho$, $t_0\leq t\leq \mathring{t}$, $i=0,1,2$,
    $$x^i|\partial_x^i u(x,t)|\leq C(p,q,l)x^{2\lambda_l+1}.$$
\end{lemma}
\proof
Rewrite the equation \eqref{ad2} as
\begin{equation}
    \partial_t u=\frac{1}{1+u'^2}u''+\frac{1}{x}P(\frac{u}{x})u'+\frac{1}{x^2}Q(\frac{u}{x})u, \label{sm13}
\end{equation}
where
\begin{equation}
    P(x)=\frac{(n-2)(1+(\mu^{-1}-\mu)x)}{(1-\mu x)(1+\mu^{-1}x)},\ Q(x)=\frac{n-2}{(1-\mu x)(1+\mu^{-1}x)}. \label{sm10}
\end{equation}
From our construction of initial value, \eqref{co2}, and \eqref{ad11}, we know if $|t_0|\ll 1$ (depending on $n,l,\rho,\beta$), then for $\frac{1}{2}\rho\leq x\leq \frac{5}{4}\rho$, $i=0,1,2,3$,
\begin{equation}
    |\partial_x^i u(x,t_0)|\leq C(n,l)\rho^{2\lambda_l+1-i}. \label{sm27}
\end{equation}
Putting \eqref{co11} (which holds for $x\geq \frac{1}{2}\rho$, by Proposition \ref{propes4}) into the equation \eqref{sm13}, we get for $\frac{1}{2}\rho\leq x\leq \frac{5}{4}\rho$, $t_0\leq t\leq \mathring{t}$,
$$|\partial_t u(x,t)|\leq \frac{C(p,q)}{\rho},$$
i.e.
$$|u(x,t)-u(x,t_0)|\leq \frac{C(p,q)}{\rho}|t-t_0|.$$
Assuming further that $|t_0|\ll 1$ (depending on $n,l,\rho$), for $(x,t)$ as above, there is (by \eqref{sm27})
$$|u(x,t)|\leq C(p,q,l)\rho^{2\lambda_l+1}.$$
Now, fix $X_*=(x_*,t_*)$ satisfying $t_0\leq t_*\leq \mathring{t},\ \frac{3}{4}\rho\leq x\leq \rho$. Using \eqref{co11} again, we get
$$|u|_{0,Q(X_*,\frac{1}{4}\rho)_{t_0}}\leq C(p,q,l)\rho^{2\lambda_l+1},\ |u'|_{0,Q(X_*,\frac{1}{4}\rho)_{t_0}}\leq C(p,q),$$
$$|\frac{1}{1+u'^2}|_{0,Q(X_*,\frac{1}{4}\rho)_{t_0}}+\rho|\frac{1}{x}P(\frac{u}{x})|_{0,Q(X_*,\frac{1}{4}\rho)_{t_0}}+\rho^2|\frac{1}{x^2}Q(\frac{u}{x})|_{0,Q(X_*,\frac{1}{4}\rho)_{t_0}}\leq C(n).$$
Also,
$$|u(\cdot,t_0)|_{3,B(x_*,\frac{1}{4}\rho)}^{(\rho)}\leq C(n,l)\rho^{2\lambda_l+1}.$$
Here $Q(X_*,R)_{t_0}=(B(x_*,R)\times (t_*-R^2,t_*))\cap \{t>t_0\}$, and the meaning of the norm is shown in the following example:
$$|u|_{2+\gamma,\Omega}^{(\rho)}=|u|_{0,\Omega}+\rho [u]_{1,\Omega}+\rho^2[u]_{2,\Omega}+\rho^{2+\gamma}[u]_{2+\gamma,\Omega},$$
where the subscript denotes a seminorm in which the exponents w.r.t. $x$ are twice those w.r.t. $t$, s.t. the norm is invariant under parabolic scaling.\\
Applying H\"{o}lder continuity estimate to \eqref{sm13}, we deduce there exists a universal constant $\gamma\in (0,1)$ s.t.
\begin{equation*}
    \rho^{\gamma}[u]_{\gamma,Q(X_*,\frac{1}{5}\rho)_{t_0}}\leq C(p,q,l)\rho^{2\lambda_l+1}.
\end{equation*}
According to the gradient H\"{o}lder continuity estimate, we may assume for the same $\gamma$,
\begin{equation*}
    \rho^{\gamma}[u']_{\gamma,Q(X_*,\frac{1}{5}\rho)_{t_0}}\leq C(p,q,l).
\end{equation*}
Thus,
$$|\frac{1}{1+u'^2}|_{\gamma,Q(X_*,\frac{1}{5}\rho)_{t_0}}^{(\rho)}+\rho|\frac{1}{x}P(\frac{u}{x})|_{\gamma,Q(X_*,\frac{1}{5}\rho)_{t_0}}^{(\rho)}+\rho^2|\frac{1}{x^2}Q(\frac{u}{x})|_{\gamma,Q(X_*,\frac{1}{5}\rho)_{t_0}}^{(\rho)}\leq C(p,q,l).$$
Applying Schauder's estimate to \eqref{sm13},
$$|u|_{2+\gamma,Q(X_*,\frac{1}{6}\rho)_{t_0}}^{(\rho)}\leq C(p,q,l)\rho^{2\lambda_l+1},$$
which obviously implies the Lemma.

\QED

Below we will eventually prove Proposition \ref{propco1}, as well as the main existence result Theorem \ref{thmco1}:
\begin{prop}
    \label{propsm4}
    If $\Lambda\gg 1$ (depending on $p,q,l$), $\rho\ll 1$, $\beta\gg 1$, $R\gg 1$ (depending on $p,q,l,\Lambda$), and $t_0\gg 1$ (depending on $p,q,l,\Lambda,\beta,R$), then Proposition \ref{propco1} holds.
\end{prop}
\proof
Actually, (\ref{propco1_1}) is exactly \eqref{es21}, and (\ref{propco1_2}) is shown in Proposition \ref{propes4}, (\ref{propes4_1}). \eqref{co11} is proved in Proposition \ref{propes4}, (\ref{propes4_2}), and \eqref{co18} is implied by \eqref{co15}, provided $\beta\gg 1$ (depending on $p,q,l$). The statement ``The profile curve of $\Sigma_t\cap (B(0,3\rho)-\bar{B}(0,\frac{1}{3}\beta(-t)^{\frac{1}{2}+\sigma_l}))$ can be parameterized by a single function as \eqref{ad1}" is a consequence of the admissible condition (\ref{ad_2}), (\ref{ad_3}), Proposition \ref{propes4}, (\ref{propes4_2}), and the $C^0$ and $C^1$ estimates in \eqref{co18}, and the statement ``The profile curve of $\Sigma_t\cap B(0,3\rho)$ can be parameterized by a single function as \eqref{ad3}" follows from the admissible condition (\ref{ad_2}), and Proposition \ref{propes4}, (\ref{propes4_2}), especially  \eqref{co11}.\\
Now, it remains to prove \eqref{co17}. By Lemma \ref{lemsm2}, if $\Lambda\gg 1$ (depending on $p,q,l$), then \eqref{co17} holds for $\frac{3}{4}\rho\leq x\leq \rho$, $t_0\leq t\leq \mathring{t}$.\\
By \eqref{ad11}, for $x\geq \sqrt{-t}$, $i=0,1,2$,
$$x^i|\partial_x^i(\frac{k}{c_l}(-t)^{\lambda_l+\frac{1}{2}}\varphi_l(\frac{x}{\sqrt{-t}}))|\leq C(n,l)x^{2\lambda_l+1}.$$
Thus by \eqref{co12}, if we choose $R\gg 1$, $\rho\ll 1$ (depending on $p,q,l,\Lambda$), then 
$$x^i|\partial_x^i u(x,t)|\leq C(n,l)x^{2\lambda_l+1},\ R\sqrt{-t}\leq x\leq \frac{3}{4}\rho,\ t_0\leq t\leq \mathring{t},\ i=0,1,2,$$
and \eqref{co17} holds for $(x,t)$ as above, provided $\Lambda\gg 1$ (depending on $n,l$).\\
Again, by \eqref{ad11}, for $y>0$, $i=0,1,2$,
$$y^i|\partial_y^i(\frac{k}{c_l}e^{-\lambda_l s}\varphi_l(y))|\leq C(n,l)e^{-\lambda_l s}(y^\alpha+y^{2\lambda_l+1}).$$
Thus by \eqref{co13}, if we choose $s_0\gg 1$ (depending on $p,q,l,\Lambda,\beta,R$), then 
$$y^i|\partial_y^i v(y,s)|\leq C(n,l)e^{-\lambda_l s}(y^\alpha+y^{2\lambda_l+1}),\ e^{-\frac{1}{2}\sigma_l s}\leq y\leq R,\ s_0\leq s\leq \mathring{s},\ i=0,1,2,$$
i.e.
$$x^i|\partial_x^i u(x,t)|\leq C(n,l)((-t)^l x^\alpha+x^{2\lambda_l+1}),\ (-t)^{\frac{1}{2}+\frac{1}{2}\sigma_l}\leq x\leq R\sqrt{-t},\ t_0\leq t\leq \mathring{t},\ i=0,1,2,$$
and \eqref{co17} holds for $(x,t)$ as above, provided $\Lambda\gg 1$ (depending on $n,l$). \\
By \eqref{mi15}, for $y>\hat{\psi}_k(0)\frac{\mu}{\sqrt{1+\mu^2}} e^{-\sigma_l s}$, $i=0,1,2$,
$$y^i|\partial_y^i(e^{-\sigma_l s}\psi_k(e^{\sigma_l s}y))|\leq C(p,q)e^{-\lambda_l s}y^\alpha.$$
Thus by \eqref{co14}, if we choose $\beta\gg 1$ (depending on $p,q,l,\Lambda$), then 
$$y^i|\partial_y^i v(y,s)|\leq C(p,q)e^{-\lambda_l s}y^\alpha,\ \frac{3}{2}\beta e^{-\sigma_l s}\leq y\leq e^{-\frac{1}{2}\sigma_l s},\ s_0\leq s\leq \mathring{s},\ i=0,1,2,$$
i.e.
$$x^i|\partial_x^i u(x,t)|\leq C(p,q)(-t)^l x^\alpha,\ \frac{3}{2}\beta(-t)^{\frac{1}{2}+\sigma_l}\leq x\leq (-t)^{\frac{1}{2}+\frac{1}{2}\sigma_l},\ t_0\leq t\leq \mathring{t},\ i=0,1,2,$$
and \eqref{co17} holds for $(x,t)$ as above, provided $\Lambda\gg 1$ (depending on $p,q$). \\
By \eqref{mi22}, for $\frac{\frac{1}{2}\beta}{\sqrt{1+\mu^2}}\leq z\leq \frac{2\beta}{\sqrt{1+\mu^2}}$, $i=0,1,2$,
$$z^i|\partial_z^i(\hat{\psi}_k(z)-\mu z)|\leq C(p,q)z^\alpha.$$
Thus by \eqref{co15}, if we choose $\beta\gg 1$ (depending on $p,q,l$), then 
$$z^i|\partial_z^i(\hat{w}(z,\tau)-\mu z)|\leq C(p,q)z^\alpha,\ \frac{\frac{1}{2}\beta}{\sqrt{1+\mu^2}}\leq z\leq \frac{2\beta}{\sqrt{1+\mu^2}},\ \tau_0\leq \tau\leq \mathring{\tau},\ i=0,1,2.$$
Using the linear transform
\begin{equation*}
    T(x,y)=(\frac{x+\mu y}{\sqrt{1+\mu^2}},\frac{-\mu x+y}{\sqrt{1+\mu^2}}),
\end{equation*}
we let
$$T(z,\hat{w}(z,\tau))=(z_1,w(z_1,\tau)),$$
then
$$w(z_1,\tau)=\frac{\hat{w}(z,\tau)-\mu z}{\sqrt{1+\mu^2}},\ w'(z_1,\tau)=\frac{\hat{w}'(z,\tau)-\mu}{\mu \hat{w}'(z,\tau)+1},\ w''(z_1,\tau)=(1+\mu^2)^{\frac{3}{2}}\frac{\hat{w}''(z,\tau)}{(\mu \hat{w}'(z,\tau)+1)^3},$$
and
$$z^i|\partial_z^iw(z_1,\tau)|\leq C(p,q)z^\alpha,\ \frac{\frac{1}{2}\beta}{\sqrt{1+\mu^2}}\leq z\leq \frac{2\beta}{\sqrt{1+\mu^2}},\ \tau_0\leq \tau\leq \mathring{\tau},\ i=0,1,2.$$
If $\beta\gg 1$ (depending on $p,q$), then $z_1\approx \sqrt{1+\mu^2}z$, and thus
$$z^i|\partial_z^iw(z,\tau)|\leq C(p,q)z^\alpha,\ \beta\leq z\leq \frac{3}{2}\beta,\ \tau_0\leq \tau\leq \mathring{\tau},\ i=0,1,2.$$
i.e.
$$x^i|\partial_x^i u(x,t)|\leq C(p,q)(-t)^l x^\alpha,\ \beta(-t)^{\frac{1}{2}+\sigma_l}\leq x\leq \frac{3}{2}\beta(-t)^{\frac{1}{2}+\sigma_l},\ t_0\leq t\leq \mathring{t},\ i=0,1,2,$$
and \eqref{co17} holds for $(x,t)$ as above, provided $\Lambda\gg 1$ (depending on $p,q$). The proof is finished.\\

\QED

\section{Vanishing Theorems of Parabolic Equations on Lawson's Cones and the Related Minimal Hypersurfaces}
In order to bound the mean curvature of the MCF solution we have constructed above near the singularity, we need a blow up argument as in \cite{sto}, which further requires some ``vanishing theorems" for solutions to a kind of linear parabolic equations like $\partial_t u=(\Delta +|A|^2)u$ on the limit spaces (Lawson's cone $C_{p,q}$ and the related minimal hypersurface $\mathcal{M}_k$), with certain growth control. Actually the right-hand side of the equation is the so-called ``Jacobi operator" of hypersurfaces, which is closely related to the stability of such minimal hypersurfaces. We present the following results, which can be proved similarly to those in \cite{sto}:
\begin{theorem}
    \label{thmva1}
    Let $u=u(|\mathbf{x}|,t)$ be a smooth, radially symmetric, ancient solution to
    $$\partial_t u=\Delta_{\mathcal{M}_k}u+|A_{\mathcal{M}_k}|^2 u,\ (\mathbf{x},t)\in \mathcal{M}_k\times (-\infty,0]$$
    for some $k>0$, where $\Delta_{\mathcal{M}_k}$, $A_{\mathcal{M}_k}$ denote the Laplacian and the second fundamental form of $\mathcal{M}_k$ respectively. If there exists $C>0$ and $\delta>0$ s.t.
    $$|u(\mathbf{x},t)|\leq C(1+|\mathbf{x}|)^{\alpha-\delta},\ (\mathbf{x},t)\in \mathcal{M}_k\times (-\infty,0],$$
    then $u\equiv 0$.
\end{theorem}

\begin{theorem}
    \label{thmva2}
     Let $u=u(|\mathbf{x}|,t)$ be a smooth, radially symmetric, ancient solution to
    $$\partial_t u=\Delta_{C_{p,q}}u+|A_{C_{p,q}}|^2 u,\ (\mathbf{x},t)\in (C_{p,q}-\{0\})\times (-\infty,0].$$
    If there exists $C>0$ and $0<\delta<n-3+2\alpha$ s.t.
    $$|u(\mathbf{x},t)|\leq C|\mathbf{x}|^{\alpha-\delta},\ (\mathbf{x},t)\in (C_{p,q}-\{0\})\times (-\infty,0],$$
    then $u\equiv 0$.
\end{theorem}

\section{Boundedness of the Mean Curvature}
In this section, we will show the mean curvature of Vel\'{a}zquez's solution $\{\Sigma_t^{\mathbf{a}}\}=\{\Sigma_t\}$, obtained in Theorem \ref{thmco1}, remains bounded up to the singular time $t=0$, provided the parameter $l$ is sufficiently large. In the outer region $|\mathbf{x}|\geq \sqrt{-t}$, the boundedness of mean curvature is in fact a direct consequence of the estimates obtained in Proposition \ref{propco1}.

\begin{prop}
    \label{propbo1}
    If $\rho\ll 1$ (depending on $p,q,l,\Lambda$), then
    $$\sup_{\begin{subarray}{c}
        \mathbf{x}\in \Sigma_t-B(0,\sqrt{-t})\\
        t_0\leq t<0
    \end{subarray}} |H_{\Sigma_t}(\mathbf{x})|<+\infty,$$
    where $H_{\Sigma_t}(\mathbf{x})$ is the mean curvature of $\Sigma_t$ at $\mathbf{x}$.
\end{prop}
\proof
First of all, by (\ref{propco1_2}) of Proposition \ref{propco1}, 
$$\sup_{\begin{subarray}{c}
        \mathbf{x}\in \Sigma_t-\bar{B}(0,2\rho)\\
        t_0\leq t<0
    \end{subarray}} |H_{\Sigma_t}(\mathbf{x})|\leq C(p,q,\rho)<+\infty.$$
To estimate the remaining part, it's not hard to compute directly that if a planar curve has the form \eqref{ad1}, then the mean curvature (w.r.t. the upward unit vector) of the corresponding $O(p)\times O(q)$ invariant hypersurface at $(\frac{x-\mu u(x,t)}{\sqrt{1+\mu^2}}\omega,\frac{\mu x+u(x,t)}{\sqrt{1+\mu^2}}\phi)$ ($\omega\in \sph^{p-1}$, $\phi\in \sph^{q-1}$), which is denoted by $H(u(x,t))$, is
$$H(u(x,t))=\frac{1}{\sqrt{1+u'^2}}(\frac{u''}{1+u'^2}+(p-1)\frac{\mu+u'}{x-\mu u}-(q-1)\frac{1-\mu u'}{\mu x+u})$$
\begin{equation}
    =\frac{1}{\sqrt{1+u'^2}}(\frac{u''}{1+u'^2}+P(\frac{u}{x})\frac{u'}{x}+Q(\frac{u}{x})\frac{u}{x^2}), \label{bo4}
\end{equation}
where $P,Q$ are defined in \eqref{sm10}. For the part of the profile curve \eqref{ad1} lying in $B^2(0,3\rho)$, if $x\geq \rho$, then \eqref{co11} gives a bound of $H(u(x,t))$. If $\frac{1}{2}\sqrt{-t}\leq x\leq \rho$, then by \eqref{co17}, we also have $|\frac{u(x,t)}{x}|\leq \frac{1}{2}\min\{\mu,\mu^{-1}\}$, provided $\rho\ll 1$ (depending on $p,q,l,\Lambda$), and thus $P(\frac{u}{x}),Q(\frac{u}{x})$ are bounded. Then \eqref{co17} also gives a bound of $H(u(x,t))$. Since the part of the profile curve \eqref{ad1} lying in $B^2(0,3\rho)-B^2(0,\sqrt{-t})$ is covered by the two cases above, the proof is complete.

\QED

To bound the mean curvature in the intermediate and tip regions, we employ a blow up argument as in \cite{sto}, which in addition requires some convergence results of certain rescaled flows. They are presented in the following two lemmas, which can be derived from the estimates in in Proposition \ref{propco2}. For two sequences of real numbers $\{a_i\}$, $\{b_i\}$, we write $a_i\ll b_i$ if $a_i=o(b_i)$ as $i\rightarrow +\infty$.

\begin{lemma}
    \label{lembo1}
    For any sequence $t_0\leq t_i<0$, $t_i\nearrow 0$, and $\Lambda_i=(-t_i)^{-\frac{1}{2}-\sigma_l}$, the sequence of flows
    $$\tilde{\Sigma}_\tau^i:=\Lambda_i \Sigma_{t_i+\frac{\tau}{\Lambda_i^2}},\ (t_0-t_i)\Lambda_i^2\leq \tau< -t_i\Lambda_i^2$$
    converges smoothly to $\mathcal{M}_k$ (a stationary MCF) in any bounded space-time region in $\R^n\times \R$.
\end{lemma}
\proof
By Proposition \ref{propco1}, $\tilde{\Sigma}_\tau^i\cap B(0,3\rho \Lambda_i)$ can be parameterized as $(x,\hat{u}_i(x,\tau))$ (see \eqref{ad3}), and $\tilde{\Sigma}_\tau^i\cap (B(0,3\rho \Lambda_i)-\bar{B}(0,\frac{1}{3}\beta(-t_i- \Lambda_i^{-2}\tau)^{\frac{1}{2}+\sigma_l}\Lambda_i))$ can be parameterized as $(\frac{x-\mu u_i(x,t)}{\sqrt{1+\mu^2}},\frac{\mu x+u_i(x,t)}{\sqrt{1+\mu^2}})$ (see \eqref{ad1}). Thus it suffices to consider the convergence of $\hat{u}_i$ and $u_i$.\\
By the definition of $\tilde{\Sigma}_\tau^i$ and \eqref{ad9}, 
$$\hat{u}_i(x,\tau)=\Lambda_i \hat{u}(\Lambda_i^{-1} x,t_i+\Lambda_i^{-2}\tau)=\Lambda_i (-t_i-\Lambda_i^{-2}\tau)^{\frac{1}{2}+\sigma_l}\hat{w}(\Lambda_i^{-1}(-t_i-\Lambda_i^{-2}\tau)^{-(\frac{1}{2}+\sigma_l)} x,\frac{1}{2\sigma_l}(-t_i-\Lambda_i^{-2}\tau)^{-2\sigma_l})$$
$$=(1-\tau(-t_i)^{2\sigma_l})^{\frac{1}{2}+\sigma_l}\hat{w}((1-\tau(-t_i)^{2\sigma_l})^{-(\frac{1}{2}+\sigma_l)}x,\frac{1}{2\sigma_l}(-t_i)^{-2\sigma_l}(1-\tau(-t_i)^{2\sigma_l})^{-2\sigma_l}).$$
According to \eqref{co15}, for any $M>0$, $\hat{w}(x,\frac{1}{2\sigma_l}(-t_i)^{-2\sigma_l}+\tau)$ converges smoothly to $\hat{\psi}_k(x)$ on $0\leq x\leq \frac{2\beta}{\sqrt{1+\mu^2}}$, $|\tau|\leq M$. Since $1-\tau(-t_i)^{2\sigma_l}$ and all its powers converge smoothly to 1 on any bounded $\tau$-interval, and $\frac{1}{2\sigma_l}(-t_i)^{-2\sigma_l}(1-\tau(-t_i)^{2\sigma_l})^{-2\sigma_l}-\frac{1}{2\sigma_l}(-t_i)^{-2\sigma_l}$ converges smoothly to $\tau$ on any bounded $\tau$-interval, for any $0<\epsilon<1<M$, $\hat{u}_i(x,\tau)$ converges smoothly to $\hat{\psi}_k(x)$ on $[0,\frac{2\beta}{\sqrt{1+\mu^2}}(1-\epsilon)]\times [-M,M]$.\\
Similarly, 
$$u_i(x,\tau)=(1-\tau(-t_i)^{2\sigma_l})^{\frac{1}{2}+\sigma_l}w((1-\tau(-t_i)^{2\sigma_l})^{-(\frac{1}{2}+\sigma_l)}x,\frac{1}{2\sigma_l}(-t_i)^{-2\sigma_l}(1-\tau(-t_i)^{2\sigma_l})^{-2\sigma_l}).$$
According to \eqref{co14}, 
$$z^{m+2r}|\partial_z^m \partial_\tau ^r (w(z,\tau)-\psi_k(z))|\leq C(p,q,l,\Lambda,m,r)\beta^{\tilde{\alpha}-\alpha}(\frac{\tau}{\tau_0})^{-\varrho}z^\alpha$$
for $\frac{3}{2}\beta\leq z\leq (2\sigma_l \tau)^{\frac{1}{4}}$, $\tau\geq \tau_0$, $m,r\geq 0$. Due to the same reason as above, $u_i(x,\tau)$ converges smoothly to $\psi_k(x)$ on $[\frac{3}{2}\beta(1+\epsilon),M]\times [-M,M]$. Therefore, the profile curve of $\tilde{\Sigma}_\tau^i$ converges smoothly to the profile curve of $\mathcal{M}_k$ in any bounded space-time region. Moreover, since the even extensions of $\hat{u}_i(\cdot,\tau)$ and $\hat{\psi}_k$ are smooth, the flow $\tilde{\Sigma}_\tau^i$ also converges to $\mathcal{M}_k$ in the desired way.

\QED

\begin{lemma}
    \label{lembo2}
    For any sequence $t_0\leq t_i<0$, $t_i\nearrow 0$, and any sequence $\{\Lambda_i\}$ s.t. $(-t_i)^{-\frac{1}{2}}\ll \Lambda_i\ll (-t_i)^{-\frac{1}{2}-\sigma_l}$, the sequence of flows
    $$\tilde{\Sigma}_\tau^i:=\Lambda_i \Sigma_{t_i+\frac{\tau}{\Lambda_i^2}},\ (t_0-t_i)\Lambda_i^2\leq \tau< -t_i\Lambda_i^2$$
    converges smoothly to $C_{p,q}$ (a stationary MCF) in any compact space-time region in $(\R^n-\{0\})\times \R$.
\end{lemma}
\proof
By Proposition \ref{propco1}, $\tilde{\Sigma}_\tau^i\cap (B(0,3\rho \Lambda_i)-\bar{B}(0,\frac{1}{3}\beta(-t_i- \Lambda_i^{-2}\tau)^{\frac{1}{2}+\sigma_l}\Lambda_i))$ can be parameterized as $(\frac{x-\mu u_i(x,t)}{\sqrt{1+\mu^2}},\frac{\mu x+u_i(x,t)}{\sqrt{1+\mu^2}})$ (see \eqref{ad1}). Thus it suffices to consider the convergence of $u_i$.\\
By \eqref{co13} and \eqref{ad11}, for $e^{-\frac{1}{2}\sigma_l s}\leq y\leq 1$, $s\geq s_1\gg s_0$, $m,r\geq 0$,
$$y^{m+2r} |\partial_y^m \partial_s^r v(y,s)|\leq C(n,l,m,r)e^{-\lambda_l s} y^\alpha.$$
Set $\bar{v}(y,s)=e^{-\sigma_l s}\psi_k(e^{\sigma_l s}y)$. According to \eqref{mi15}, in the domain of $\bar{v}$, $y\geq \hat{\psi}_k(0)\frac{\mu}{\sqrt{1+\mu^2}}e^{-\sigma_l s}$, there holds for all integers $m,r\geq 0:$ (note $k\approx 1$):
\begin{equation}
    |\partial_y^m \partial_s^r \bar{v}(y,s)|\leq C(p,q,m,r)e^{-\lambda_l s}y^{\alpha-m}. \label{sm39}
\end{equation}
Now \eqref{co14}, \eqref{sm39} tell us, for $\frac{3}{2}\beta e^{-\sigma_l s}\leq y\leq e^{-\frac{1}{2}\sigma_l s}$, $s\geq s_0$,
\begin{equation}
    y^{m+2r} |\partial_y^m \partial_s^r v(y,s)|\leq C(p,q,l,\Lambda,m,r)e^{-\lambda_l s} y^\alpha. \label{bo1}
\end{equation}
Thus, \eqref{bo1} holds for $\frac{3}{2}\beta e^{-\sigma_l s}\leq y\leq 1$, $s\geq s_1$. By the definition of $\tilde{\Sigma}_\tau^i$ and \eqref{ad6}, 
$$u_i(x,\tau)=\Lambda_i u(\Lambda_i^{-1} x,t_i+\Lambda_i^{-2}\tau)=\Lambda_i (-t_i-\Lambda_i^{-2}\tau)^{\frac{1}{2}}v(\Lambda_i^{-1}(-t_i-\Lambda_i^{-2}\tau)^{-\frac{1}{2}} x,-\ln(-t_i-\Lambda_i^{-2}\tau))$$
$$=(-t_i\Lambda_i^2-\tau)^{\frac{1}{2}}v((-t_i\Lambda_i^2-\tau)^{-\frac{1}{2}} x,-\ln(-t_i-\Lambda_i^{-2}\tau)).$$
Therefore, if $\frac{3}{2}\beta (-t_i-\Lambda_i^{-2}\tau)^{\frac{1}{2}+\sigma_l}\Lambda_i\leq x\leq (-t_i\Lambda_i^2-\tau)^{\frac{1}{2}}$, $\tau\geq \Lambda_i^2(-e^{-s_1}-t_i)$, 
$$x^{m+2r}|\partial_x^m \partial_\tau^r u_i(x,\tau)|\leq C(p,q,l,\Lambda,m,r) (-t_i-\Lambda_i^{-2}\tau)^{(1-\alpha)(\frac{1}{2}+\sigma_l)}\Lambda_i^{1-\alpha}x^\alpha.$$
Since $(-t_i)^{-\frac{1}{2}}\ll \Lambda_i\ll (-t_i)^{-\frac{1}{2}-\sigma_l}$, on any bounded $\tau$-interval, $(-t_i-\Lambda_i^{-2}\tau)^{(1-\alpha)(\frac{1}{2}+\sigma_l)}\Lambda_i^{1-\alpha}\rightarrow 0$ uniformly, and $(-t_i-\Lambda_i^{-2}\tau)^{\frac{1}{2}+\sigma_l}\Lambda_i\rightarrow 0$, $(-t_i\Lambda_i^2-\tau)^{\frac{1}{2}}\rightarrow +\infty$ uniformly. Also, $\Lambda_i^2(-e^{-s_1}-t_i)\rightarrow -\infty$, which implies for any $M>1$, $u_i(x,\tau)\rightarrow 0$ smoothly on $[M^{-1},M]\times [-M,M]$, and the flow $\tilde{\Sigma}_\tau^i$ also converges to $C_{p,q}$ in the desired way.

\QED

The next proposition bounds the mean curvature in the intermediate and tip regions:
\begin{prop}
    \label{propbo2}
    If $-\alpha<a<1-\alpha$, and 
    \begin{equation}
        \lambda_l(1-\frac{a}{1-\alpha})-\frac{1}{2}\geq 0, \label{bo6}
    \end{equation}
    then
    $$\sup_{\begin{subarray}{c}
        \mathbf{x}\in \Sigma_t\cap \bar{B}(0,\sqrt{-t})\\
        t_0\leq t<0
    \end{subarray}} (1+(-t)^{-(\frac{1}{2}+\sigma_l)}|\mathbf{x}|)^a |H_{\Sigma_t}(\mathbf{x})|<+\infty.$$
\end{prop}
\proof
Suppose our assertion is not true, then there exists a sequence $t_0\leq T_i\nearrow 0$ s.t.
$$M_i:=\sup_{\begin{subarray}{c}
       \mathbf{x}\in \Sigma_t\cap \bar{B}(0,\sqrt{-t})\\
        t_0\leq t\leq T_i
    \end{subarray}} (1+(-t)^{-(\frac{1}{2}+\sigma_l)}|\mathbf{x}|)^a |H_{\Sigma_t}(\mathbf{x})|\rightarrow +\infty,$$
and thus we can find $t_i\in [t_0,T_i]$ and $\mathbf{x}_i\in \Sigma_{t_i}\cap \bar{B}(0,\sqrt{-t_i})$ s.t. 
\begin{equation}
    (1+(-t_i)^{-(\frac{1}{2}+\sigma_l)}|\mathbf{x}_i|)^a |H_{\Sigma_{t_i}}(\mathbf{x}_i)|=M_i. \label{bo5}
\end{equation}
After passing to a subsequence, there are several possibilities of the behabior of the sequence $(\mathbf{x}_i,t_i)$ as $i\rightarrow +\infty$:
\begin{enumerate}
    \item $|\mathbf{x}_i|=O((-t_i)^{\frac{1}{2}+\sigma_l})$.
    \item $(-t_i)^{\frac{1}{2}+\sigma_l}\ll |\mathbf{x}_i|\ll \sqrt{-t_i}$.
    \item $|\mathbf{x}_i|\sim \sqrt{-t_i}$ (i.e. $|\mathbf{x}_i|=O(\sqrt{-t_i})$ and $\sqrt{-t_i}=O(|\mathbf{x}_i|)$).
\end{enumerate}
When (1) happens, define a sequence of rescaled flows as
$$\tilde{\Sigma}_\tau^i:=\Lambda_i \Sigma_{t_i+\frac{\tau}{\Lambda_i^2}},\ (t_0-t_i)\Lambda_i^2\leq \tau\leq 0,$$
where $\Lambda_i=(-t_i)^{-\frac{1}{2}-\sigma_l}$. Then
$$M_i=\sup_{\begin{subarray}{c}
       \mathbf{x}\in \Sigma_t\cap \bar{B}(0,\sqrt{-t})\\
        t_0\leq t\leq t_i
    \end{subarray}} (1+(-t)^{-(\frac{1}{2}+\sigma_l)}|\mathbf{x}|)^a |H_{\Sigma_t}(\mathbf{x})|$$
$$=\sup_{\begin{subarray}{c}
       \mathbf{y}\in \tilde{\Sigma}_\tau^i\cap \bar{B}(0,\sqrt{-t_i\Lambda_i^2-\tau})\\
        (t_0-t_i)\Lambda_i^2\leq \tau\leq 0
    \end{subarray}} (1+(-t_i-\Lambda_i^{-2}\tau)^{-(\frac{1}{2}+\sigma_l)}\Lambda_i^{-1}|\mathbf{y}|)^a \Lambda_i |H_{\tilde{\Sigma}_\tau^i}(\mathbf{y})|.$$
Now, define functions $u_i:\tilde{\Sigma}_\tau^i\rightarrow \R$ as
$$u_i(\mathbf{y},\tau)=\frac{\Lambda_i}{M_i}H_{\tilde{\Sigma}_\tau^i}(\mathbf{y})$$
(the sign of $H$ can be chosen, e.g. corresponding to the outer unit normal vector). According to \cite{eh}, Lemma 1.1, (iv), along the mean curvature flow (whose velocity equals the mean curvature), $u_i$ satisfies
$$\partial_\tau u_i=\Delta_{\tilde{\Sigma}_\tau^i}u_i+|A_{\tilde{\Sigma}_\tau^i}|^2 u_i,$$
and
$$(1+(-t_i-\Lambda_i^{-2}\tau)^{-(\frac{1}{2}+\sigma_l)}\Lambda_i^{-1}|\mathbf{y}|)^a |u_i(\mathbf{y},\tau)|\leq 1,\ \mathbf{y}\in \tilde{\Sigma}_\tau^i\cap \bar{B}(0,\sqrt{-t_i\Lambda_i^2-\tau}),\ (t_0-t_i)\Lambda_i^2\leq \tau\leq 0,$$
$$(1+|\mathbf{y}_i|)^{a}|u_i(\mathbf{y}_i,0)|=1 (\mathbf{y}_i=\Lambda_i \mathbf{x}_i).$$
Since $\Lambda_i^{-1}\ll \sqrt{-t_i}$, as $i\rightarrow +\infty$, $(t_0-t_i)\Lambda_i^2\rightarrow -\infty$, $\sqrt{-t_i\Lambda_i^2-\tau}\rightarrow +\infty$, $(-t_i-\Lambda_i^{-2}\tau)^{-(\frac{1}{2}+\sigma_l)}\Lambda_i^{-1}\rightarrow 1$ uniformly on compact $\tau$-intervals, and in particular $\{u_i\}$ are uniformly bounded in any bounded space-time region. By Lemma \ref{lembo1}, the sequence of flows $\{\tilde{\Sigma}_\tau^i\}$ converges to the stationary flow $\mathcal{M}_k$ locally smoothly, and thus, using Schauder's estimate, we may assume (by passing to a subsequence if necessary) that $\{u_i\}$ converges locally smoothly to a smooth ancient solution $u_\infty$ to
$$\partial_\tau u_\infty=\Delta_{\mathcal{M}_k}u_\infty+|A_{\mathcal{M}_k}|^2 u_\infty,$$
defined on $\mathcal{M}_k\times (-\infty,0]$, satisfying
\begin{equation}
    |u_\infty(\mathbf{y},\tau)|\leq (1+|\mathbf{y}|)^{-a} \label{bo2}
\end{equation}
on this domain. Since $\mathbf{y}_i=\Lambda_i \mathbf{x}_i$ are uniformly bounded, we may also assume $\mathbf{y}_i\rightarrow \mathbf{y}_\infty\in \mathcal{M}_k$, with
$$|u_\infty(\mathbf{y}_\infty,0)|=(1+|\mathbf{y}_\infty|)^{-a}>0.$$
But \eqref{bo2} and Theorem \ref{thmva1} imply $u_\infty \equiv 0$ (note that $-a<\alpha$), which is impossible.\\
When (2) happens, define a sequence of rescaled flows as
$$\tilde{\Sigma}_\tau^i:=|\mathbf{x}_i|^{-1}\Sigma_{t_i+|\mathbf{x}_i|^2 \tau},\ (t_0-t_i)|\mathbf{x}_i|^{-2}\leq \tau\leq 0.$$
Then
$$M_i=\sup_{\begin{subarray}{c}
       \mathbf{x}\in \Sigma_t\cap \bar{B}(0,\sqrt{-t})\\
        t_0\leq t\leq t_i
    \end{subarray}} (1+(-t)^{-(\frac{1}{2}+\sigma_l)}|\mathbf{x}|)^a |H_{\Sigma_t}(\mathbf{x})|$$
$$=\sup_{\begin{subarray}{c}
       \mathbf{y}\in \tilde{\Sigma}_\tau^i\cap \bar{B}(0,\sqrt{-t_i|\mathbf{x}_i|^{-2}-\tau})\\
        (t_0-t_i)|\mathbf{x}_i|^{-2}\leq \tau\leq 0
    \end{subarray}} (1+(-t_i-|\mathbf{x}_i|^2\tau)^{-(\frac{1}{2}+\sigma_l)}|\mathbf{x}_i|\ |\mathbf{y}|)^a |\mathbf{x}_i|^{-1} |H_{\tilde{\Sigma}_\tau^i}(\mathbf{y})|$$
$$\geq \sup_{\begin{subarray}{c}
       \mathbf{y}\in \tilde{\Sigma}_\tau^i\cap \bar{B}(0,\sqrt{-t_i|\mathbf{x}_i|^{-2}-\tau})\\
        (t_0-t_i)|\mathbf{x}_i|^{-2}\leq \tau\leq 0
    \end{subarray}} (\Lambda_i^{-1}(-t_i-|\mathbf{x}_i|^2\tau)^{-(\frac{1}{2}+\sigma_l)})^a \Lambda_i^a |\mathbf{x}_i|^{a-1}|\mathbf{y}|^a |H_{\tilde{\Sigma}_\tau^i}(\mathbf{y})|,$$
where $\Lambda_i=(-t_i)^{-\frac{1}{2}-\sigma_l}$. Now, define functions $u_i:\tilde{\Sigma}_\tau^i\rightarrow \R$ as
$$u_i(\mathbf{y},\tau)=\frac{\Lambda_i^a |\mathbf{x}_i|^{a-1}}{M_i}H_{\tilde{\Sigma}_\tau^i}(\mathbf{y}).$$
Again, along the mean curvature flow (whose velocity equals the mean curvature), $u_i$ satisfies
$$\partial_\tau u_i=\Delta_{\tilde{\Sigma}_\tau^i}u_i+|A_{\tilde{\Sigma}_\tau^i}|^2 u_i,$$
and
$$(\Lambda_i^{-1}(-t_i-|\mathbf{x}_i|^2\tau)^{-(\frac{1}{2}+\sigma_l)})^a |\mathbf{y}|^a |u_i(\mathbf{y},\tau)|\leq 1,\ \mathbf{y}\in \tilde{\Sigma}_\tau^i\cap \bar{B}(0,\sqrt{-t_i|\mathbf{x}_i|^{-2}-\tau}),\ (t_0-t_i)|\mathbf{x}_i|^{-2}\leq \tau\leq 0,$$
$$(1+|\mathbf{x}_i|^{-1}\Lambda_i^{-1})^{a}|u_i(\mathbf{y}_i,0)|=1 (\mathbf{y}_i=|\mathbf{x}_i|^{-1} \mathbf{x}_i).$$
Since $\Lambda_i^{-1}\ll |\mathbf{x}_i|\ll \sqrt{-t_i}$, as $i\rightarrow +\infty$, $(t_0-t_i)|\mathbf{x}_i|^{-2}\rightarrow -\infty$, $(1+|\mathbf{x}_i|^{-1}\Lambda_i^{-1})^{a}\rightarrow 1$, $\sqrt{-t_i|\mathbf{x}_i|^{-2}-\tau}\rightarrow +\infty$, $(\Lambda_i^{-1}(-t_i-|\mathbf{x}_i|^2\tau)^{-(\frac{1}{2}+\sigma_l)})^a\rightarrow 1$ uniformly on compact $\tau$-intervals, and in particular $\{u_i\}$ are uniformly bounded in any compact space-time region in $(\R^n-\{0\})\times (-\infty,0]$. By Lemma \ref{lembo2}, the sequence of flows $\{\tilde{\Sigma}_\tau^i\}$ converges to the stationary flow $C_{p,q}$ locally smoothly, and thus, using Schauder's estimate, we may assume (by passing to a subsequence if necessary) that $\{u_i\}$ converges locally smoothly to a smooth ancient solution $u_\infty$ to
$$\partial_\tau u_\infty=\Delta_{C_{p,q}}u_\infty+|A_{C_{p,q}}|^2 u_\infty,$$
defined on $(C_{p,q}-\{0\})\times (-\infty,0]$, satisfying
\begin{equation}
    |u_\infty(\mathbf{y},\tau)|\leq |\mathbf{y}|^{-a} \label{bo3}
\end{equation}
on this domain. Since $u_i,u_\infty$ are radially symmetric, for all $\mathbf{y}_\infty\in C_{p,q}$, $|\mathbf{y}_\infty|=1$,
$$|u_\infty(\mathbf{y}_\infty,0)|=1>0.$$
But \eqref{bo3} and Theorem \ref{thmva2} imply $u_\infty \equiv 0$ (note that $\alpha-(n-3+2\alpha)\leq \alpha-1<-a<\alpha$), which is impossible.\\
When (3) happens, then there exists $0<\epsilon<1$ s.t. for all $i$ large enough, $\mathbf{x}_i=(\frac{x_i-\mu u(x_i,t_i)}{\sqrt{1+\mu^2}}\omega_i,\frac{\mu x_i+u(x_i,t_i)}{\sqrt{1+\mu^2}}\phi_i)$, with $\omega_i\in \sph^{p-1}$, $\phi_i\in \sph^{q-1}$, $\epsilon\sqrt{-t_i}\leq x_i\leq \sqrt{-t_i}$. By the formula \eqref{bo4} and the estimate \eqref{co17}, 
$$|H_{\Sigma_{t_i}}(\mathbf{x}_i)|\leq C(-t_i)^{\lambda_l-\frac{1}{2}},$$
where $C>0$ is independent of $i$. Thus,
$$(1+(-t_i)^{-(\frac{1}{2}+\sigma_l)}|\mathbf{x}_i|)^a |H_{\Sigma_{t_i}}(\mathbf{x}_i)|\leq C(-t_i)^{\lambda_l-\frac{1}{2}-a\sigma_l}\leq C<+\infty,$$
since $\lambda_l-\frac{1}{2}-a\sigma_l\geq 0$, which contradicts to \eqref{bo5}.\\
Because all the possible cases are impossible, our assertion follows.

\QED

\begin{cor}
\label{corbo1}
    If $n\geq 9$, $l\geq 3$, or $n=8$, $l\geq 4$, then
    $$\sup_{\begin{subarray}{c}
        \mathbf{x}\in \Sigma_t\cap \bar{B}(0,\sqrt{-t})\\
        t_0\leq t<0
    \end{subarray}} |H_{\Sigma_t}(\mathbf{x})|<+\infty.$$
\end{cor}
\proof
The condition \eqref{bo6} is equivalent to
$$a\leq (1-\alpha)(1-\frac{1}{2\lambda_l}).$$
Since we also require $-\alpha<a<1-\alpha$, such $a$ exists iff
$$-\alpha<(1-\alpha)(1-\frac{1}{2\lambda_l})\Leftrightarrow\lambda_l>\frac{1-\alpha}{2}\Leftrightarrow l>1-\alpha,$$
which is true iff $n\geq 9$, $l\geq 3$, or $n=8$, $l\geq 4$. Choosing $a$ satisfying the conditions above and noting $(1+(-t)^{-(\frac{1}{2}+\sigma_l)}|\mathbf{x}|)^a\geq 1$, this corollary is a direct consequence of Proposition \ref{propbo2}.

\QED

The main Theorem \ref{thmin1} follows from Theorem \ref{thmco1}, Proposition \ref{propbo1}, and Corollary \ref{corbo1}.

\begin{remark}
    \label{rembo1}
    Using the blow up argument in Proposition \ref{propbo2}, one can prove for any $\epsilon>0$,
    $$\sup_{\mathbf{x}\in \Sigma_t,t_0\leq t<0}(-t)^{\frac{1}{2}-\sigma_l+\epsilon} |H_{\Sigma_t}(\mathbf{x})|<+\infty$$
    if $l=2$, or
    $$\sup_{\mathbf{x}\in \Sigma_t,t_0\leq t<0}(-t)^{\epsilon} |H_{\Sigma_t}(\mathbf{x})|<+\infty$$
    if $n=8$, $l=3$. Actually, in order to derive a contradiction, the assertion in Proposition \ref{propbo2} should be
    $$\sup_{\begin{subarray}{c}
        \mathbf{x}\in \Sigma_t\cap \bar{B}(0,\sqrt{-t})\\
        t_0\leq t<0
    \end{subarray}} (-t)^{\frac{1}{2}-\sigma_l+\epsilon}(1+(-t)^{-(\frac{1}{2}+\sigma_l)}|\mathbf{x}|)^a |H_{\Sigma_t}(\mathbf{x})|<+\infty$$
     with $-\alpha< a\leq (1-\alpha)(1-\lambda_l^{-1}(\sigma_l-\epsilon))$ if $l=2$, or
     $$\sup_{\begin{subarray}{c}
        \mathbf{x}\in \Sigma_t\cap \bar{B}(0,\sqrt{-t})\\
        t_0\leq t<0
    \end{subarray}} (-t)^{\epsilon}(1+(-t)^{-(\frac{1}{2}+\sigma_l)}|\mathbf{x}|)^a |H_{\Sigma_t}(\mathbf{x})|<+\infty$$
    with $-\alpha< a\leq (1-\alpha)(1-\lambda_l^{-1}(\frac{1}{2}-\epsilon))$ if $n=8$, $l=3$.
\end{remark}

\section*{Declarations}
\begin{itemize}
\item Conflict of interest: The author declares that he has no conflict of interest.

\item Data availability: No datasets were generated or analysed during the current study.

\end{itemize}


\begin{bibdiv}
\begin{biblist}

\bib{lw}{article}{
     author    = {Li, Haozhao},
     author    = {Wang, Bing},
     title     = {The Extension Problem of the Mean Curvature Flow (I)},
     journal   = {Invent. math.},
     volume    = {218},
     year      = {2019},
     pages     = {721\ndash 777},
   }
\bib{ac}{article}{
     author    = {A. Cooper, Andrew},
     title     = {A Characterization of the Singular Time of the Mean Curvature Flow},
     journal   = {Proceedings of the American Mathematical Society},
     volume    = {139},
     year      = {2011},
     pages     = {2933\ndash 2942},
   }

\bib{hui}{article}{
     author    = {Huisken, Gerhard},
     title     = {Flow by Mean Curvature of Convex Surfaces into Spheres},
     journal   = {J. Differential Geometry},
     volume    = {20},
     year      = {1984},
     pages     = {237\ndash 266},
   }
 \bib{lie}{book}{
    AUTHOR = {Lieberman, Gary M.},
     TITLE = {Second Order Parabolic Differential Equations},
 PUBLISHER = {World Scientific Publishing Co. Pte. Ltd.},
  YEAR = {1996},
 }
 \bib{dav}{article}{
     author    = {Davini, Andrea},
     title     = {On Calibrations for Lawson's Cones},
     journal   = {Rend. Sem. Mat. Univ. Padova},
     volume    = {111},
     year      = {2004},
     pages     = {55\ndash 70},
   }
 \bib{vel}{article}{
     author    = {Vel\'{a}zquez, J.J.L.},
     title     = {Curvature Blow-up in Perturbations of Minimal Cones Evolving by Mean Curvature Flow},
     journal   = {Annali della Scuola Normale Superiore di Pisa, Classe di Scienze},
     volume    = {21},
     year      = {1994},
     pages     = {595\ndash 628},
   }
 \bib{gs}{article}{
     author    = {Guo, Siao-Hao},
     author    = {Sesum, Natasa},
     title     = {Analysis of Vel\'{a}zquez’s Solution to the Mean Curvature Flow with a Type II Singularity},
     journal   = {Communications in Partial Differential Equations},
     volume    = {43},
     year      = {2018},
     pages     = {185\ndash 285},
   }
 \bib{sto}{article}{
     author    = {Stolarski, Maxwell},
     title     = {Existence of Mean Curvature Flow Singularities with Bounded Mean Curvature},
     journal   = {Duke Math. J.},
     volume    = {172(7)},
     year      = {2023},
     pages     = {1235\ndash 1292},
   }
 \bib{eh}{article}{
     author    = {Ecker, Klaus},
     author    = {Huisken, Gerhard},
     title     = {Interior Estimates for Hypersurfaces Moving by Mean Curvature},
     journal   = {Invent. Math.},
     volume    = {105},
     year      = {1991},
     pages     = {547\ndash 569},
   }
 \bib{hv}{article}{
     author    = {Herrero, M.A.},
     author    = {Vel\'{a}zquez, J.J.L.},
     title     = {A Blow-up Result for Semilinear Heat Equations in the Supercritical Case},
     journal   = {unpublished},
   }
 \bib{bd}{article}{
     author    = {Betancor, J.J.},
     author    = {De León-Contreras, M.},
     title     = {Parabolic Equations Involving Bessel Operators and Singular Integrals},
     journal   = {Integral Equations and Operator Theory},
     volume    = {90},
     year      = {2018},
     number    = {18},
   }
 \bib{wat}{book}{
    AUTHOR = {Watson, G.N.},
     TITLE = {A Treatise on the Theory of Bessel Functions},
 PUBLISHER = {Cambridge at the University Press},
  YEAR = {1922},
 }
  \bib{tem}{book}{
    AUTHOR = {Temme, N.M.},
     TITLE = {Asymptotic Methods for Integrals},
 PUBLISHER = {World Scientific Publishing Co. Pte. Ltd.},
  YEAR = {2015},
 }
  \bib{gp}{book}{
    AUTHOR = {Galindo, Alberto},
    AUTHOR = {Pascual, Pedro},
     TITLE = {Quantum Mechanics I},
 PUBLISHER = {Springer Berlin, Heidelberg},
  YEAR = {1990},
 }
  \bib{wan}{article}{
     author    = {Wang, Mu-Tao},
     title     = {The Mean Curvature Flow Smoothes Lipschitz
Submanifolds},
     journal   = {Communications in
Analysis and Geometry},
     volume    = {12(3)},
     year      = {2004},
     pages     = {581\ndash 599},
   }
   \bib{man}{book}{
    AUTHOR = {Mantegazza, Carlo},
     TITLE = {Lecture Notes on Mean Curvature Flow},
 PUBLISHER = {Springer Basel AG},
  YEAR = {2011},
 }

\end{biblist}
\end{bibdiv}
\end{document}